\documentclass[aop,MSNbibl,citesort,seceqn,dvips]{arximspdf}
\usepackage{mathbh}
\usepackage{mathrsfs}
\usepackage{graphicx}

\doi{10.1214/11-AOP657}
\volume{40}
\issue{4}
\pubyear{2012}
\firstpage{1577}
\lastpage{1635}

\makeatletter
\newcommand\R{{\mathbb R}}
\newcommand\C{{\mathbb C}}
\newcommand\Z{{\mathbb Z}}
\newcommand\N{{\mathbb N}}
\renewcommand\P{{\mathbb P}}
\newcommand\E{{\mathbb E}}
\renewcommand\H{\mathfrak{H}}
\newcommand\HH{\mathfrak{H}}
\newcommand\1{{\mathbh 1}}
\newcommand\tensor{\otimes}
\newcommand\supp{\operatorname{supp}}
\newcommand\del{\partial}

\newcommand{\ff}{\varphi}
\newcommand{\CC}{{\mathbb C}}
\newcommand\sh{\sharp}

\newcommand{\cont}[1]{\stackrel{#1}{\frown}}

\newtheorem{theorem}{Theorem}[section]

\newtheorem{proposition}[theorem]{Proposition}
\newproclaim{definition}[theorem]{Definition}
\newtheorem{corollary}[theorem]{Corollary}
\newtheorem{lemma}[theorem]{Lemma}
\newproclaim{remark}[theorem]{Remark}
\makeatother

\begin{document}
\begin{frontmatter}

\title{Wigner chaos and the fourth moment}
\runtitle{Wigner chaos and the fourth moment}

\begin{aug}
\author[A]{\fnms{Todd} \snm{Kemp}\corref{}\thanksref{a1}\ead[label=e1]{tkemp@math.ucsd.edu}},
\author[B]{\fnms{Ivan} \snm{Nourdin}\thanksref{a2}\ead[label=e2]{inourdin@gmail.com}},
\author[C]{\fnms{Giovanni} \snm{Peccati}\ead[label=e3]{giovanni.peccati@gmail.com}}\\
\and
\author[D]{\fnms{Roland} \snm{Speicher}\thanksref{a4}\ead[label=e4]{speicher@math.uni-sb.de}}

\runauthor{Kemp, Nourdin, Peccati and Speicher}
\affiliation{UCSD, Universit\'e Henri Poincar\'e, Universit\'e du
Luxembourg and Universit\"at~des Saarlandes}
\address[A]{T. Kemp\\
Department of Mathematics\\
University of California, San Diego\\
La Jolla, California 92093-0112\\
USA\\\printead{e1}} 
\address[B]{I. Nourdin\\
Universit\'{e} de Lorraine\\
Institut Elie Cartan de Math\'{e}matiques\\
Vandoeuvre-l\`{e}s-Nancy, F-54506\\
France\\
\printead{e2}}
\address[C]{G. Peccati\\Universit\'e du Luxembourg\\
Facult\'{e} des Sciences\\
6, rue Richard Coudenhove-Kalergi\\
L-1359 Luxembourg\\
\printead{e3}}
\address[D]{R. Speicher\\
Fachrichtung Mathematik\\
Universit\"{a}t des Saarlandes\hspace*{45pt}\\
Postfach 151150\\
66041 Saarbr\"{u}cken\\
Germany\\
\printead{e4}}
\end{aug}
\thankstext{a1}{Supported in part by NSF Grants DMS-07-01162 and
DMS-10-01894.}
\thankstext{a2}{Supported in part by (French) ANR grant
``Exploration des Chemins Rugueux.''}
\thankstext{a4}{Supported in part by a Discovery grant from NSERC.}

\received{\smonth{9} \syear{2010}}
\revised{\smonth{2} \syear{2011}}

%
\begin{abstract}
We prove that a normalized sequence of multiple Wigner integrals (in a
fixed order of free Wigner chaos) converges in law to the standard
semicircular distribution if and only if the corresponding sequence of
fourth moments converges to $2$, the fourth moment of the semicircular
law. This extends to the free probabilistic, setting some recent
results by Nualart and Peccati on characterizations of central limit
theorems in a fixed order of Gaussian Wiener chaos. Our proof is
combinatorial, analyzing the relevant noncrossing partitions that
control the moments of the integrals. We can also use these techniques
to distinguish the first order of chaos from all others in terms of
\mbox{distributions}; we then use tools from the free Malliavin calculus to
give quantitative bounds on a distance between different orders of
chaos. When applied to highly symmetric kernels, our results yield a
new transfer principle, connecting central limit theorems in free
Wigner chaos to those in Gaussian Wiener chaos. We use this to prove a
new free version of an important classical theorem, the Breuer--Major theorem.
\end{abstract}

%

\end{frontmatter}

\section{Introduction and background} \label{sectionintroduction}

Let $(W_t)_{t\ge0}$ be a standard one-dimen\-sional Brownian motion, and
fix an
integer $n\geq1$. For every deterministic (Lebesgue) square-integrable
function $f$ on
$\mathbb{R}_+^n$, we denote by $I_n^W(f)$ the $n$th (multiple)
Wiener--It\^o stochastic
integral of $f$ with respect to $W$ (see, e.g., \cite{Janson,Kuo,nualartbook,PecTaqbook} for definitions; here and in the
sequel $\R_+$ refers to the nonnegative half-line $[0,\infty)$).
Random variables such as $I^W_n(f)$ play a~fundamental role in modern stochastic
analysis, the key fact being that every square-integrable functional
of\vadjust{\goodbreak}
$W$ can be
uniquely written as an infinite orthogonal sum of symmetric Wiener--It\^o
integrals of
increasing orders. This feature, known as the \textit{Wiener--It\^o chaos
decomposition}, yields an explicit representation of the isomorphism
between the space of square-integrable functionals of $W$ and the
symmetric Fock space
associated with $L^2(\mathbb{R_+})$. In particular, the Wiener chaos is
the starting point of the powerful \textit{Malliavin calculus of
variations} and its many applications in theoretical and applied
probability (see again \cite{Janson,nualartbook} for an introduction
to these topics). We recall that the collection of all random variables
of the type $I^W_n(f)$, where~$n$ is a~fixed integer, is customarily
called the~$n$th \textit{Wiener chaos} associated with~$W$. Note that the
first Wiener chaos is just the Gaussian space spanned by~$W$.

The following result, proved in \cite{nunugio}, yields a very
surprising condition under which
a sequence $I^W_n(f_k)$ converges in distribution, as $k\rightarrow
\infty$, to a~Gaussian random variable.
[In this statement, we assume as given an underlying probability space
$(X, \mathcal{F}, \P)$,
with the symbol $\E$ denoting expectation with respect to
$\P$.]

\begin{theorem}[(Nualart, Peccati)] \label{thmNP1} Let $n\ge2$ be an
integer, and let $(f_k)_{k\in\N}$ be a sequence of symmetric functions
(cf. Definition~\ref{defmirror} below) in $L^2(\R_+^n)$, each with
$n!\|f_k\|_{L^2(\R_+^n)}=1$. The following statements are equivalent:
\begin{longlist}[(2)]
\item[(1)] The fourth moment of the stochastic integrals $I^W_n(f_k)$
converge to $3$.
\[
\lim_{k\to\infty} \E(I^W_n(f_k)^4) = 3.
\]
\item[(2)] The random variables $I^W_n(f_k)$ converge in distribution
to the standard normal law $N(0,1)$.
\end{longlist}
\end{theorem}
 Note that the Wiener chaos of order $n\ge2$ does not contain
any Gaussian random variables, cf. \cite{Janson}, Chapter~6. Since the
fourth moment of the normal $N(0,1)$ distribution is equal to $3$, this
Central Limit Theorem shows that, within a fixed order of chaos and as
far as normal approximations are concerned, second and fourth moments
alone control all higher moments of distributions.

\begin{remark} The Wiener isometry shows that the second moment
of~$I_n^W(f)$ is equal to $n!\|f\|_{L^2}^2$, and so Theorem~\ref{thmNP1}
could be stated intrinsically in terms of random variables in a fixed
order of Wiener chaos. Moreover, it could be stated with the a
priori weaker assumption that $\E(I_n^W(f_k)^2)\to\sigma^2$ for some
$\sigma>0$, with the results then involving $N(0,\sigma^2)$ and fourth
moment~$3\sigma^4$, respectively. We choose to rescale to variance $1$
throughout most of this paper.
\end{remark}

Theorem~\ref{thmNP1} represents a drastic simplification of the
so-called ``method of moments and cumulants'' for normal approximations
on a Gaussian space, as described, for example, in \cite{Major,Sur};
for a detailed in-depth treatement of these techniques in the arena of
Wiener chaos, see the forthcoming book~\cite{PecTaqbook}.\vadjust{\goodbreak} We refer the
reader to the survey~\cite{NP-survey} and the forthcoming
monograph~\cite{NP-book} for an introduction to several applications of Theorem
\ref{thmNP1} and its many ramifications, including power variations of
stochastic processes, limit theorems in homogeneous spaces, random
matrices and polymer fluctuations. See in particular~\cite{NP-PTRF,Noupecrei3,NO} for approaches to Theorem~\ref{thmNP1} based
respectively on Malliavin calculus and Stein's method, as well as
applications to universality results for nonlinear statistics of
independent random variables.

In the recent two decades, a new probability theory known as \textit{free
probability} has gained momentum due to its extremely powerful
contributions both to its birth subject of operator algebras and to
random matrix theory; see, for example, \cite{AGZ,HiaiPetz,NicaSpeicherBook,VDM}. Free probability theory offers a
new kind of independence between random variables, \textit{free
independence}, that is, modeled on the free product of groups rather
than tensor products; it turns out to succinctly describe the
relationship between eigenvalues of large random matrices with
independent entries. In free probability, the central limit
distribution is the Wigner semicircular law [cf. equation (\ref{eqsemicircle})], further demonstrating the link to random matrices.
\textit{Free Brownian motion}, discussed in Section~\ref{sectionfBM} below, is
a (noncommutative) stochastic process whose increments are freely
independent and have semicircular distributions. Essentially, one
should think of free Brownian motion as Hermitian random matrix-valued
Brownian motion in the limit as matrix dimension tends to infinity;
see, for example, \cite{BCG} for a~detailed analysis of the related
large deviations.

If $(S_t)_{t\ge0}$ is a free Brownian motion, the construction of the
Wiener--It\^o integral can be mimicked to construct the so-called
\textit{Wigner stochastic integral} (cf. Section~\ref{sectionWignerintegral}) $I^S_n(f)$ of a deterministic function $f\in L^2(\R_+^n)$.
The noncommutativity of $S_t$ gives $I_n^S$ different properties; in
particular, it is no longer sufficient to restrict to the class of
symmetric $f$. Nevertheless, there is an analogous theory of \textit{Wigner chaos} detailed in \cite{BianeSpeicher}, including many of the
powerful tools of Malliavin calculus in free form. The main theorem of
the present paper is the following precise analog of the central limit
Theorem~\ref{thmNP1} in the free context.\vspace*{-3pt}

\begin{theorem} \label{thm4thmomentsemicircle} Let $n\ge2$ be an
integer, and let $(f_k)_{k\in\N}$ be a sequence of mirror symmetric
functions (cf. Definition~\ref{defmirror}) in $L^2(\R_+^n)$, each
with $\|f_k\|_{L^2(\R_+^n)} = 1$. The following statements are equivalent:
\begin{longlist}[(2)]
\item[(1)] The fourth moments of the Wigner stochastic integrals
$I^S_n(f_k)$ converge to~$2$.
\[
\lim_{k\to\infty} \E(I^S_n(f_k)^4) = 2.
\]
\item[(2)] The random variables $I^S_n(f_k)$ converge in law to the
standard semicircular distribution $S(0,1)$ [cf. equation (\ref{eqsemicircle})] as $k\to\infty$.\vspace*{-3pt}
\end{longlist}
\end{theorem}

\begin{remark} The expectation $\E$ in Theorem~\ref{thm4thmomentsemicircle}(1) must be properly interpreted in the free context; in
Section~\ref{sectionNCprob}\vadjust{\goodbreak} we will discuss the right framework (of a
trace $\E=\varphi$ on the von Neumann algebra generated by the free
Brownian motion). We will also make it clear what is meant by the law
of a~noncommutative random variable like $I^S_n(f_k)$.
\end{remark}

\begin{remark} \label{rklawimpliesmoments} Since the fourth moment
of the standard semicircular distribution is $2$, (2) nominally implies
(1) in Theorem~\ref{thm4thmomentsemicircle} since convergence in
distribution implies convergence of moments (modulo growth
constraints); the main thrust of this paper is the remarkable reverse
implication. The mirror symmetry condition on $f$ is there merely to
guarantee that the stochastic integral $I_n^S(f)$ is indeed a
self-adjoint operator; otherwise, it has no law to speak of (cf.
Section~\ref{sectionNCprob}).
\end{remark}

Our proof of Theorem~\ref{thm4thmomentsemicircle} is through the
method of moments which, in the context of the Wigner chaos, is
elegantly formulated in terms of \textit{noncrossing} pairings and
partitions. While, on some level, the combinatorics of partitions can
be seen to be involved in any central limit theorem, our present proof
is markedly different from the form of the proofs given in \cite{Noupecrei3,NO,nunugio}. All relevant technology is discussed in
Sections~\ref{sectionNCprob}--\ref{sectionNCpartitions} below;
further details on the method of moments in free probability theory can
be found in the book~\cite{NicaSpeicherBook}.

As a key step toward proving Theorem~\ref{thm4thmomentsemicircle},
but of independent interest and also completely analogous to the
classical case, we prove the following characterization of the fourth
moment condition in terms of standard integral contraction operators on
the kernels of the stochastic integrals (as discussed at length in
Section~\ref{sectionWignerintegral} below).

\begin{theorem} \label{thm4thmomentcont}
Let $n$ be a natural number, and let $(f_k)_{k\in\N}$ be a sequence
of functions in $L^2(\R_+^n)$, each with $\|f_k\|_{L^2(\R_+^n)} = 1$.
The following statements are equivalent:
\begin{longlist}[(2)]
\item[(1)] The fourth absolute moments of the stochastic integrals
$I^S_n(f_k)$ converge to~$2$.
\[
\lim_{k\to\infty} \E(|I^S_n(f_k)|^4) = 2.
\]
\item[(2)] All nontrivial contractions (cf. Definition~\ref{defcontraction}) of $f_k$ converge to $0$: for each $p=1,2,\ldots,n-1$,
\[
\lim_{k\to\infty} f_k\cont{p} f_k^\ast= 0 \qquad\mbox{in }
L^2(\R_+^{2n-2p}).
\]
\end{longlist}
\end{theorem}

While different orders of Wiener chaos have disjoint classes of laws,
it is (at the present time) unknown if the same holds for the Wigner
chaos. As a~first result in this direction, the following important
corollary to Theorem~\ref{thm4thmomentcont} allows us to distinguish
the laws of Wigner integrals in the first order of chaos from all
higher orders.

\begin{corollary} \label{cornosemicircular} Let $n\ge2$ be an
integer, and consider a nonzero mirror symmetric function $f\in L^2(\R
_+^n)$. Then the Wigner integral $I^S_n(f)$ satisfies $\E[I_n^S(f)^4] >
2\E[I^S_n(f)^2]^2$. In particular, the distribution of the Wigner
integral $I^S_n(f)$ cannot be semicircular.
\end{corollary}

Combining these results with those in \cite{NP-PTRF,Noupecrei3,NO,nunugio}, we can state the following
\textit{Wiener--Wigner transfer principle} for translating results between the
classical and free chaoses.

\begin{theorem} \label{thmtransfer}
Let $n\ge2$ be an integer, and let $(f_k)_{k\in\N}$ be a sequence of
fully symmetric (cf. Definition~\ref{defmirror}) functions in
$L^2(\R
_+^n)$. Let $\sigma>0$ be a~finite constant. Then, as $k\to\infty$:
\begin{longlist}[(2)]
\item[(1)] $\E [I^W_n(f_k)^2 ] \to n!\sigma^2$ if and
only if
$\E [I^S_n(f_k)^2 ]\to\sigma^2$.
\item[(2)] If the asymptotic relations in (1) are verified, then
$I^W_n(f_k)$ converges in law to a normal random variable $N(0,n!\sigma
^2)$ if and only if $I^S_n(f_k)$ converges in law to a semicircular
random variable $S(0,\sigma^2)$.
\end{longlist}
\end{theorem}

Theorem~\ref{thmtransfer} will be shown by combining Theorems~\ref{thm4thmomentsemicircle} and~\ref{thm4thmomentcont} with the findings
of \cite{nunugio}; the transfer principle allows us to easily prove yet
unknown free versions of important classical results, such as the
Breuer--Major theorem (Corollary~\ref{corBreuerMajor} below).

\begin{remark} It is important to note that the transfer principle
Theorem~\ref{thmtransfer} requires the strong assumption that the
kernels $f_k$ are \textit{fully symmetric} in both the classical and free
cases. While this is no loss of generality in the Wiener chaos, it
applies to only a small subspace of the Wigner chaos of orders $3$ or higher.
\end{remark}

Corollary~\ref{cornosemicircular} shows that the semicircular law is
not the law of any stochastic integral of order higher than $1$. We are
also able to prove some sharp quantitative estimates for the distance
to the semicircular law. The key estimate, using Malliavin calculus, is
as follows: it is a free probabilistic analog of \cite{NP-PTRF}, Theorem~3.1. We state it here in less generality than we prove it in
Section~\ref{sectionMalliavinestimates}.

\begin{theorem} \label{theoremMalliavinestimate} Let $S$ be a
standard semicircular random variable [cf. equation (\ref{eqsemicircle})]. Let $F$ have a finite Wigner chaos expansion;
that is,
$F = \sum_{n=1}^N I^S_n(f_n)$ for some mirror symmetric functions
$f_n\in L^2(\R_+^n)$ and some finite $N$. Let $\mathcal{C}_2$ and
$\mathscr{I}_2$ be as in Definition~\ref{defFourierclass}. Then
%
\begin{eqnarray} \label{eqWasserstein1intro} d_{\mathcal
{C}_2}(F,S)&\equiv&\mathop{\mathop{\sup}_{h\in\mathcal{C}_2}}_{\mathscr{I}_2(h)\le
1} |\E
[h(F)] - \E[h(S)]|\nonumber
\\[-8pt]
\\[-8pt] &\le&\frac{1}{2}\E\tensor\E \biggl( \biggl|\int
_0^\infty
\nabla_t(N_0^{-1}F)\sh(\nabla_t F)^\ast\,dt - 1\tensor1
\biggr|
\biggr).
\nonumber
\end{eqnarray}
\end{theorem}

 The Malliavin calculus operators $\nabla$ and $N_0$ and the
product $\sh$ on tensor-product-valued \textit{biprocesses} are defined
below in Section~\ref{sectionMalliavin}, where we also describe all
the relevant structure, including why the free
Cameron--Gross--Malliavin derivative $\nabla_t F$ of a random variable
$F$ takes values in the tensor product $L^2(\R_+)\tensor L^2(\R_+)$.
The class $\mathcal{C}_2$ is somewhat smaller than the space of
Lipschitz functions, and so the metric $d_{\mathcal{C}_2}$ on the
left-hand side of equation~(\ref{eqWasserstein1}) is, a priori, weaker
than the Wasserstein metric. This distance does metrize convergence in
law, however.\vspace*{-3pt}

\begin{remark} \label{remarksmartpath} The key element in the proof
of Theorem~\ref{theoremMalliavinestimate} is to measure the distance between $F$ and $S$ by means
of a procedure close to the so-called \textit{smart path method}, as
popular in Spin Glasses; cf. \cite{Talagrand}. In this technique, one
assumes that $F$ and $S$ are independent, and then assesses the
distance between their laws by controlling the variations of the
mapping $t\mapsto\mathbb{E}[h(\sqrt{1-t}F+\sqrt{t}S)]$ (where $h$
is a
suitable test function) over the interval $[0,1]$. As shown below, our
approach to the smart path method requires that we replace $\sqrt{t}S$
by a free Brownian motion $S_t$ (cf. Section~\ref{sectionfBM})
\textit{freely independent} from $F$, so that we can use the free stochastic
calculus to proceed with our estimates.\vspace*{-3pt}
\end{remark}

Using Theorem~\ref{theoremMalliavinestimate}, we can prove the
following sharp quantitative bound for the distance from any double
Wigner integral to the semicircular law.\vspace*{-2pt}

\begin{corollary} \label{corWasserstein} Let $f\in L^2(\R_+^2)$ be
mirror-symmetric and normalized $\|f\|_{L^2(\R_+^n)}=1$, let $S$ be a
standard semicircular random variable and let $d_{\mathcal{C}_2}$ be
defined as in equation (\ref{eqWasserstein1intro}). Then
%
\begin{equation} \label{eqWasserstein2} d_{\mathcal{C}_2}(I^S_2(f),S)
\le\frac{1}{2}\sqrt{\frac{3}{2}} \sqrt{\E[I^S_2(f)^4]-2}.\vspace*{-2pt}
\end{equation}
\end{corollary}

 In principle, equation (\ref{eqWasserstein1intro}) could be
used to give quantitative estimates like equation (\ref{eqWasserstein2}) for any order of Wigner chaos. However, the analogous techniques
from the classical literature heavily rely on the full symmetry of the
function $f$; in the more general mirror symmetric case required in the
Wigner chaos, such estimates are, thus far, beyond our reach.

The remainder of this paper is organized as follows. Sections~\ref{sectionNCprob} through~\ref{sectionNCpartitions} give (concise)
background and notation for the free probabilistic setting, free
Brownian motion and its associated stochastic integral the Wigner
integral and the relevant class of partitions (noncrossing pairings)
that control moments of these integrals. Section~\ref{sectionCLTs} is
devoted to the proofs of Theorems~\ref{thm4thmomentsemicircle} and
\ref{thm4thmomentcont} along with Corollary~\ref{cornosemicircular} and Theorem~\ref{thmtransfer}. In Section~\ref{sectionMalliavin}, we collect and summarize all of the tools of free
stochastic calculus and free Malliavin calculus needed to prove the
quantitative results of Section~\ref{lastsection}; this final
section
is devoted to the proofs of Theorem~\ref{theoremMalliavinestimate}
(in Section~\ref{sectionMalliavinestimates}) and Corollary~\ref{corWasserstein}
(in Section~\ref{sectiondistanceestimate}), along with
an abstract list of equivalent forms of our central limit theorem in
the second Wigner chaos. Finally,\vadjust{\goodbreak} \hyperref[appendix]{Appendix} contains the
proof of Theorem~\ref{theorempolynomialsdense}, an important
technical approximation tool needed for the proof of Theorem~\ref{theoremMalliavinestimate}
but also of independent interest.\vspace*{-3pt}

\subsection{Free probability} \label{sectionNCprob}
A \textit{noncommutative probability space} is a complex linear algebra $\mathscr
{A}$ equipped with an involution (like the adjoint operation $X\mapsto
X^\ast$ on matrices) and a unital linear functional $\ff\dvtx
\mathscr
{A}\to\C$. The standard classical example is $\mathscr{A} = L^\infty
(\Omega,\mathcal{F},\P)$ where $\mathcal{F}$ is a $\sigma$-field of
subset of $\Omega$, and $\P$ is a probability measure on $\mathcal{F}$;
in this case the involution is complex conjugation and $\varphi$ is
expectation with respect to $\P$. One can identify $\mathcal{F}$ from
$\mathscr{A}$ through the idempotent elements which are the indicator
functions $\1_E$ of events $E\in\mathcal{F}$, and so this terminology
for a probability space contains the same information as the usual one.
Another relevant example that is actually noncommutative is given by
\textit{random matrices}; here $\mathscr{A} = L^\infty(\Omega,\mathcal
{F},\P ; M_d(\C))$, $d\times d$-matrix-valued random variables, where
the involution is matrix adjoint and the natural linear functional $\ff
$ is given by $\ff(X) = \frac{1}{d}\E\operatorname{Tr}(X)$. Both of these
examples only deal with bounded random variables, although this can be
extended to random variables with finite moments without too much effort.

The pair $(L^\infty(\Omega,\mathcal{F},\P),\E)$ has a lot of analytic
structure not present in many noncommutative probability spaces; we
will need these analytic tools in much of the following work. We assume
that $\mathscr{A}$ is a von Neumann algebra, an algebra of operators on
a (separable) Hilbert space, closed under adjoint and weak convergence.
Moreover, we assume that the linear functional $\ff$ is weakly
continuous, positive [meaning $\ff(X)\ge0$ whenever $X$ is a
nonnegative element of $\mathscr{A}$; i.e., whenever $X=YY^\ast$ for
some $Y\in\mathscr{A}$], faithful [meaning that if $\ff(YY^\ast)=0$,
then $Y=0$] and tracial, meaning that $\ff(XY)=\ff(YX)$ for all
$X,Y\in
\mathscr{A}$, even though in general $XY\ne YX$. Such a $\ff$ is called
a \textit{trace} or \textit{tracial state}. Both of the above examples
(bounded random variables and bounded random matrices) satisfy these
conditions. A von Neumann algebra equipped with a tracial state is
typically called a \textit{(tracial) $W^\ast$-probability space}. Some of
the theorems in this paper require the extra structure of a $W^\ast
$-probability space, while others hold in a general abstract
noncommutative probability space. To be safe, we generally assume the
$W^\ast$-setting in what follows. Though we do not explicitly specify
traciality in the proceeding, we will always assume $\ff$ is a trace.

In a $W^\ast$-probability space, we refer to the self-adjoint elements
of the algebra as \textit{random variables}. Any random variable has a
\textit{law} or \textit{distribution} defined as follows: the law of $X\in
\mathscr{A}$ is the unique Borel probability measure~$\mu_X$ on~$\R$
with the same moments as $X$; that is, such that
\[
\int_{\R} t^n \mu_X(dt) = \ff(X^n), \qquad n=0,1,\ldots.
\]
The existence and uniqueness of $\mu_X$ follow from the positivity of
$\ff$; see~\cite{NicaSpeicherBook}, Propositions~3.13. Thus, in
general, noncommutative probability, the method of moments and
cumulants plays a central role.\vadjust{\goodbreak}

In this general setting, the notion of \textit{independence} of events is
harder to pin down. Voiculescu introduced a general noncommutative
notion of independence in \cite{Voiculescu} which has, of late, been
very important both in operator algebras and in random matrix theory.
Let $\mathscr{A}_1,\ldots,\mathscr{A}_n$ be unital subalgebras of
$\mathscr{A}$. Let $X_1,\ldots, X_m$ be elements chosen from among the
$\mathscr{A}_i$'s such that, for $1\le j<m$, $X_j$ and $X_{j+1}$ do not
come from the same $\mathscr{A}_i$, and such that $\ff(X_j)=0$ for each
$j$. The subalgebras $\mathscr{A}_1,\ldots,\mathscr{A}_n$ are said to
be \textit{free} or \textit{freely independent} if, in this circumstance,
$\ff
(X_1X_2\cdots X_n) = 0$. Random variables are called freely independent
if the unital algebras they generate are freely independent. By
centering moments it is easy to check that, in the case that all the
indices are distinct, this is the same as classical independence
expressed in terms of moments. For example, if $X,Y$ are freely
independent they satisfy $\ff [(X^n-\ff(X^n))(Y^m-\ff
(Y^m)) ]
= 0$, which reduces to $\ff(X^nY^m) = \ff(X^n)\ff(Y^m)$. But if there
are repetitions of indices among the (generally noncommutative) random
variables, freeness is much more complicated than classical
independence; for example, if $X,Y$ are free, then $\ff(\mathit{XYXY}) = \ff
(X^2)\ff(Y)^2 + \ff(X)^2\ff(Y^2) - \ff(X)^2\ff(Y)^2$.
Nevertheless, if
$X,Y$ are freely independent, then their joint moments are determined
by the moments of $X$ and $Y$ separately. Indeed, the law of the random
variable $X+Y$ is determined by (and can be calculated using the
Stieltjes transforms of) the laws of $X$ and $Y$ separately. It was
later discovered by Voiculescu~\cite{Voiculescu2} and others that pairs
of random matrices with independent entries are (asymptotically) freely
independent in terms of expected trace; this has led to powerful new
tools for analyzing the density of eigenvalues of random matrices.

The notion of \textit{conditioning} is also available in free probability.
\begin{definition}\label{defconditionalexpectation}
Let $(\mathscr{A},\ff)$ be a $W^\ast$-probability space, and let
$\mathscr{B}\subseteq\mathscr{A}$ be a unital $W^\ast$-subalgebra.
There is a \textit{conditional expectation} map $\ff[ \cdot |\mathscr
{B}]$ from $\mathscr{A}$ onto $\mathscr{B}$. It is characterized by
the property
%
\begin{equation} \label{eqconditionalexpectation}
\ff[XY] = \ff [X\ff[Y|\mathscr{B}] ] \qquad\mbox{for all }X\in\mathscr{B}, Y\in\mathscr{A}.
\end{equation}
Conditional expectation has the following properties:
\begin{longlist}[(3)]
\item[(1)] $\ff[ \cdot |\mathscr{B}]$ is weakly continuous and
completely positive;
\item[(2)] $\ff[ \cdot |\mathscr{B}]$ is a contraction (in operator
norm) and preserves the identity;
\item[(3)] If $Y\in\mathscr{A}$ and $X,Z\in\mathscr{B}$, then $\ff
[\mathit{XYZ}|\mathscr{B}] = X\ff[Y|\mathscr{B}]Z$.
\end{longlist}
If $X\in\mathscr{A}$, then we denote by $\ff[ \cdot |X]$ the
conditional expectation onto the unital von Neumann subalgebra of
$\mathscr{A}$ generated by $X$.
\end{definition}

 Such conditional expectations were introduced in \cite{Takesaki}
 [where properties (1)--(3) were proved]. As one should
expect, if $X$ and $Y$ are free, then~$\ff[Y|X] =\allowbreak \ff(Y)$, as in the
classical case. Many analogs of classical probabilistic
constructions
(such as martingales) are well-defined in free probability, using
Definition~\ref{defconditionalexpectation}. See, for example, \cite{Biane2} for a discussion of
\textit{free L\'evy processes}.\vadjust{\goodbreak}

\subsection{Free Brownian motion} \label{sectionfBM} The (centred)
\textit{semicirclular distribution} (or Wigner law) $S(0,t)$ is the
probability distribution
%
\begin{equation} \label{eqsemicircle} S(0,t)(dx) = \frac{1}{2\pi t}
\sqrt{4t-x^2}\,dx, \qquad|x|\le2\sqrt{t}.
\end{equation}
Since this distribution is symmetric about $0$, its odd moments are all
$0$. Simple calculation shows that the even moments are given by
(scaled) \textit{Catalan numbers}: for nonnegative integers $m$,
\[
\int_{-2\sqrt{t}}^{2\sqrt{t}} x^{2m} S(0,t)(dx) = C_m t^m,
\]
where $C_m = \frac{1}{m+1}{2m\choose m}$. In particular, the second
moment (and variance) is $t$ while the fourth moment is $2t^2$.

A \textit{free Brownian motion} $S=(S_t)_{t\ge0}$ is a noncommutative
stochastic process; it is a one-parameter family of self-adjoint
operators $S_t$ in a $W^\ast$-probability space $(\mathscr{A},\ff)$,
with the following defining characteristics:
\begin{longlist}[(2)]
\item[(0)] $S_0 = 0$;
\item[(1)] For $0<t_1<t_2<\infty$, the law of $S_{t_2}-S_{t_1}$ is the
semicircular distribution of variance $t_2-t_1$;
\item[(2)] For all $n$ and $0<t_1<t_2<\cdots<t_n<\infty$, the
increments $S_{t_1},S_{t_2}-S_{t_1}$, $S_{t_3}-S_{t_2}$, \ldots,
$S_{t_n}-S_{t_{n-1}}$ are freely independent.
\end{longlist}
The freeness of increments can also be expressed by saying that
$S_{t_2}-S_{t_1}$ is free from $\mathcal{S}_{t_1}$ whenever
$t_2>t_1\ge
0$; here $\mathcal{S}_t$ is the von Neumann algebra generated by $\{
S_s \dvtx 0\le s\le t\}$. In particular, it follows easily that $\ff
[S_{t_2}|\mathcal{S}_{t_1}] = S_{t_1}$ for $t_2\ge t_1\ge0$, so free
Brownian motion is a martingale.

There are at least two good ways to construct a~free Brownian motion~$S$.
The first involves the \textit{free (Boltzman) Fock space} $\mathscr
{F}_0(\mathfrak{H})$ constructed on a~Hilbert space $\mathfrak{H}$:
$\mathscr{F}_0(\mathfrak{H}) \equiv\bigoplus_{n=0}^\infty\mathfrak
{H}^{\tensor n}$ where the direct-sum and tensor products are Hilbert
space operations, and $\mathfrak{H}^{\tensor0}$ is defined to be a
one-dimensional complex space with a distinguished unit basis vector
called the \textit{vacuum} $\Omega$ (not to be confused with the state
space of a probability space). Given any vector $h\in\mathfrak{H}$,
the \textit{creation operator} $a^{\dagger}(h)$ on $\mathscr
{F}_0(\mathfrak
{H})$ is defined by left tensor-product with $h$: $a^\dagger(h) \psi=
h\tensor\psi$. Its adjoint $a(h)$ is the \textit{annihilation operator},
whose action on an $n$-tensor is given by $a(h) h_1\tensor\cdots
\tensor h_n = \langle h,h_1\rangle h_2\tensor\cdots\tensor h_n$
[and
$a(h)\Omega= 0$]. The creation and annihilation operators are thus
raising and lowering operators. Their sum $X(h) = a^\dagger(h) + a(h)$
is a~self-adjoint operator known as the \textit{field operator} in the
direction $h$. Let~$\mathcal{S}(\mathfrak{H})$ denote the von Neumann
algebra generated by $\{X(h) ; h\in\mathfrak{H}\}$, a~(small)
subset
of all bounded operators on the Fock space $\mathscr{F}_0(\mathfrak
{H})$. The \textit{vacuum expectation state} $\ff(Y) = \langle Y\Omega
,\Omega\rangle_{\mathscr{F}_0(\mathfrak{H})}$ is a tracial state on
$\mathcal{S}(\mathfrak{H})$. Now, take the special case $\mathfrak
{H}=L^2(\R_+)$; then $S_t = X(\1_{[0,t]})$ is a free Brownian motion
with respect to $(\mathcal{S}(\mathfrak{H}),\ff)$.

\begin{remark} This construction of Brownian motion can also be done in
the classical case, replacing the free Fock space with the symmetric
(Bosonic) Fock space; for this line of thought see \cite{Parthasarathy}. Although it is abstract, it is directly related to
concrete constructions in the Wigner, and Wiener, chaos. Note: when
$\mathfrak{H}=L^2(\R_+)$, $\mathfrak{H}^{\tensor n}$ may be identified
with $L^2(\R_+^n)$, and it is these kernels we will work with
throughout most of this paper.\vspace*{-2pt}
\end{remark}

A second, more appealing (if less direct) construction of free Brownian
motion uses random matrices. Let $W_t^d$ be a $d\times d$ complex
Hermitian matrix all of whose entries above the main diagonal are
independent complex standard Brownian motions. Set $S^d_t =
d^{-1/2}W_t^d$. Then the ``limit as $d\to\infty$'' of~$S_t^d$ is a free
Brownian motion. This limit property holds in the sense of moments, as
follows: equip the algebra $\mathcal{S}^d$ generated by $\{S_t^d ;
t\in\R_+\}$ with the tracial state $\ff_d = \frac{1}{d}\E \operatorname
{Tr}$. Then if $P = P(X_1,X_2,\ldots,X_k)$ is any polynomial in $k$
noncommuting indeterminates, and $t_1,\ldots,t_k\in\R_+$, then
\[
\lim_{d\to\infty}\ff_d [P(S^d_{t_1},\ldots,S^d_{t_k})
] = \ff
 [P(S_{t_1},\ldots,S_{t_k}) ],
\]
where $S = (S_t)_{t\ge0}$ is a free Brownian motion. So, at least in
terms of moments, we may think of free Brownian motion as
``infinite-dimensional matrix-valued Brownian motion.''\vspace*{-2pt}

\begin{remark} The algebra $\mathcal{S}^d$ of random matrices described
above is not a von Neumann algebra in the standard sense, since its
elements do not have finite matrix norms in the standard sup metric.
The Gaussian tails of the entries guarantee, however, that mixed matrix
moments of all orders are finite, which is all that is needed to make
sense of the standard notion of convergence in noncommutative
probability theory.\vspace*{-2pt}
\end{remark}

\subsection{The Wigner integral} \label{sectionWignerintegral} In
this section we largely follow \cite{BianeSpeicher}; related
discussions and extensions can be found in \cite{Ansh1,Ansh2,Ansh3}.
Taking a note from Wiener and It\^o, we define a stochastic integral
associated with free Brownian motion in the usual manner. Let $S$ be a
free Brownian motion, and let $f\in L^2(\R_+^n)$ be an off-diagonal
rectangular indicator function, taking the form $f = \1_{[s_1,t_1]
\times\cdots\times[s_n,t_n]}$, where the intervals $[s_1,t_1]$,
\ldots, $[s_n,t_n]$ are pairwise disjoint. The \textit{Wigner integral}
$I^S_n(f)$ is defined to be the product operator $I_n^S(f) =
(S_{t_1}-S_{s_1})\cdots(S_{t_n}-S_{s_n})$. Extend $I^S_n$ linearly
over the set of all off-diagonal step-functions, which is dense in
$L^2(\R_+^n)$. The freeness of the increments of $S$ yield the simple
Wigner isometry
%
\begin{equation} \label{eqWignerisometry} \ff [I^S_n(g)^\ast
I^S_n(f) ] = \langle f,g \rangle_{L^2(\R_+^n)}.
\end{equation}
In other words, $I^S_n$ is an isometry from the space of off-diagonal
step functions into the Hilbert space of operators generated by the
free Brownian motion $S$, equipped with the inner product $\langle X, Y
\rangle_\varphi= \ff [Y^\ast X ]$. This means~$I_n^S$ extends
to an isometry from the closure, which is the full space $L^2(\R_+^n)$,
thus fully defining the Wigner\vadjust{\goodbreak} integral. If $f$ is any function in
$L^2(\R_+^n)$, we may write
\[
I_n^S(f) = \int f(t_1,\ldots,t_n)\, dS_{t_1}\cdots \,dS_{t_n}.
\]
This stands in contrast to the classical Gaussian Wiener integral,
which we shall denote $I_n^W$:
\[
I_n^W(f) = \int f(t_1,\ldots,t_n)\, dW_{t_1}\cdots \,dW_{t_n}.\vspace*{-2pt}
\]

\begin{remark} This construction long post-dates Wigner's work. The
terminology was invented in \cite{BianeSpeicher} as a humorous nod to
the fact that Wigner's semicircular law plays the Central Limit role
here, and the similarity between the names \textit{Wigner} and \textit{Wiener}.\vspace*{-2pt}
\end{remark}

\begin{remark} \label{rknonsymmetry} This is the same as It\^o's
construction of the multiple Wiener integral in classical Wiener--It\^o
chaos. Note, however, that the increments $S_{t_1}-S_{s_1}, \ldots,
S_{t_n}-S_{s_n}$ \textit{do not commute}. Hence, unlike for the Wigner
integral, permuting the variables of $f$ generally changes the value
of~$I_n^S(f)$.\looseness=-1\vspace*{-2pt}
\end{remark}

 The image of the $n$-fold Wigner integral $I^S_n$ on all of
$L^2(\R_+^n)$ is called the $n$th order of \textit{Wigner chaos} or \textit{free chaos}. It is easy to calculate that different orders of chaos are
orthogonal from one another (in terms of the trace inner product); this
also follows from contraction and product formulas below. The
noncommutative $L^2$-space generated by $(S_t)_{t\ge0}$ is the
orthogonal sum of the orders of Wigner chaos; this is the free analog
of the Wiener chaos decomposition.\vspace*{-2pt}

\begin{remark} \label{rksemicircularfamily} The first Wigner chaos,
the image of $I_1^S$, is a centred \textit{semicircular family} in the
sense of \cite{NicaSpeicherBook}, Definition~8.15, exactly as the first
Wigner chaos is a centred Gaussian family. In particular, In the first
order of Wigner chaos, the law of any random variable is semicircular
$S(0,t)$ for some variance $t>0$.\vspace*{-2pt}
\end{remark}

We are generally interested only in self-adjoint elements of a given
order of chaos. Taking note of Remark~\ref{rknonsymmetry}, we have
%
\begin{eqnarray}\label{eqintadjoint}
 I_n^S(f)^\ast
&=& \biggl(\int f(t_1,\ldots,t_n) \,
dS_{t_1}\cdots \,dS_{t_n} \biggr)^\ast\nonumber
\\[-2pt]
&=& \int\overline{f(t_1,\ldots, t_n)} \, dS_{t_n}\cdots \,dS_{t_1} =
\int
\overline{f(t_n,\ldots,t_1)}\, dS_{t_1}\cdots \,dS_{t_n}\\[-2pt]
 &=& I^S_n(f^\ast),
\nonumber
\end{eqnarray}
where $f^\ast(t_1,\ldots,t_n) = \overline{f(t_n,\ldots,t_1)}$. This
prompts a definition.\vspace*{-2pt}

\begin{definition}\label{defmirror} Let $n$ be a natural number,
and let $f$ be a function in $L^2(\R_+^n)$.\vadjust{\goodbreak}
\begin{longlist}[(3)]
\item[(1)] The \textit{adjoint} of $f$ is the function $f^\ast
(t_1,\ldots
,t_n) = \overline{f(t_n,\ldots,t_1)}$.
\item[(2)] $f$ is called \textit{mirror symmetric} if $f=f^\ast$; that
is, if $f(t_1,\ldots,t_n) = \overline{f(t_n,\ldots,t_1)}$ for almost
all $t_1,\ldots,t_n\ge0$ with respect to the product Lebesgue measure.
\item[(3)] $f$ is called \textit{fully symmetric} if it is real-valued
and, for any permutation~$\sigma$ in the symmetric group $\Sigma_n$,
$f(t_1,\ldots,t_n) = f(t_{\sigma(1)},\ldots,t_{\sigma(n)})$ for almost
all $t_1,\ldots,t_n\ge0$ with respect to the product Lebesge measure.\vspace*{-3pt}
\end{longlist}
\end{definition}

 Thus an element $I^S_n(f)$ of the $n$th Wigner chaos is self
adjoint iff $f$ is mirror symmetric. Note, in the classical
Gaussian Wiener chaos, it is typical to consider only kernels that are
fully symmetric, since if $\tilde{f}$ is constructed from $f$ by
permuting its arguments, then $I^W_n(f) = I^W_n(\tilde{f})$. This
relation does \textit{not} hold for $I^S_n$.\vspace*{-3pt}

\begin{remark} The calculation in equation (\ref{eqintadjoint}) may
seem nonrigorous. A~more pedantic writing would do the calculation
first for an off-diagonal rectangular indicator function $f= \1
_{[s_1,t_1] \times\cdots\times[s_n,t_n]}$, in which case the adjoint
is merely $[(S_{t_1}-S_{s_1})\cdots(S_{t_n}-S_{s_n})]^\ast=
(S_{t_n}-S_{s_n})\cdots(S_{t_1}-S_{s_1})$ since~$S_t$ is self adjoint;
extending (sesqui)linearly and completing yields the full result. This
is how statements like $(dS_{t_1}\cdots \,dS_{t_n})^\ast= dS_{t_n}\cdots
\,dS_{t_1}$ should be interpreted throughout this paper.\vspace*{-3pt}
\end{remark}

Contractions are an important construction in Wigner and Wiener chaos;
we briefly review them now.\vspace*{-3pt}

\begin{definition}\label{defcontraction}
Let $n,m$ be natural numbers, and let $f\in L^2(\R^n_+)$ and $g\in
L^2(\R^m_+)$.
Let $p\le\min\{n,m\}$ be a natural number.
The \textit{$p$th contraction} $f\cont{p} g$ of $f$ and $g$ is the
$L^2(\R
^{n+m-2p}_+)$ function defined by nested
integration of the middle $p$ variables in $f\tensor g$
\begin{eqnarray*}
f\cont{p} g (t_1,\ldots,t_{n+m-2p}) &=& \int_{\R_+^{p}} f(t_1,\ldots
,t_{n-p},s_1,\ldots,s_p)\\[-2pt]
&&\hphantom{\int_{\R_+^{p}}}{}\times g(s_p,\ldots,s_1,t_{n-p+1},\ldots
,t_{n+m-2p})\,
ds_1\cdots \,ds_p.\vspace*{-3pt}
\end{eqnarray*}
\end{definition}

 Notice that when $p=0$, there is no integration, just the
products of $f$ and $g$ with disjoint arguments; in other words,
$f\cont
{0}g = f\tensor g$.\vspace*{-3pt}

\begin{remark} \label{rkcontractionsnotassociative} It is easy to
check that the operation $\cont{p}$ is not generally associative.\vspace*{-3pt}
\end{remark}

\begin{remark} \label{rkcontractionsvsnested} In \cite{NP-PTRF,Noupecrei3,NO,nunugio}
as well as standard references like~\cite{NP-survey,NP-book,nualartbook},
contractions are usually defined
as follows:
\begin{eqnarray*}
\hspace*{-0.5pt}f\tensor_p f(t_1,\ldots,t_{n+m-2p}) &=& \int_{\R_+^p}(t_1,\ldots
,t_{n-p},s_1,\ldots,s_p)\\[-2pt]
&&\hspace*{-0.5pt}\hphantom{\int_{\R_+^p}}{}\times g(t_{n-p+1},\ldots,t_{n+m-2p},s_1,\ldots
,s_p)\,
ds_1\cdots \,ds_p.\vadjust{\goodbreak}
\end{eqnarray*}
Notice that this operation is related to our nested contraction $\cont
{p}$ as follows:
\[
f\tensor_p \overline{g^\ast} (t_1,\ldots,t_{n-p},t_{n+m-2p},\ldots
,t_{n-p+1}) = f\cont{p}g(t_1,\ldots,t_{n+m-2p}).
\]
In other words, up to reordering of variables, the two operations are
the same. In particular, if $f,g$ are fully symmetric, then $f\cont{p}
g$ and $f\tensor_p g$ have the same symmetrizations. This will be
relevant to Theorem~\ref{thmtransfer} below.
\end{remark}

The following lemma records two useful facts about contractions and
adjoints; the proof is easy calculation.
\begin{lemma} \label{lemmacontadj} Let $n,m$ be natural numbers, and
let $f\in L^2(\R_+^n)$ and $g\in L^2(\R_+^m)$.
\begin{longlist}[(2)]
\item[(1)] If $p\le\min\{n,m\}$ is a natural number, then $(f\cont{p}
g)^\ast= g^\ast\cont{p}f^\ast$.
\item[(2)] If $n=m$, then the constant $f\cont{n} g$ satisfies
$f\cont
{n} g = g\cont{n} f = \langle f,g^\ast\rangle_{L^2(\R^n)}$.
\end{longlist}
\end{lemma}

 Contractions provide a useful tool for decomposing products
of stochastic integrals, in precise analogy to the classical context.
The following is \cite{BianeSpeicher}, Proposition~5.3.3.

\begin{proposition}[(Biane--Speicher)] \label{propproductcont} Let
$n,m$ be natural numbers, and let $f\in L^2(\R_+^n)$ and $g\in L^2(\R
_+^m)$. Then
%
\begin{equation} \label{eqproductcont} I^S_n(f)\cdot I^S_m(g) = \sum
_{p=0}^{\min\{n,m\}} I^S_{n+m-2p}(f\cont{p}g).
\end{equation}
\end{proposition}

\begin{remark} In the Gaussian Wiener chaos, a similar though more
complicated product formula holds.
\[
I^W_n(f)\cdot I^W_m(g) = \sum_{p=0}^{\min\{n,m\}}p!\pmatrix
{n\cr p}\pmatrix
{m\cr p} I^W_{n+m-2p}(f\cont{p}g).
\]
It is common for formulas from classical probability to have free
probabilistic analogs with simpler forms, usually with binomial
coefficients removed. This can be understood in terms of the relevant
(smaller) class of partitions that control moments in the theory, as we
discuss in Section~\ref{sectionNCpartitions} below.
\end{remark}

\subsection{Noncrossing partitions} \label{sectionNCpartitions}
\begin{figure}

\includegraphics{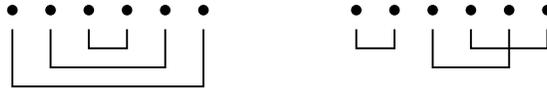}

\caption{Two pairings of $[6]=\{
1,2,3,4,5,6\}$. The first (totally-nested) pairing is noncrossing,
while the second is not.}\label{figpairingsof[6]}
\end{figure}

Proposition~\ref{propproductcont} shows that contractions are
involved in the algebraic structure of the space of stochastic
integrals. Since contractions involve integrals pairing different
classes of indices, general moments of stochastic integrals are best
understood in terms of a more abstract description of these pairings.
For convenience, we write $[n]$ to represent the set $[n]\equiv\{
1,2,\ldots,n\}$ for any positive integer $n$. If $n$ is even, then a
\textit{pairing} or \textit{matching} of $[n]$ is a partition of $[n]$ into
$n/2$ disjoint subsets each of size $2$. For example, $\{\{1,6\},\{2,5\}
,\{3,4\}\}$ and $\{\{1,2\},\{3,5\},\{4,6\}\}$ are two pairings of
$[6]=\{1,2,3,4,5,6\}$. It is convenient to represent such pairings
graphically, as in Figure~\ref{figpairingsof[6]}.

%

 It will be convenient to allow for more general partitions in
the sequel. A \textit{partition} of $[n]$ is (as the name suggests) a
collection of mutually disjoint nonempty subsets $B_1,\ldots,B_r$ of
$[n]$ such that $B_1\sqcup\cdots\sqcup B_r = [n]$. The subsets are
called the \textit{blocks} of the partition. By convention we order the
blocks by their least elements; that is, $\min B_i < \min B_j$ iff
$i<j$. The set of all partitions on $[n]$ is denoted $\mathscr{P}(n)$,
and the subset of all pairings is $\mathscr{P}_2(n)$.

\begin{definition}\label{defNC} Let $\pi\in\mathscr{P}(n)$ be a
partition of $[n]$. We say $\pi$ has a~\textit{crossing} if there are two
distinct blocks $B_1,B_2$ in $\pi$ with elements \mbox{$x_1,y_1\in B_1$} and
$x_2,y_2\in B_2$ such that $x_1<x_2<y_1<y_2$. (This is demonstrated in
Figure~\ref{figpairingsof[6]}.)

 If $\pi\in\mathscr{P}(n)$ has no crossings, it is said to be
a \textit{noncrossing partition}. The set of noncrossing partitions of
$[n]$ is denoted $NC(n)$. The subset of noncrossing pairings is denoted
$NC_2(n)$.
\end{definition}

The reader is referred to \cite{NicaSpeicherBook} for an extremely
in-depth discussion of the algebraic and enumerative properties of the
lattices $NC(n)$. For our purposes, we present only those structural
features that will be needed in the analysis of Wigner integrals.

\begin{figure}

\includegraphics{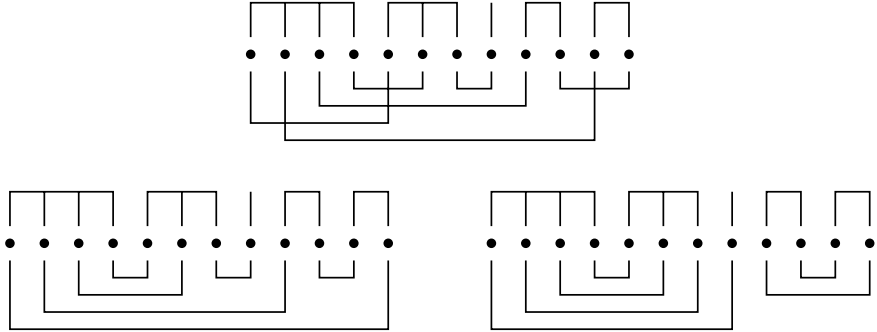}

\caption{The partition $4\tensor3\tensor
1\tensor
2\tensor2$ is drawn above the dots; below are three pairings that
respect it. The two bottom pairings are in $NC_2(4\tensor3\tensor
1\tensor2\tensor2)$.}\label{figrespect}
\end{figure}

\begin{definition}\label{defrespect}
$\!\!\!$Let $n_1,\ldots,n_r$ be
positive integers with \mbox{$n=n_1+\cdots+n_r$}. The set $[n]$ is then
partitioned accordingly as $[n] = B_1\sqcup\cdots\sqcup B_r$ where $B_1
= \{1,\ldots,n_1\}$, $B_2 = \{n_1+1,\ldots,n_1+n_2\}$, and so forth
through $B_r = \{n_1+\cdots+n_{r-1}+1,\ldots,n_1+\cdots+n_r\}$. Denote
this partition as $n_1\tensor\cdots\tensor n_r$.

 Say that a pairing $\pi\in\mathscr{P}_2(n)$ \textit{respects
$n_1\tensor\cdots\tensor n_r$} if no block of $\pi$ contains more than
one element from any given block of $n_1\tensor\cdots\tensor n_r$.
(This is demonstrated in Figure~\ref{figrespect}.) The set of such
respectful pairings is denoted $\mathscr{P}_2(n_1\tensor\cdots
\tensor
n_r)$. The set of noncrossing pairings that respect $n_1\tensor\cdots
\tensor n_r$ is denoted $NC_2(n_1\tensor\cdots\tensor n_r)$.
\end{definition}

 Partitions $n_1\tensor\cdots\tensor n_r$ as described in
Definition~\ref{defrespect} are called \textit{interval partitions},
since all of their blocks are intervals. Figure~\ref{figrespect} gives
some examples of respectful pairings.

%

\begin{remark} The same definition of \textit{respectful} makes perfect
sense for more general partitions, but we will not have occasion to use
it for anything but pairings. However, see
Remark~\ref{rksupinf}.\vspace*{-2pt}
\end{remark}

\begin{remark} \label{rknonflat} Consider the partition $n_1\otimes
\cdot\cdot\cdot\otimes n_r = \{B_1,\ldots,B_r\}$, as well as a pairing
$\pi\in\mathscr{P}_2(n) $, where $n = n_1+\cdots+n_r$. In the
classical literature about Gaussian subordinated random fields (cf.
\cite{PecTaqbook}, Chapter~4, and the references therein) the pair
$(n_1\otimes\cdot\cdot\cdot\otimes n_r, \pi)$ is represented
graphically as follows: (i) draw the blocks $B_1,\ldots,B_r$ as superposed
rows of dots (the $i$th row containing exactly $n_i$ dots,
$i=1,\ldots,r$), and (ii) join two dots with an edge if and only if the
corresponding two elements constitute a block of $\pi$. The graph thus
obtained is customarily called a \textit{Gaussian diagram}. Moreover, if
$\pi$ respects $n_1\tensor\cdots\tensor n_r$ according to Definition
\ref{defrespect}, then the Gaussian diagram is said to be \textit{nonflat}, in the sense that all its edges join different horizontal
lines, and therefore are not \textit{flat}, that is, not horizontal. The
noncrossing condition is difficult to discern from the Gaussian diagram
representation, which is why we do not use it here; therefore the
nonflat terminology is less meaningful for us, and we prefer the
intuitive notation from Definition~\ref{defrespect}.\vspace*{-2pt}
\end{remark}

One more property of pairings will be necessary in the proceeding analysis.\vspace*{-2pt}

\begin{definition}\label{defconnected} Let $n_1,\ldots,n_r$ be
positive integers, and let $\pi\in\mathscr{P}_2(n_1\tensor\cdots
\tensor n_r)$. Let $B_1$, $B_2$ be two blocks in $n_1\tensor\cdots
\tensor n_r$. Say that \textit{$\pi$ links $B_1$ and $B_2$} if there is a
block $\{i,j\}\in\pi$ such that $i\in B_1$ and $j\in B_2$.

 Define a graph $C_\pi$ whose vertices are the blocks of
$n_1\tensor\cdots\tensor n_r$; $C_\pi$ has an edge between $B_1$ and
$B_2$ iff $\pi$ links $B_1$ and $B_2$. Say that \textit{$\pi$ is
connected} with respect to $n_1\tensor\cdots\tensor n_r$ (or that
\textit{$\pi$ connects the blocks of $n_1\tensor\cdots\tensor n_r$}) if the
graph $C_\pi$ is connected.

 Denote by $NC_2^c(n_1\tensor\cdots\tensor n_r)$ the set of
noncrossing pairings that both respect and connect $n_1\tensor\cdots
\tensor n_r$.\vadjust{\goodbreak}
\end{definition}

 For example, the second partition in Figure~\ref{figrespect}
is in $NC_2^c(4\tensor3\tensor1\tensor2\tensor2)$, while the third
is not. The interested reader may like to check that $NC_2(4\tensor
3\tensor1\tensor2\tensor2)$ has $5$ elements, and all are connected
except the third example in Figure~\ref{figrespect}.\vspace*{-2pt}

\begin{remark} \label{rksupinf} For a positive integer $n$, the set
$NC(n)$ of noncrossing partitions on $[n]$ is a lattice whose partial
order is given by reverse refinement. The top element $1_n$ is the
partition $\{\{1,\ldots,n\}\}$ containing only one block; the bottom
element $0_n$ is $\{\{1\},\ldots,\{n\}\}$ consisting of $n$ singletons.
The conditions of Definitions~\ref{defrespect} and~\ref{defconnected}
can be described elegantly in terms of the lattice operations meet
$\wedge$ (i.e., inf) and join $\vee$ (i.e., sup). If $n=n_1+\cdots
+n_r$, then $\pi\in NC_2(n)$ respects $n_1\tensor\cdots\tensor n_r$ if
and only if $\pi\wedge(n_1\tensor\cdots\tensor n_r) = 0_n$; $\pi$
connects the blocks of $n_1\tensor\cdots\tensor n_r$ if and only if
$\pi
\vee(n_1\tensor\cdots\tensor n_r) = 1_n$.\vspace*{-2pt}
\end{remark}

\begin{remark} \label{rkdecomposition} Given $n_1,\ldots,n_r$ and a
respectful noncrossing pairing $\pi\in NC_2(n_1\tensor\cdots\tensor
n_r)$, there is a unique decomposition of the full index set $[n]$,
where $n=n_1+\cdots+n_r$, into subsets $D_1,\ldots,D_m$ of the blocks
of $n_1\tensor\cdots\tensor n_r$, such that the restriction of $\pi$ to
each $D_i$ connects the blocks of $D_i$. These $D_i$ are the vertices
of the graph $C_\pi$ grouped according to connected components of the
graph. For example, in the third pairing in Figure~\ref{figrespect},
the decomposition has two components, $D_1 = 4\tensor3\tensor1$ and
$D_2 = 2\tensor2$. To be clear, this notation is slightly misleading
since the $2\tensor2$ in this case represents indices $\{9,10\},\{
11,12\}$, not $\{1,2\},\{3,4\}$; we will be a little sloppy about this
to make the following much more readable.\vspace*{-2pt}
\end{remark}

There is a close connection between respectful noncrossing pairings and
expectations of products of Wigner integrals. To see this, we first
introduce an action of pairings on functions.\vspace*{-2pt}

\begin{definition} \label{defpairingintegral} Let $n$ be an even
integer, and let $\pi\in\mathscr{P}_2(n)$. Let $f\dvtx \R_+^n\to
\C$ be
measurable. The \textit{pairing integral} of $f$ with respect to $\pi$,
denoted $\int_\pi f$, is defined (when it exists) to be the constant
\[
\int_\pi f = \int f(t_1,\ldots,t_n) \prod_{\{i,j\}\in\pi} \delta
(t_i-t_j)\,dt_1\cdots \,dt_n.\vspace*{-2pt}
\]
\end{definition}

 For example, given the second pairing $\pi=\{\{1,2\},\{3,5\}
,\{4,6\}\}$ in Figure~\ref{figpairingsof[6]},
\[
\int_\pi f = \int_{\R_+^3} f(r,r,s,t,s,t)\,dr\,ds\,dt.\vspace*{-2pt}
\]

\begin{remark} The operation $\int_\pi$ is not well defined on
$L^2(\R
_+^n)$; for example, if $n=2$ and $\pi=\{\{1,2\}\}$, then $\int_\pi f$
is finite if and only if $f$ is the kernel\vadjust{\goodbreak} of a trace class
Hilbert--Schmidt operator on $L^2(\R_+)$. However, it is easy to see
that $\int_\pi f$ is well-defined whenever $f$ is a tensor product of
functions, and $\pi$ respects the interval partition induced by this
tensor product (cf. Lemma~\ref{lemmaCS}). (This is one of the reasons
why one should interpret multiple stochastic integrals as integrals on
product spaces without diagonals, since integrals on diagonals are in
general not defined.) This is precisely the case we will deal with in
all of the following.\vspace*{-2pt}
\end{remark}

Note that a contraction $f\cont{p}g$ can be interpreted in terms of a
pairing integral, using a \textit{partial pairing}, that is, one that
pairs only a subset of the indices. If $f\in L^2(\R_+^n)$ and $g\in
L^2(\R_+^m)$, and $p\le\min\{n,m\}$ is a natural number, then
\[
f\cont{p}g = \int_{\tau_p} f\tensor g,
\]
where $\tau_p$ is the partial pairing $\{\{n,n+1\},\{n-1,n+2\},\ldots
,\{
n-p+1,n+p\}\}$ of $[n+m]$.

%
\begin{figure}

\includegraphics{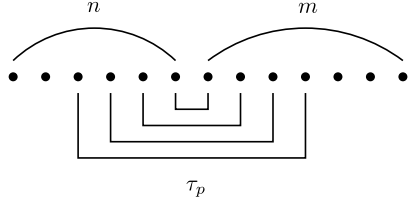}

\caption{A partial pairing $\tau_p$ of
$[n+m]$ corresponding to a $p$-contraction; here $n=6$, $m=7$, and $p=4$.}
\label{figpartialpairingcont}\vspace*{-2pt}
\end{figure}

%
\begin{figure}[b]
\vspace*{-2pt}
\includegraphics{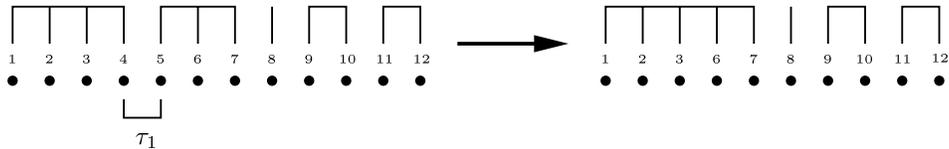}

\caption{The partial pairing $\tau
_1$ acts on the left on $4\tensor3\tensor1\tensor2\tensor2$, joining
the first two blocks and deleting the middle indices, to produce the
partition $5\tensor1\tensor2\tensor2$. The indices are labeled to make
the action clearer.}\label{figpartialpairingaction}
\end{figure}

The partial contraction pairings $\tau_p$ provide a useful
decomposition of the set of all respectful noncrossing pairings, in the
following sense. Let $n_1,\ldots,n_r$ be positive integers. If $p\le
\min\{n_1,n_2\}$, the partial pairing $\tau_p$ acts (on the left) on
the partition $n_1\tensor n_2\tensor n_3\tensor\cdots\tensor n_r$ to
produce the partition $(n_1+n_2-2p)\tensor n_3\tensor\cdots\tensor
n_r$. That is, $\tau_p$ joins the first two blocks of $n_1\tensor
\cdots
\tensor n_r$ and deletes the paired indices to produce a new interval
partition. This is demonstrated in Figure~\ref{figpartialpairingaction}.

\begin{figure}

\includegraphics{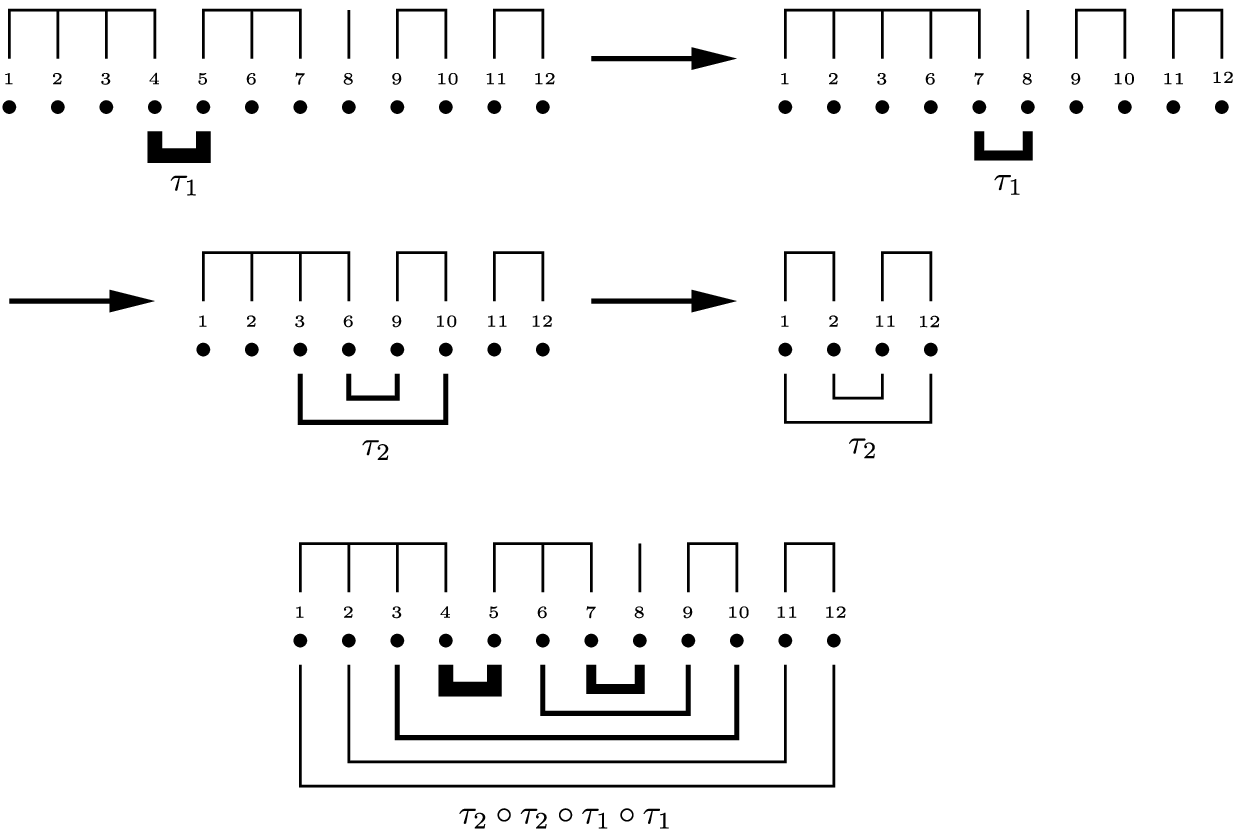}%
\vspace*{-2pt}
\caption{The composition $\tau
_2\circ
\tau_2\circ\tau_1\circ\tau_1$ produces a noncrossing pairing that
respects $4\tensor3\tensor1\tensor2\tensor2$.}\label{figpairingdecomposition}\vspace*{-2pt}
\end{figure}

Considered as such a function, we may then compose partial contraction
pairings. For example,\vadjust{\goodbreak} following Figure~\ref{figpartialpairingaction},
 we may act again with $\tau_1$ on $5\tensor1\tensor2\tensor2$
to yield $4\tensor2\tensor2$; then with $\tau_2$ to get $2\tensor2$;
and finally $\tau_2$ maps this partition to the empty partition.
Stringing these together gives a respectful pairing of the original
interval partition, which we denote $\tau_2\circ\tau_2\circ\tau
_1\circ
\tau_1$. Figure~\ref{figpairingdecomposition} displays this composition.

%

To be clear: we start from the left and then do the partial pairing
$\tau_p$ between the first and second block; after this application,
the (rest of the) first and second blocks are treated as a single
block. This is still the case if $p=0$; here there are no paired
indices, but the action of $\tau_0$ records the fact that, for further
discussion, the first two blocks are now connected. An example is given
in Figure~\ref{figtau0} below, where the action of $\tau_0$ is
graphically represented by a dashed line.

%
\begin{figure}

\includegraphics{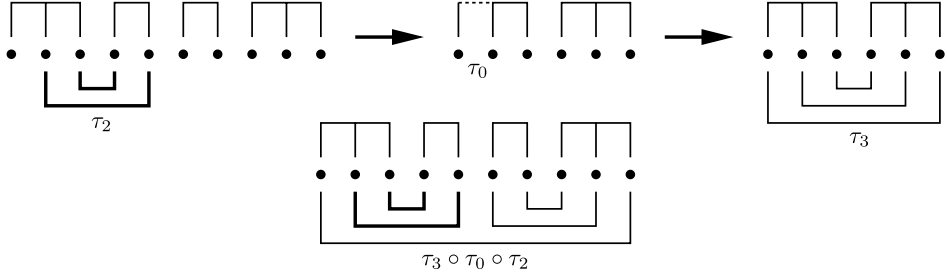}%
\vspace*{-2pt}
\caption{The pairing $\pi=\{\{1,10\},\{2,5\},\{
3,4\}
,\{6,9\},\{7,8\}\}$ respects the interval partition $3\tensor2\tensor
2\tensor3$. Its decomposition is given by $\pi= \tau_3\circ\tau
_0\circ
\tau_2$.}\label{figtau0}\vspace*{-2pt}
\end{figure}

With this convention, further $\tau_p$ may act only on the first two
blocks, which results in a \textit{unique} decomposition of any respectful
pairing into partial contractions, as the next lemma makes clear.\vspace*{-2pt}


\begin{lemma} \label{lemmacontractiondecomposition}
Let $n_1,\ldots,n_r$ be positive integers, and let $\pi\in
NC_2(n_1\tensor\cdots\tensor n_r)$. There is a unique sequence of
partial contractions $\tau_{p_1},\ldots,\tau_{p_{r-1}}$ such that
$\pi
= \tau_{p_{r-1}}\circ\cdots\circ\tau_{p_1}$.\vspace*{-2pt}
\end{lemma}


\begin{pf} Any noncrossing pairing must contain an interval $\{
i,i+1\}$; cf.~\cite{NicaSpeicherBook}, Remark~9.2(2). Hence, since
$\pi
$ respects $n_1\tensor\cdots\tensor n_r=\{B_1,\ldots,B_r\}$, there
must\vadjust{\goodbreak}
be two \textit{adjacent} blocks linked by $\pi$. Let $j\in[k]$ be the
smallest index for which $B_j,B_{j+1}$ are connected by $\pi$; hence
all of the blocks $B_1,\ldots,B_j$ pair among the blocks
$B_{j+1},\ldots
,B_r$. Note that any partition that satisfies this constraint and also
respects the coarser interval partition $(n_1+\cdots+n_j)\tensor
n_{j+1}\tensor\cdots\tensor n_r$ is automatically in $NC_2(n_1\tensor
\cdots\tensor n_r)$. In other words, we can begin by decomposing $\pi=
\pi'\circ(\tau_0)^{j-1}$, where $\pi'\in NC_2((n_1+\cdots
+n_j)\tensor
n_{j+1}\tensor\cdots\tensor n_r)$ links the first and second blocks of
this interval partition. By construction, this $j$ is unique.

Let $n_0=n_1+\cdots+n_j$, so $\pi'$ links $\{1,\ldots,n_0\}$ with $\{
n_0+1,\ldots,n_0+n_{j+1}\}$. It follows that $\{n_0,n_0+1\}\in\pi'$:
for if $n_0$ pairs with some element $n_0+i$ with $i\ge2$, then
$n_0+1,\ldots,n_0+i-1$ cannot pair anywhere 
without introducing crossings. Following these lines, an easy induction
shows that there is some $p\in[\min\{n_0,n_{j+1}\}]$ such that the
pairs $\{n_0,n_0+1\},\{n_0-1,n_0+2\},\ldots,\{n_0-p+1,n_0+p\}$ are in
$\pi'$, while all indices $1,\ldots,n_0-p$ and $n_0+p+1,\ldots
,n_0+n_{j+1}$ pair outside $[n_0+n_{j+1}]$. In other words, $\pi' =
\pi
''\circ\tau_p$ for some noncrossing pairing $\pi''$ that respects
$(n_0-p)\tensor(n_{j+1}-p)\tensor n_3\tensor\cdots\tensor n_r$. What's
more, since $p$ was chosen maximally so that there are no further
pairings in the blocks $(n_0-p)\tensor(n_{j+1}-p)$, these two may be
treated as a single block, and $\pi''$ is only constrained to be in
$NC_2(n_0+n_{j+1}-2p,n_3,\ldots,n_r)$. Since $p>0$, the lemma now
follows by a finite induction; uniqueness results from the left-most
choice of $j$ and maximal choice of $p$ at each stage.\vspace*{-2pt}
\end{pf}

By carefully tracking the proof of Lemma~\ref{lemmacontractiondecomposition},
we can give a complete description of the class of
respectful pairings in terms of their decompositions.\vspace*{-2pt}

\begin{lemma} \label{lemmacontractionenumeration} Let $n_1,\ldots
,n_r$ be positive integers. The class $NC_2(n_1\tensor\cdots\tensor
n_r)$ is equal to the set of compositions $\tau_{p_{r-1}}\circ\cdots
\circ\tau_{p_1}$ where $(p_1,\ldots,p_{r-1})$ satisfy the inequalities
%
\begin{eqnarray} \label{eqrespectfulineq}
 0&\le& p_1\le\min\{n_2,n_1\}, \nonumber\\
0&\le& p_k\le\min\{n_{k+1},n_1+\cdots+n_k-2p_1-\cdots-2p_{k-1}\},\\
\eqntext{1<k<r-1,
2 (p_1+\cdots+p_{r-1}) =n_1+\cdots+n_r.}
\end{eqnarray}
\end{lemma}

 Inequalities (\ref{eqrespectfulineq}) in Lemma~\ref{lemmacontractionenumeration} successively guarantee that the partial
contractions $\tau_{p_k}$ in the decomposition of $\pi$ only contract
elements from within two adjacent blocks; the final equality is to
guarantee that all indices are paired in the end. Since every
respectful pairing has a contraction decomposition, and each
contraction decomposition satisfying inequalities~(\ref{eqrespectfulineq}) is respectful
(a fact which follows from an easy induction),
these inequalities define $NC_2(n_1\tensor\cdots\tensor n_r)$. This
completely combinatorial description would be the starting point for an
enumeration of the class of respectful pairings; however, even in the
case $n_1=\cdots=n_r$, the enumeration appears to be extremely difficult.

We conclude this section with a proposition that demonstrates the
efficacy of pairing integrals and noncrossing pairings in the analysis
of Wigner integrals.

\begin{proposition} \label{propmomentsofproducts} Let $n_1,\ldots
,n_r$ be positive integers, and suppose $f_1,\ldots,\allowbreak f_r$ are functions
with $f_i\in L^2(\R_+^{n_i})$ for $1\le i\le r$. The expectation $\ff$
of the\vspace*{1.5pt} product of Wigner integrals $I_{n_1}^S(f_1)\cdots
I_{n_r}^S(f_r)$ is given by
%
\begin{equation} \label{eqmomentsofproducts} \ff [
I_{n_1}^S(f_1)\cdots I_{n_r}^S(f_r) ] = \sum_{\pi\in
NC_2(n_1\tensor\cdots\tensor n_r)} \int_\pi f_1\tensor\cdots
\tensor
f_r.
\end{equation}
\end{proposition}

\begin{remark} This result has been used in the literature (e.g., to
prove~\cite{BianeSpeicher}, Theorem~5.3.4), but it appears to have a
folklore status in that a~proof has not been written down. The
following proof is an easy application of Proposition~\ref{propproductcont}, together with Lemma
\ref{lemmacontractionenumeration}.
\end{remark}

\begin{pf} By iterating equation (\ref{eqproductcont}), we arrive at
the following unwieldy expression. (For readability, we have hidden the
explicit dependence of the Wigner integral $I^S_n$ on the number of
variables $n$ in its argument.)
%
\begin{equation} \label{eqiteratedproductcont}
I^S(f_1)\cdots I^S(f_r) = \sum_{p_{r-1}} \cdots\sum_{p_1} I^S
\bigl( \bigl(\cdots \bigl( (f_1\cont{p_1}f_2 )\cont
{p_2}f_3
\bigr)\cdots \bigr)\cont{p_{r-1}}f_r \bigr),\hspace*{-35pt}
\end{equation}
where $p_1,\ldots,p_{r-1}$ range over the set specified by the first
two inequalities in equation (\ref{eqrespectfulineq}). (This is the
range of the $p_k$ for the same reason that those inequalities specify
the range of the $p_k$ for contraction decompositions: the first two
inequalities in (\ref{eqrespectfulineq}) merely guarantee that
contractions are performed, successively, only between two adjacent
blocks of $n_1\tensor\cdots\tensor n_r$.) Note: following Remark~\ref{rkcontractionsnotassociative}, the order the contractions are
performed in equation (\ref{eqiteratedproductcont}) is important.

Taking expectation in equation (\ref{eqiteratedproductcont}), note
that most terms have $\ff=0$ since any nontrivial stochastic integral
is centred (as it is orthogonal to constants in the $0$th order of
chaos). Hence, the only terms that contribute to the sum are those for
which the iterated contractions pair all indices of the functions; that
is, the sum is over those $p_1,\ldots,p_{r-1}$ for which $2(p_1+\cdots
+p_{r-1}) = n_1+\cdots+n_r$, so that the stochastic integral $I^S$ in
the sum is $I^S_0$. Since such a trivial stochastic integral is just
the identity on the constant function inside, this shows that
\[
\ff [ I^S(f_1)\cdots I^S(f_r) ] = \sum_{p_{r-1}}\cdots
\sum
_{p_1} \bigl( \bigl(\cdots \bigl( (f_1\cont{p_1}f_2
)\cont
{p_2}f_3 \bigr)\cdots \bigr)\cont{p_{r-1}}f_r \bigr),
\]
where the sum is over those $p_1,\ldots,p_{r-1}$ satisfying the same
inequalities mentioned above, along with the condition \mbox{$2(p_1+\cdots
+p_{r-1}) = n_1+\cdots+n_r$}; that is, the $p_k$ satisfy inequalities
(\ref{eqrespectfulineq}). Each such iterated contraction integral
corresponds to a pairing integral of $f_1\tensor\cdots\tensor f_r$ in
the obvious fashion,
\[
 \bigl( \bigl(\cdots \bigl( (f_1\cont{p_1}f_2 )\cont
{p_2}f_3
\bigr)\cdots \bigr)\cont{p_{r-1}}f_r \bigr) = \int_{\tau
_{p_{r-1}}\circ\cdots
\circ\tau_{p_1}} f_1\tensor\cdots\tensor f_r.
\]
Lemma~\ref{lemmacontractionenumeration} therefore completes the
proof.\vspace*{-2pt}
\end{pf}

\begin{remark} Another proof of Proposition~\ref{propmomentsofproducts} can be achieved using a random matrix approximation to the
free Brownian motion, as discussed in Section~\ref{sectionfBM}. The
starting point is the classical counterpoint to Proposition~\ref{propmomentsofproducts} \cite{Janson}, Theorem~7.33, which states that the
expectation of a product of Wiener integrals is a similar sum of
pairing integrals over respectful (i.e., nonflat) pairings, but in
this case crossing pairings must also be included. Modifying this
formula for matrix-valued Brownian motion, and controlling the leading
terms in the limit as matrix size tends to infinity using the so-called
``genus expansion,'' leads to equation (\ref{eqmomentsofproducts}).
The (quite involved) details are left to the interested reader.\vspace*{-2pt}
\end{remark}

\section{Central limit theorems} \label{sectionCLTs}

 We begin by proving Theorem~\ref{thm4thmomentcont}, which
we restate here for convenience.\vspace*{-2pt}

\begin{theoremcont*} Let $n$ be a natural number, and let $(f_k)_{k\in
\N
}$ be a sequence of functions
in $L^2(\R_+^n)$, each with $\|f_k\|_{L^2(\R_+^n)} = 1$. The following
statements are equivalent:
\begin{longlist}[(2)]
\item[(1)] The fourth absolute moments of the stochastic integrals
$I^S_n(f_k)$ converge to~$2$.
\[
\lim_{k\to\infty} \ff(|I^S_n(f_k)|^4) = 2.
\]
\item[(2)] All nontrivial contractions of $f_k$ converge to $0$. For
each $p=1,2,\ldots,\allowbreak n-1$,
\[
\lim_{k\to\infty} f_k\cont{p} f_k^\ast= 0 \qquad\mbox{in }
L^2(\R_+^{2n-2p}).\vspace*{-2pt}
\]
\end{longlist}
\end{theoremcont*}

\begin{pf} The expression $|I^S_n(f_k)|^4$ is short-hand for
$[I_n^S(f_k)\cdot I_n^S(f_k)^\ast]^2$. Since [according to equation
(\ref{eqintadjoint})] $I_n^S(f_k)^\ast= I_n^S(f_k^\ast)$, this is
a\vadjust{\goodbreak}
product of Wigner integrals, to which we will apply Proposition~\ref{propproductcont}. First,
%
\begin{equation} \label{eq4thmomentcont1} I_n^S(f_k)\cdot
I_n^S(f_k^\ast) = \sum_{p=0}^n I_{2n-2p}^S(f_k\cont{p}f_k^\ast).
\end{equation}
The Wigner integrals on the right-hand side of equation (\ref{eq4thmomentcont1}) are in different orders of chaos, and hence are
orthogonal (with respect to the $\ff$-inner product). Thus, we can expand
\begin{eqnarray*}
\ff(|I_n^S(f_k)|^4) &=& \ff \bigl[\bigl(I_n^S(f_k)\cdot
I_n^S(f_k^\ast)\bigr)^2 \bigr] \\
&=& \langle I_n^S(f_k)\cdot I_n^S(f_k^\ast) , I_n^S(f_k)\cdot
I_n^S(f_k^\ast)\rangle_\ff\\
&=& \sum_{p=0}^n \langle I_{2n-2p}^S(f_k\cont{p}f_k^\ast) ,
I_{2n-2p}^S(f_k\cont{p}f_k^\ast) \rangle_\ff,
\end{eqnarray*}
where in the second equality we have used the fact that
$I_n^S(f_k)\cdot I_n^S(f_k^\ast)$ is self adjoint. Now employing the
Wigner isometry [equation (\ref{eqWignerisometry})], this yields
%
\begin{equation} \label{eq4thmomentcont2} \ff(|I_n^S(f_k)|^4) =
\sum
_{p=0}^n \langle f_k\cont{p} f_k^\ast, f_k\cont{p} f_k^\ast\rangle
_{L^2(\R_+^{2n-2p})}.
\end{equation}
Consider first the two boundary terms in the sum in equation (\ref{eq4thmomentcont2}). When $p=n$, we have
\[
f_k\cont{n}f_k^\ast= \langle f_k,f_k\rangle_{L^2(\R_+^n)} = 1,
\]
according to Lemma~\ref{lemmacontadj}(2) and the assumption that
$f_k$ is normalized in $L^2$. On the other hand, when $p=0$, the
contraction $f_k\cont{0}f_k$ is just the tensor product $f\tensor
f^\ast
$, and we have
\[
\langle f_k\tensor f_k^\ast,f_k\tensor f_k^\ast\rangle_{L^2(\R_+^{2n})}
= \langle f_k,f_k\rangle_{L^2(\R_+^n)}\langle f_k^\ast,f_k^\ast
\rangle
_{L^2(\R_+^n)}=1.
\]
(Both terms in the product are equal to $\|f_k\|_{L^2}^2=1$, following
Definition~\ref{defmirror} of $f_k^\ast$.) Equation (\ref{eq4thmomentcont2}) can therefore be rewritten as
%
\begin{equation} \label{eq4thmomentcont3} \ff(|I_n^S(f_k)|^4) = 2 +
\sum_{p=1}^{n-1} \| f_k\cont{p} f_k^\ast \|_{L^2(\R
_+^{2n-2p})}^2.
\end{equation}
Thus, the statement that the limit of $\ff(|I_n^S(f_k)|^4)$ equals $2$
is equivalent to the statement that the limit of the sum on the
right-hand side of equation~(\ref{eq4thmomentcont3}) is $0$. This is
a sum of nonnegative terms, and so each of the terms must have limit
$0$. This completes the proof.
\end{pf}

 Corollary~\ref{cornosemicircular} now follows quite easily.

\begin{corollarynosemicirc*} Let $n\ge2$ be an integer, and consider a
nonzero mirror symmetric function $f\in L^2(\R_+^n)$. Then the Wigner
integral $I_n^S(f)$\vadjust{\goodbreak} satisfies $\ff[I_n^S(f)^4] > 2\ff[I_n^S(f)^2]^2$.
In particular, the distribution of the Wigner integral $I^S_n(f)$
cannot be semicircular.
\end{corollarynosemicirc*}

\begin{pf} By rescaling, we may assume that $\|f\|_{L^2(\R_+^n)} =
1$; in this case, equation (\ref{eq4thmomentcont3}) shows that $\ff
[I_n^S(f)^4]\ge2\ff[I_n^S(f)^2]^2$. To achieve a contradiction, we
assume that $\ff[I_n^S(f)^4] = 2\ff[I_n^S(f)^2]^2 = 2$ [which would be
the case if $I_n^S(f)$ were semicircular]. Then the constant sequence
$f_k=f$ for all $k$ satisfies condition (1) of Theorem~\ref{thm4thmomentcont}; hence, for $1\le p\le n-1$,
\[
f\cont{p}f^\ast= \lim_{k\to\infty} f_k\cont{p}f_k^\ast= 0 \qquad
\mbox
{in } L^2(\R_+^{2n-2p}).
\]
Take, for example, $p=n-1$. Let $g\in L^2(\R_+)$, so that $g\tensor
g^\ast\in L^2(\R_+^2)$. Then we may calculate the inner product
\begin{eqnarray*} &&\langle f \cont{n-1}f^\ast,g\tensor g^\ast\rangle
_{L^2(\R_+^2)}\\
 && \qquad = \int[f \cont{n-1}f^\ast](s,t)\overline
{[g\tensor
g^\ast](s,t)}\,ds\,dt \\
&& \qquad = \int \biggl(\int f(s,s_2,\ldots,s_n)f^\ast(s_n,\ldots,s_2,t)\,
ds_2\cdots \,ds_n \biggr)\overline{g(s)}g(t)\,ds\,dt \\
&& \qquad = \int g^\ast(s) f(s,s_2,\ldots,s_n) \cdot\overline{g^\ast(t)
f(t,s_2,\ldots,s_n)}\, ds\,dt\,ds_2\cdots \,ds_n \\
&& \qquad = \| g^\ast\cont{1} f\|_{L^2(\R_+^{n-1})}^2.
\end{eqnarray*}
%
By assumption, $f \cont{n-1}f^\ast= 0$, and so we have $g^\ast\cont
{1} f = 0$ for all $g\in L^2(\R_+)$. That is, for almost all
$s_2,\ldots
,s_n\in\R_+$,
\[
\int_0^\infty\overline{g(s)}f(s,s_2,\ldots,s_n)\,ds = 0.
\]
For fixed $s_2,\ldots,s_n$ for which this holds, taking $g$ to be the
function $g(s) = f(s,s_2,\ldots,s_n)$ yields that $f(s,s_2,\ldots,s_n)
= 0$ for almost all $s$. Hence, $f=0$ almost surely. This contradicts
the normalization $\|f\|_{L^2(\R_+^n)} = 1$.
\end{pf}


We now proceed towards the proof of Theorem~\ref{thm4thmomentsemicircle}. First,
we state a~technical result that will be of use.

\begin{lemma} \label{lemmaCS} Let $n_1,\ldots,n_r$ be positive
integers, and let $f_i\in L^2(\R_+^{n_i})$ for $1\le i\le r$. Let $\pi$
be a pairing in $\mathscr{P}_2(n_1\tensor\cdots\tensor n_r)$. Then
\[
 \biggl|\int_\pi f_1\tensor\cdots\tensor f_r \biggr| \le\|f_1\|
_{L^2(\R
_+^{n_1})} \cdots\|f_r\|_{L^2(\R_+^{n_r})}.
\]
\end{lemma}

\begin{pf} This follows by iterated application of the
Cauchy--Schwarz inequality along the pairs in $\pi$. It is proved as
\cite{Janson}, Lemma~7.31.\vadjust{\goodbreak}
\end{pf}

The following proposition shows that contractions control all important
pairing integrals.

\begin{proposition} \label{propvanishingpairingintegrals} Let $n$ be
a positive integer. Consider a sequence~$(f_k)_{k\in\N}$ with $f_k\in
L^2(\R_+^n)$ for all $k$, such that:
\begin{longlist}[(3)]
\item[(1)] $f_k=f_k^\ast$ for all $k$;
\item[(2)] there is a constant $M>0$ such that $\|f_k\|_{L^2(\R
_+^n)}\le M$ for all $k$;
\item[(3)] for each $p=1,2,\ldots,n-1$,
\[
\lim_{k\to\infty} f_k\cont{p} f_k^\ast= 0 \qquad\mbox{in }
L^2(\R
_+^{2n-2p}).
\]
\end{longlist}
Let $r\ge3$, and let $\pi$ be a connected noncrossing pairing that
respects $n^{\tensor r}$: $\pi\in NC_2^c(n^{\tensor r})$; cf.
Definitions~\ref{defrespect} and~\ref{defconnected}. Then
\[
\lim_{k\to\infty} \int_\pi f_k^{\tensor r} = 0.
\]
\end{proposition}

\begin{figure}

\includegraphics{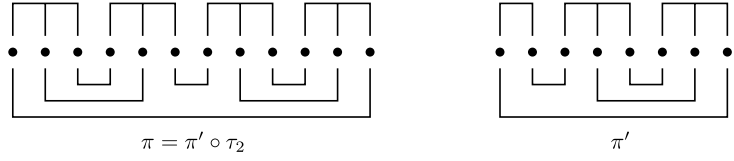}

\caption{A pairing $\pi\in
NC_2^c(3^{\tensor4})$, with the first step in its contraction
decomposition (per Lemma~\protect\ref{lemmacontractiondecomposition}).}\label{figexampleFubini}
\end{figure}

\begin{pf} Begin by decomposing $\pi= \tau_{p_{r-1}}\circ\cdots
\circ
\tau_{p_1}$ following Lemma~\ref{lemmacontractiondecomposition}.
There must be some nonzero $p_i$; to simplify notation, we assume that
$p_1>0$. (Otherwise we may perform a cyclic rotation and relabel
indices from the start.) Note also that, since $\pi$ connects the
blocks of $n^{\tensor r}$ and $r>2$, it follows that $p_1<n$: else the
first two blocks $\{1,\ldots,n\}$ and $\{n+1,\ldots,2n\}$ would form a
connected component in the graph $C_\pi$ from Definition~\ref{defconnected}, so $C_\pi$ would not be connected. Set $\pi' = \tau
_{p_k}\circ\cdots\circ\tau_{p_2}$, so that $\pi= \pi'\circ\tau
_{p_1}$. Then
(as in the proof of Proposition~\ref{propmomentsofproducts}) it
follows that
%
\begin{equation} \label{eqpairingintegral1} \int_\pi f_k^{\tensor r}
= \int_{\pi'} (f_k\cont{p_1} f_k)\tensor f_k^{\tensor(r-2)}.
\end{equation}
To make this clear, an example is given in Figure~\ref{figexampleFubini}, with the corresponding iterations of the integral in equation (\ref{eqFubini}).
%
\begin{eqnarray} \label{eqFubini}
\int_\pi f^{\tensor4} &=& \int_{\R_+^6}
f(t_1,t_2,t_3)f(t_3,t_2,t_4)f(t_4,t_5,t_6)f(t_6,t_5,t_1)\,
dt_1\,dt_2\,dt_3\,dt_4\,dt_5\,dt_6 \nonumber
\\
&=& \int_{\R_+^4} (f\cont{2}f)(t_1,t_4)f(t_4,t_5,t_6)f(t_6,t_5,t_1)\,
dt_1\,dt_4\,dt_5\,dt_6\\
 &=& \int_{\pi'} (f\cont{2} f)\tensor f^{\tensor2}.
\nonumber
\end{eqnarray}

 Employing Lemma~\ref{lemmaCS}, we therefore have
%
\begin{eqnarray}
 \biggl|\int_\pi f_k^{\tensor r} \biggr| &=& \biggl|\int_{\pi'}
(f_k\cont
{p_1} f_k)\tensor f_k^{\tensor(r-2)} \biggr| \nonumber\\
&\le&\|f_k\cont{p_1}f_k\|_{L^2(\R_+^{2n-2p})} \cdot\|f_k\|_{L^2(\R
_+^n)}^{r-2} \\
&\le&\|f_k\cont{p_1}f_k\|_{L^2(\R_+^{2n-2p})}\cdot M^{r-2},\nonumber
\end{eqnarray}
using assumption (2) in the proposition. But from assumptions (1) and~(3),
$\|f_k\cont{p_1}f_k\|_{L^2(\R_+^{2n-2p})}\to0$. The result
follows.
\end{pf}

%

We can now prove the main theorem of the paper, Theorem~\ref{thm4thmomentsemicircle}, which we restate here for convenience.

\begin{theoremlaw*} Let $n\ge2$ be an integer, and let $(f_k)_{k\in\N}$
be a sequence of mirror symmetric functions in $L^2(\R_+^n)$, each with
$\|f_k\|_{L^2(\R_+^n)} = 1$. The following statements are equivalent:
\begin{longlist}[(2)]
\item[(1)] The fourth moments of the stochastic integrals $I^S_n(f_k)$
converge to $2$.
\[
\lim_{k\to\infty} \ff(I^S_n(f_k)^4) = 2.
\]
\item[(2)] The random variables $I^S_n(f_k)$ converge in law to the
standard semicircular distribution $S(0,1)$ as $k\to\infty$.
\end{longlist}
\end{theoremlaw*}

\begin{pf}
As pointed out in Remark~\ref{rklawimpliesmoments}, the implication
(2)${}\Longrightarrow{}$(1) is essentially elementary: we need only demonstrate
uniform tail estimates. In fact, the laws $\mu_k$ of $I^S_n(f_k)$ are
all uniformly compactly-supported: by \cite{BianeSpeicher}, Theorem~5.3.4 (which is a version of the Haagerup inequality,
cf. \cite{Haagerup}), any Wigner integral satisfies
\[
 \| I_n^S(f) \| \le(n+1)\|f\|_{L^2(\R_+^n)}.
\]
Since all the functions $f_k$ are normalized in $L^2$, it follows that
$\supp \mu_k\subseteq[-n-1,n+1]$ for all $k$. Since the semicircle
law is also supported in this interval, we may approximate the function
$x\mapsto x^4$ by a $C_c(\R)$ function that agrees with it on all the
supports, and hence convergence in distribution of $\mu_k$ to the
semicircle law implies convergence of the fourth moments by definition.

We will use Proposition~\ref{propvanishingpairingintegrals},
together with Proposition~\ref{propmomentsofproducts}, to prove the
remarkable reverse implication. Since $S(0,1)$ is compactly supported,
it is enough to verify that the moments of $I^S_n(f_k)$ converge to the
moments of $S(0,1)$, as described following equation (\ref{eqsemicircle}).
Since $I^S_n(f_k)$ is orthogonal to the constant $1$ in
the first order of chaos, $I^S_n(f_k)$ is centred; the Wigner isometry
of equation (\ref{eqWignerisometry})\vadjust{\goodbreak} yields that the second moment of
$I^S_n(f_k)$ is constantly~$1$ due to normalization. Therefore, take
$r\ge3$. Proposition~\ref{propmomentsofproducts} yields
that\looseness=-1
%
\begin{equation} \label{eqmainproof1} \ff [I^S_n(f_k)^r
] =
\sum_{\pi\in NC_2(n^{\tensor r})} \int_\pi f_k^{\tensor r}.
\end{equation}\looseness=0
Following Remark~\ref{rkdecomposition}, any $\pi\in NC_2(n^{\tensor
r})$ can be (uniquely) decomposed into a disjoint union of connected
pairings $\pi= \pi_1\sqcup\cdots\sqcup\pi_m$ with $\pi_i\in
NC_2^c(n^{\tensor r_i})$ for some $r_i$'s with $r_1+\cdots+r_m = r$.
Since the decomposition respects the partition $n^{\tensor r}$, the
pairing integrals decompose as products.
%
\begin{equation} \label{eqmainproof2} \int_\pi f_k^{\tensor r} =
\prod_{i=1}^m \int_{\pi_i} f_k^{\tensor r_i}.
\end{equation}
Assumption (1) in this theorem implies, by Theorem~\ref{thm4thmomentcont}, that $f_k\cont{p} f_k^\ast\to0$ in $L^2$ for each $p\in\{
1,\ldots,n-1\}$. Therefore, from Proposition~\ref{propvanishingpairingintegrals}, it follows that for each of the decomposed
connected pairings $\pi_i$ with $r_i\ge3$, the corresponding pairing
integral $\int_{\pi_i} f_k^{\tensor r_i}$ converges to $0$ in $L^2$.
Since the number of factors $m$ in the product is bounded above by $r$
(which does not grow with $k$), this demonstrates that equation (\ref{eqmainproof1}) really expresses the limiting $r$th moment as a sum over
a small subset of $NC_2(n^{\tensor r})$. Let $NC_2^2(n^{\tensor r})$
denote the set of those respectful pairings $\pi$ such that, in the
decomposition $\pi=\pi_1\sqcup\cdots\sqcup\pi_m$, each $r_i=2$, in
other words, such that the connected components of the graph $C_\pi$
each have two vertices. Thus we have shown that\looseness=-1
%
\begin{equation} \label{eqmainproof3} \lim_{k\to\infty}\ff
[I^S_n(f_k)^r ] = \sum_{\pi\in NC_2^2(n^{\tensor r})} \lim
_{k\to
\infty}\int_\pi f_k^{\tensor r}.
\end{equation}\looseness=0
Note: if each $r_i=2$ and $r=r_1+\cdots+r_m$, then $r=2m$ is even. In
other words, if $r$ is odd, then $NC_2^2(n^{\tensor r})$ is empty, and
we have proved that all limiting odd moments of $I^S_n(f_k)$ are $0$.
If $r=2m$ is even, on the other hand, then the factors $\pi_i$ in the
decomposition of $\pi$ can each be thought of as $\pi_i\in
NC_2(n\tensor n)$. The reader may readily check that the only
noncrossing pairing that respects $n\tensor n$ is the totally nested
pairing $\pi_i = \{\{n,n+1\},\{n-1,n+2\},\ldots,\{1,2n\}\}$ in Figure
\ref{figpairingsof[6]}. Thus, utilizing the mirror symmetry of~$f_k$,
\[
\int_{\pi_i} f_k\tensor f_k = \int_{\pi_i} f_k\tensor f_k^\ast= \|
f_k\|
_{L^2(\R_+^n)}^2 = 1.
\]
Therefore, equation (\ref{eqmainproof3}) reads
%
\begin{equation} \label{eqmainproof4} \lim_{k\to\infty}\ff
[I^S_n(f_k)^{2m} ] = \sum_{\pi\in NC_2^2(n^{\tensor 2m})} 1 = |NC_2^2(n^{\tensor2m}) |.
\end{equation}
In each tensor factor of $n^{\tensor2m}$, all edges of each pairing in
$\pi$ act as one unit (since they pair in a uniform nested fashion as
described above); this sets up a bijection $NC_2^2(n^{\tensor2m})
\cong NC_2(2m)$. The set of noncrossing pairings of $[2m]$ is
well\vadjust{\goodbreak}
known to be enumerated by the Catalan number $C_m$ (cf.~\cite{NicaSpeicherBook}, Lemma~8.9), which is the $2m$th moment of $S(0,1)$; see
the discussion following equation (\ref{eqsemicircle}). This completes
the proof.
\end{pf}

Next we prove the Wigner--Wiener transfer principle, Theorem~\ref{thmtransfer}, restated below.

\begin{theoremtransfer*}
Let $n\ge2$ be an integer, and let $(f_k)_{k\in\N}$ be a sequence of
fully symmetric functions in $L^2(\R_+^n)$. Let $\sigma>0$ be a finite
constant. Then, as $k\to\infty$:
\begin{longlist}[(2)]
\item[(1)] $\E [I^W_n(f_k)^2 ] \to n!\sigma^2$ if and
only if
$\ff [I^S_n(f_k)^2 ]\to\sigma^2$;
\item[(2)] if the asymptotic relations in (1) are verified, then
$I^W_n(f_k)$ converges in law to a normal random variable $N(0,n!\sigma
^2)$ if and only if $I^S_n(f_k)$ converges in law to a semicircular
random variable $S(0,\sigma^2)$.
\end{longlist}
\end{theoremtransfer*}

\begin{pf} Point (1) is a simple consequence of the Wigner isometry
of equation (\ref{eqWignerisometry}), stating that for fully symmetric
$f\in L^2(\R_+^n)$, $\ff [I^S_n(f)^2 ] = \|f\|_2^2$ (since $f$
is fully symmetric, $f=f^\ast$ in particular), together with the
classical Wiener isometry which states that $\E [I^W_n(f)^2 ]
= n!\|f\|_2^2$. For point~(2), by renormalizing $f_k$ we may apply
Theorems~\ref{thm4thmomentsemicircle} and~\ref{thm4thmomentcont}
to see that $I_n^S(f_k)$ converges\vspace*{-1.5pt} to $S(0,1)$ in law if and only if
the contractions $f_k\cont{p}f_k^\ast= f_k\cont{p}f_k$ converge\vspace*{-2pt} to $0$
in $L^2$ for $p=1,2,\ldots, n-1$. Since $f$ is fully symmetric, these
nested contractions $f_k\cont{p}f_k$ are the same as the contractions
$f\tensor_p f$ in \cite{nunugio} (cf. Remark~\ref{rkcontractionsvsnested}),
and the main theorems in that paper show that these
contractions tend to $0$ in $L^2$ if and only if the Wiener integrals
$I^W_n(f_k)$ converge in law to a normal random variable, with variance
$n!$ due to our normalization. This completes the proof.
\end{pf}

As an application, we prove a free analog of the Breuer--Major theorem
for stationary vectors. This classical theorem can be stated as follows.

\begin{theorem*}[(Breuer--Major theorem)] Let $(X_k)_{k\in\Z}$ be a
doubly-infinite sequence of (jointly Gaussian) standard normal random
variables, and let $\rho(k) = \E(X_0X_k)$ denote the covariance
function. Suppose there is an integer $n\ge1$ such that $\sum_{k\in
\Z}
|\rho(k)|^n <\infty$. Let $H_n$ denote the $n$th Hermite polynomial,
\[
H_n(x) = (-1)^n e^{x^2/2}\,\frac{d^n}{dx^n}e^{-x^2/2}.
\]
($\{H_n \dvtx n\ge0\}$ are the monic orthogonal polynomials associated
to the law $N(0,1)$.) Then the sequence
\[
V_m = \frac{1}{\sqrt{m}}\sum_{k=0}^{m-1} H_n(X_k) \mathop
{\longrightarrow}^{\mathrm{law}}
 N(0,n!\sigma^2) \qquad \mbox{as }
m\to
\infty,
\]
where $\sigma^2 = \sum_{k\in\Z}\rho(k)^n$.
\end{theorem*}

See, for example, the preprint \cite{NPP} for extensions and
quantitative improvements of this theorem. Note that the Hermite
polynomial $H_n$ is related to Wiener integrals as follows: if
$(W_t)_{t\ge0}$ is a standard Brownian motion, then $W_1$ is a $N(0,1)$
variable, and
\[
H_n(W_1) = I^W_n\bigl(\1_{[0,1]^n}\bigr).
\]
(See, e.g., \cite{Kuo}.) The function $\1_{[0,1]^n}$ is fully
symmetric. On the other hand, if $(S_t)_{t\ge0}$ is a free Brownian
motion, then
\[
I^S_n\bigl(\1_{[0,1]^n}\bigr) = U_n(S_1),
\]
where $U_n$ is the $n$th Chebyshev polynomial of the second kind,
defined (on $[-2,2]$) by
%
\begin{equation} \label{eqTchebyshevdef}
U_n(2\cos\theta) = \frac{\sin((n+1)\theta)}{\sin\theta}.
\end{equation}
($\{U_n\dvtx n\ge0\}$ are the monic orthogonal polynomials associated to
the law $S(0,1)$; see \cite{BianeSpeicher,VDM}.) Hence, the
Wigner--Wiener transfer principle Theorem~\ref{thmtransfer}
immediately yields the following \textit{free Breuer--Major theorem}.

\begin{corollary} \label{corBreuerMajor} Let $(X_k)_{k\in\Z}$ be a
doubly-infinite semicircular system random variables $S(0,1)$, and let
$\rho(k)=\ff(X_0X_k)$ denote the covariance function with $X_0$.
Suppose there is an integer $n\ge1$ such that $\sum_{k\in\Z} |\rho
(k)|^n <\infty$. Then the sequence
\[
V_m = \frac{1}{\sqrt{m}}\sum_{k=0}^{m-1} U_n(X_k) \mathop
{\longrightarrow}^{\mathrm{law}} S(0,\sigma^2) \qquad \mbox{as } m\to
\infty,
\]
where $\sigma^2 = \sum_{k\in\Z}\rho(k)^n$.
\end{corollary}

\section{Free stochastic calculus} \label{sectionMalliavin}

In this section, we briefly outline the definitions and properties of
the main players in the free Malliavin calculus. We closely follow
\cite{BianeSpeicher}. The ideas that led to the development of stochastic
analysis in this context can be traced back to \cite{KumSpeicher};
\cite{BianeSpeicher2} provides an important application to the theory of
free entropy.

\subsection{Abstract Wigner space} \label{sectionabstractWignerspace} As in Nualart's treatise \cite{nualartbook}, we first set up the
constructs of the Malliavin calculus in an abstract setting, then
specialize to the case of stochastic integrals. As discussed in Section
\ref{sectionfBM}, the free Brownian motion is canonically constructed
on the free Fock space $\mathscr{F}_0(\H)$ over a separable Hilbert
space $\H$. Refer to the algebra $\mathcal{S}(\H)$ [generated by the
field variables $X(h)$ for $h\in\H$], endowed with the vacuum
expectation state $\ff$, as an \textit{abstract Wigner space}. While
$\mathcal{S}(\H)$ consists of operators \textit{on}~$\mathscr{F}_0(\H)$,
it can be identified as a subset of the Fock space due to the following fact.

\begin{proposition} \label{propFockiso} The function
%
\begin{eqnarray}\label{eqFockspaceisom}
 \mathcal{S}(\H)&\to&\mathscr{F}_0(\H) ,\nonumber
\\[-8pt]
\\[-8pt]
Y&\mapsto& Y\Omega
\nonumber
\end{eqnarray}
is an injective isometry. It extends to an isometric isomorphism from
the noncommutative $L^2$-space $L^2(\mathcal{S}(\H),\ff)$ onto
$\mathscr
{F}_0(\H)$.
\end{proposition}

In fact, the action of the map in equation (\ref{eqFockspaceisom}) can
be explicitly written in terms of Chebyshev polynomials [introduced in
equation (\ref{eqTchebyshevdef})]. If $\{h_i\}_{i\in\N}$ is an
orthonormal basis for $\H$, $k_1,k_2,\ldots,k_r$ are indices with
$k_j\ne k_{j+1}$ for $1\le j<r$, and $n_1,\ldots,n_r$ are positive
integers, then
%
\begin{equation} \label{eqTchebyshev}\quad
U_{n_1}(X(h_{k_1}))\cdots
U_{n_r}(X(h_{k_r}))\Omega= h_{k_1}^{\tensor n_1}\tensor\cdots\tensor
h_{k_r}^{\tensor n_r}\in\mathscr{F}_0(\H).
\end{equation}
[This is the precise analogue of the classical theorem with $X( \cdot
)$ an isonormal Gaussian process and the $U_n$ replaced by Hermite
polynomials $H_n$; in the classical case the tensor products are all
symmetric, hence the disjoint neighbors condition on the indices
$k_1,\ldots,k_r$ is unnecessary.] Hence, in order to define a gradient
operator (an analogue of the Cameron--Gross--Malliavin derivative) on the
abstract Wigner space $\mathcal{S}(\H)$, we may begin by defining it on
the Fock space $\mathscr{F}_0(\H)$.

\subsection{Derivations, the gradient operator, and the divergence
operator} \label{sectiongraddiv}

In free probability, the notion of a derivative is replaced by a
\textit{free difference quotient}, which generalizes the following
construction. Let $u\dvtx \R\to\C$ be a~$C^1$ function. Then define a
function $\del u\dvtx \R\times\R\to\C$ by
%
\begin{equation} \label{eqfreediffq} \del u (x,y) =
\cases{\displaystyle
\frac{u(x)-u(y)}{x-y}, &\quad $x\ne y $,\cr\displaystyle
 u'(x), &\quad $x=y$.}
\end{equation}
The function $\del u$ is continuous on $\R^2$ since $u$ is $C^1$. This
operation is a~\textit{derivation} in the following sense (as the reader
may readily verify): if $u,v\in C^1(\R)$ then
%
\begin{equation} \label{eqderivation} \del(uv)(x,y) = u(x)\del v(x,y)
+ \del u(x,y) v(y).
\end{equation}
Hence, $\del u \in L^2_{\mathrm{loc}}(\R^2)\cong L^2_{\mathrm{loc}}(\R)\tensor
L^2_{\mathrm{loc}}(\R)$. In other words, we can think of $\del$ as a map
%
\begin{equation} \label{eqderivationtensorrange} \del\dvtx C^1(\R
)\to L^2_{\mathrm{loc}}(\R)\tensor L^2_{\mathrm{loc}}(\R).
\end{equation}
If we restrict $\del$ to polynomials $u\in\C[X]$ in a single
indeterminate, then $\del u\in\C[X,Y]$, polynomials in two (commuting)
variables, and the same isomorphism yields $\C[X,Y]=\C[X]\tensor\C[X]$.
The action of $\del$ can be succinctly expressed here as
%
\begin{eqnarray} \label{eqdelmonomial}
 \del\dvtx \C[X]&\to&\C[X]\tensor\C[X],\nonumber
 \\[-8pt]
 \\[-8pt]
X^n&\mapsto&\sum_{j=1}^n X^{j-1}\tensor X^{n-j}.
\nonumber
\end{eqnarray}
 The operator $\del$ is called the \textit{canonical derivation}.
In the context of equation~(\ref{eqdelmonomial}), the derivation
property is properly expressed as follows:
%
\begin{equation} \label{eqderivation2} \del(AB) = (A\tensor1)\cdot
\del B + \del A\cdot(1\tensor B).
\end{equation}
 It is not hard to check that $\del$ is, up to scale, the
unique such derivation which maps $\C[X]$ into $\C[X]\tensor\C[X]$
(i.e., the only derivations on $\R$ are multiples of the usual
derivative). This uniqueness fails, of course, in higher dimensions.

Free difference quotients are noncommutative multivariate
generalizations of this operator $\del$ (acting, in particular, on
noncommutative polynomials). The definition follows.
\begin{definition}\label{deffreedifferencequotient} Let
$\mathscr
{A}$ be a unital von Neumann algebra, and let $X\in\mathscr{A}$. The
\textit{free difference quotient} $\del_X$ in the direction $X$ is the
unique derivation [cf. equation (\ref{eqderivation2})] with the
property that $\del_X (X) = 1\tensor1$.
\end{definition}
(There is a more general notion of free difference quotients
relative to a~subalgebra, but we will not need it in the present
paper.) Free difference quotients are central to the analysis of free
entropy and free Fisher information (cf. \cite{Voiculescu4,Voiculescu3}). The operator $\del$ plays the role of the
derivative in the version of It\^o's formula that holds for the
stochastic integrals discussed below in Section~\ref{sectionfreestochasticintegration}; cf. \cite{BianeSpeicher}, Proposition~4.3.2.
We will use $\del$ and $\del_X$, and their associated calculus (cf.
\cite{Voiculescu4}), in the calculations in Section~\ref{sectionMalliavinestimates}. We mention them here to point out a
counter-intuitive property of derivations in free probability: their
range is a tensor-product space.

Returning to abstract Wigner space, we now proceed to define a free
analog of the Cameron--Gross--Malliavin derivative in this context; it
will be modeled on the behavior (and hence tensor-product range space)
of the derivation $\del$.

\begin{definition} \label{defgradient} The \textit{gradient operator}
$\nabla\dvtx \mathscr{F}_0(\H)\to\mathscr{F}_0(\H)\tensor\H
\tensor
\mathscr{F}_0(\H)$ is densely defined as follows: $\nabla\Omega= 0$,
and for vectors $h_1,\ldots,h_n\in\H$,
%
\begin{equation} \label{eqgradient}
 \qquad \nabla(h_1\tensor\cdots\tensor h_n) = \sum_{j=1}^n (h_1\tensor
\cdots
\tensor h_{j-1})\tensor h_j \tensor(h_{j+1}\tensor\cdots\tensor h_n),
\end{equation}
where $h_1\tensor\cdots\tensor h_{j-1}\equiv\Omega$ when $j=1$ and
$h_{j+1}\tensor\cdots\tensor h_n \equiv\Omega$ when $j=n$. In
particular, $\nabla h = \Omega\tensor h \tensor\Omega$.

The \textit{divergence operator} $\delta\dvtx \mathscr{F}_0(\H)\tensor
\H
\tensor\mathscr{F}_0(\H)\to\mathscr{F}_0(\H)$ is densely defined as
follows: if $h_1,\ldots,h_n$ and $g_1,\ldots,g_m$ and $h$ are in $\H
$, then
%
\begin{eqnarray} \label{eqdivergence}
&&\delta \bigl( (h_1\tensor
\cdots
\tensor h_n)\tensor h \tensor(g_1\tensor\cdots\tensor g_m) \bigr)\nonumber
\\[-8pt]
\\[-8pt]
&& \qquad =
h_1\tensor\cdots\tensor h_n\tensor h\tensor g_1\tensor\cdots\tensor
g_m.
\nonumber
\end{eqnarray}
\end{definition}

 These actions, on first glance, look trivial; the important
point is the range of $\nabla$ and the domain of $\delta$ are tensor
products, and so the placement of the parentheses in equations (\ref{eqgradient}) and (\ref{eqdivergence}) is very important. When we
reinterpret $\nabla,\delta$ in terms of their action on stochastic
integrals, they will seem more natural and familiar.

The operator $N_0\equiv\delta\nabla\dvtx \mathscr{F}_0(\H)\to
\mathscr
{F}_0(\H)$ is the \textit{free Ornstein--Uhlenbeck operator} or \textit{free
number operator}; cf. \cite{Biane1}. Its action on an $n$-tensor is
given by $N_0 (h_1\tensor\cdots\tensor h_n) = n h_1\tensor\cdots
\tensor h_n$. In particular, the free Ornstein--Uhlenbeck operator,
densely defined on its natural domain, is invertible on the orthogonal
complement of the vacuum vector. This will be important in Section~\ref{lastsection}. It is easy to describe the domains $\mathscr{D}(N_0)$
and $\mathscr{D}(N_0^{-1})$; we will delay these descriptions until
Section~\ref{sectiondivgrad2}.

Definition~\ref{defgradient} defines $\nabla,\delta$ on domains
involving the algebraic Fock\break space~$\mathscr{F}_{\mathrm{alg}}(\H)$
(consisting of finitely-terminating sums of tensor products of vectors
in $\H$). It is then straightforward to show that they are closable
operators, adjoint to each other. The preimage of $\mathscr
{F}_{\mathrm
{alg}}(\H)$ under the isomorphism of equation~(\ref{eqFockspaceisom})
is actually contained in $\mathcal{S}(\H)$: Equation (\ref{eqTchebyshev}) shows that it consists of noncommutative polynomials in
variables $\{X(h),h\in\H\}$. Denote this space as~$\mathcal
{S}_{\mathrm
{alg}}(\H)$. We will concern ourselves primarily with the actions
of~$\nabla,\delta$ on this polynomial algebra (as is typical in the
classical setting as well). Note, we actually identify $\mathcal
{S}_{\mathrm{alg}}(\H)$ as a subset of $\mathscr{F}(\H)$ via
Proposition~\ref{propFockiso}, therefore using the same symbols
$\nabla,\delta$ for the conjugated actions of these Fock space
operators. Under this isomorphism, the full domain $\mathscr{D}(\nabla
)$ is the closure of~$\mathcal{S}_{\mathrm{alg}}(\H)$; similarly,
$\mathscr{D}(N_0)$ and $\mathscr{D}(N_0^{-1})$ have~$\mathcal
{S}_{\mathrm{alg}}(\H)$ (minus constants in the latter case) as a core.

\begin{proposition} \label{propderivation} The gradient operator
$\nabla\dvtx \mathcal{S}_{\mathrm{alg}}(\H)\to\mathcal
{S}_{\mathrm
{alg}}(\H)\tensor\H\tensor\mathcal{S}_{\mathrm{alg}}(\H)$ is a
derivation.
%
\begin{equation} \label{eqgradientderivation} \nabla(AB) = A\cdot
(\nabla B) + (\nabla A)\cdot B, \qquad A,B\in\mathcal{S}_{\mathrm
{alg}}(\H).
\end{equation}
\end{proposition}

 In equation (\ref{eqgradientderivation}), the left and right
actions of $\mathcal{S}_{\mathrm{alg}}(\HH)$ are the obvious ones
$A\cdot(U\tensor h\tensor V) = (AU)\tensor h\tensor V$ and $(U\tensor
h\tensor V)\cdot B = U\tensor h\tensor(VB)$. This is the same
derivation property as in equation (\ref{eqderivation2}). In
particular, iterating equation (\ref{eqgradientderivation}) yields the formula
%
\begin{equation} \label{eqgradientderivation2}
\nabla(X(h_1)\cdots X(h_n) )\! = \!\sum_{j=1}^n X(h_1)\cdots
X(h_{j-1})\!\tensor\! h_j\!\tensor\! X(h_{j+1})\cdots X(h_n).\hspace*{-40pt}
\end{equation}
When $n=1$, equation (\ref{eqgradientderivation2}) says $\nabla X(h) =
1\tensor h\tensor1$, which matches the classical gradient operator (up
to the additional tensor product with $1$).

As shown in \cite{BianeSpeicher}, both operators $\nabla$ and $\delta$
are densely defined and closable operators, both with respect to the
$L^2(\ff)$ [or $L^2(\ff\tensor\ff)$] topology and the weak operator
topology. It is most convenient to work with them on the dense domains
given in terms of $\mathcal{S}_{\mathrm{alg}}$.

We now state the standard integration by parts formula. First, we need
an appropriate pairing between the range of $\nabla$ and $\H$, which is
given by the linear extension of the following.
%
\begin{eqnarray} \label{eqpairing}
 \langle \cdot , \cdot \rangle_\H\dvtx \bigl(\mathcal
{S}_{\mathrm
{alg}}(\H)\tensor\H\tensor\mathcal{S}_{\mathrm{alg}}(\H) \bigr)
\times\H&\to&\mathcal{S}_{\mathrm{alg}}(\H)\tensor\mathcal
{S}_{\mathrm
{alg}}(\H), \nonumber
\\[-8pt]
\\[-8pt]
 \langle A\tensor h_1\tensor B, h_2\rangle_\H&=& \langle h_1,h_2\rangle
 A\tensor B.
\nonumber
\end{eqnarray}
In the special case $\H=L^2(\R_+)$ to which we soon restrict, this
pairing is quite natural; see equation (\ref{eqpairing2}) below. The
next proposition appears as~\cite{BianeSpeicher}, Lemma~5.2.2.\vspace*{-2pt}

\begin{proposition}[(Biane, Speicher)] \label{propintegrationbyparts}
If $Y\in\mathcal{S}_{\mathrm{alg}}(\H)$ and $h\in\H$,
%
\begin{equation} \label{eqintegrationbyparts} \ff\tensor\ff
(\langle\nabla Y,h\rangle_\H ) = \ff \bigl(Y\cdot X(h) \bigr).\vspace*{-2pt}
\end{equation}
\end{proposition}

\begin{remark} Since $\langle\nabla Y,h\rangle_\H$ is in the tensor
product $\mathcal{S}_{\mathrm{alg}}(\H)\tensor\mathcal{S}_{\mathrm
{alg}}(\H)$, its expectation must be taken with respect to the product
measure $\ff\tensor\ff$.\vspace*{-2pt}
\end{remark}

\subsection{Free stochastic integration and biprocesses} \label{sectionfreestochasticintegration} We now specialize to the case $\H= L^2(\R
_+)$. In this setting, we have already studied well the field variables $X(h)$.
%
\begin{equation} \label{eqX=I} X(h) = I_1^S(h) = \int h(t)\,dS_t.
\end{equation}
[Equation (\ref{eqX=I}) follows easly from the construction $S_t =
X(\1
_{[0,t]})$ of free Brownian motion.] To improve readability, we refer
to the polynomial algebra $\mathcal{S}_{\mathrm{alg}}(L^2(\R_+))$
simply as $\mathcal{S}_{\mathrm{alg}}$; therefore, since $S_t = X(\1
_{[0,t]})$, $\mathcal{S}_{\mathrm{alg}}$ contains all (noncommutative)
polynomial functions of free Brownian motion. The gradient $\nabla$
maps $\mathcal{S}_{\mathrm{alg}}$ into $\mathcal{S}_{\mathrm
{alg}}\tensor L^2(\R_+) \tensor\mathcal{S}_{\mathrm{alg}}$. It is
convenient to identify the range space in the canonical way with
vector-valued $L^2$-functions.
\[
\mathcal{S}_{\mathrm{alg}}\tensor L^2(\R_+)\tensor\mathcal
{S}_{\mathrm{alg}}
\cong L^2 (\R_+; \mathcal{S}_{\mathrm{alg}}\tensor\mathcal
{S}_{\mathrm{alg}} ).
\]
That is, for $Y\in\mathcal{S}_{\mathrm{alg}}$, we may think of
$\nabla
Y$ as a function. As usual, for $t\ge0$, denote $(\nabla Y)(t) =
\nabla
_t Y$. Thus, $\nabla Y$ is a noncommutative stochastic process taking
values in the tensor product $\mathcal{S}_{\mathrm{alg}}\tensor
\mathcal
{S}_{\mathrm{alg}}$.\vspace*{-2pt}

\begin{definition} \label{defbiprocess} Let $(\mathscr{A},\ff)$ be
a $W^\ast$-probability space. A \textit{biprocess} is a~stochastic process
$t\mapsto U_t\in\mathscr{A}\tensor\mathscr{A}$. For $1\le p\le
\infty$,
say $U$ is an \textit{$L^p$ biprocess}, $U\in\mathscr{B}_p$, if the norm
%
\begin{equation} \label{eqBpnorm} \|U\|_{\mathscr{B}_p}^2 = \int
_0^\infty\|U_t\|_{L^p(\mathscr{A}\tensor\mathscr{A},\ff\tensor\ff
)}^2\,
dt\vadjust{\goodbreak}
\end{equation}
is finite. (When $p=\infty$ the inside norm is just the operator norm
of $U_t$ in $\mathscr{A}\tensor\mathscr{A}$.)

 Let $\{\mathscr{A}_t \dvtx t\ge0\}$ be a filtration of
subalgebras of $\mathscr{A}$; say that $U$ is \textit{adapted} if $U_t\in
\mathscr{A}_t\tensor\mathscr{A}_t$ for all $t\ge0$.

 A biprocess is called \textit{simple} if it is of the form
%
\begin{equation} \label{eqsimplebiprocess} U = \sum_{j=1}^n
A_j\tensor B_j \1_{[t_{j-1},t_j)},
\end{equation}
where $0=t_0<t_1<\cdots<t_n$ and $A_j,B_j$ are in the algebra
$\mathscr
{A}$. The simple biprocess in equation (\ref{eqsimplebiprocess}) is
adapted if and only if $A_j,B_j\in\mathscr{A}_{t_{j-1}}$ for $1\le
j\le
n$. The closure of the space of simple biprocesses in $\mathscr{B}_p$
is denoted~$\mathscr{B}_p^a$, the space of \textit{$L^p$ adapted biprocesses}.
\end{definition}

\begin{remark} Customarily, our algebra $\mathscr{A}$ will contain a
free Brownian motion $S=(S_t)_{t\ge0}$, and we will consider only
filtrations $\mathscr{A}_t$ such that $S_s\in\mathscr{A}_t$ for
$s\le
t$. Thus, when we say a process or biprocess is adapted, we typically
mean with respect to the free Brownian filtration.
\end{remark}

 So, if $Y\in\mathcal{S}_{\mathrm{alg}}$, then $\nabla Y$ is
a biprocess. Since $\mathcal{S}_{\mathrm{alg}}$ consists of polynomials
in free Brownian motion, it is not too hard to see that $\nabla Y\in
\mathscr{B}_p$ for any $p\ge1$ (cf. \cite{BianeSpeicher}, Proposition~5.2.3). Note that the pairing of equation (\ref{eqpairing}), in the case $\H=L^2(\R_+)$, amounts to the following. If
$U\in\mathscr{B}_2$ is an $L^2$ biprocess and $h\in L^2(\R_+)$, then
%
\begin{equation} \label{eqpairing2} \langle U, h\rangle_{L^2(\R
_+)} =
\int_{\R_+} U_t \overline{h(t)}\,dt.
\end{equation}

We now describe a generalization of the Wigner integral $\int h(t)\,
dS_t$ to allow ``random'' integrands; moreover, we will allow
integrands that are not only processes but \textit{biprocesses}. (If $X_t$
is a process, then $X_t\tensor1$ is a biprocess, so we develop the
theory only for biprocesses.)

\begin{definition}\label{defstochasticintegralofsimplebiprocess} Let $U = \sum_{j=1}^n A_j\tensor B_j \1_{[t_{j-1},t_j)}$ be
a simple biprocess, and let $S=(S_t)_{t\ge0}$ be a free Brownian
motion. The \textit{stochastic integral} of $U$ with respect to $S$ is
defined to be
%
\begin{equation} \label{eqstochasticintegral} \int U_t\sh \,dS_t =
\sum
_{j=1}^n A_j(S_{t_j}-S_{t_{j-1}})B_j.
\end{equation}
\end{definition}

\begin{remark} \label{rksharp} The $\sh$-sign is used to denote the
action of $U_t$ on both the left and the right of the Brownian
increment. In general, we use it to denote the action of $\mathscr
{A}\tensor\mathscr{A}$ on $\mathscr{A}$ by $(A\tensor B)\sh C = \mathit{ACB}$;
more generally, for any vector space $\mathscr{X}$, it denotes the
action of $\mathscr{A}\tensor\mathscr{A}$ on $\mathscr{A}\tensor
\mathscr
{X}\tensor\mathscr{A}$ by $(A\tensor B)\sh(C\tensor X\tensor D) =
(AC)\tensor X\tensor(DB)$. Since the second tensor factor of $\mathscr
{A}$ acts on the right rather than the left, it might be more accurate
to describe $\sh$ as an action of $\mathscr{A}\tensor\mathscr
{A}^{\mathrm{op}}$, where the \textit{opposite algebra} $\mathscr
{A}^{\mathrm{op}}$ is equal to $\mathscr{A}$ as a set but has the
reversed product.
\end{remark}

\begin{remark} Let $U$ be a simple biprocess as in equation (\ref{eqsimplebiprocess}). If~$A_j$ are constant multiples of the identity, and
$B_j=1$, then the stochastic integral in Definition~\ref{defstochasticintegralofsimplebiprocess} reduces to the Wigner integral, $\int
U_t\sh \,dS_t = I^S_1(h)$ where $h = \sum_{j=1}^n A_j \1
_{[t_{j-1},t_j)}$.
\end{remark}

Let $U$ be an adapted simple biprocess. A standard calculation,
utilizing the freeness of the increments of $(S_t)_{t\ge0}$, yields
the general \textit{Wigner--It\^o isometry},
%
\begin{equation} \label{eqWigner-Itoisometry}
 \biggl\|\int U_t\sh \,dS_t \biggr\|_{L^2(\mathscr{A},\ff)} = \| U
 \|_{\mathscr{B}_2}.
\end{equation}
This isometry therefore extends the definition of the stochastic
integral to all of $\mathscr{B}_2^a$ by a density argument (since
simple biprocesses are dense in $\mathscr{B}_2^a$).

\subsection{An It\^o formula} \label{sectionItoformula}

There is a rich theory of free stochastic differential equations based
on the stochastic integral of Definition~\ref{defstochasticintegralofsimplebiprocess} (cf. \cite{fSDE1,fSDE2,fSDE3}) which mirror
classical processes (like the Ornstein--Uhlenbeck process) in the free
world, and \cite{fSDE4} which uses free SDEs for an important
application to random matrix ensembles and operator algebras. The
stochastic calculus in this context is based on a free version of the
It\^o formula, \cite{BianeSpeicher}, Proposition~4.3.4. It involves the
derivation $\del$ in place of the first order term; in order to
describe the appropriate It\^o correction term, we need the following
definition.\looseness=-1

\begin{definition} \label{defLaplacian} Let $\mu$ be a probability
measure on $\R$ all of whose moments are finite. Define the operator
$\Delta_\mu\dvtx \C[X]\to\C[X]$ on polynomials as follows:
%
\begin{equation} \label{eqLaplacian} \Delta_\mu h(x) = 2\,\frac
{d}{dx}\int_{\R} \del h(x,y) \mu(dy).
\end{equation}
\end{definition}

 The It\^o formula in our context applies to It\^o processes
of the form $M_t = M_0+\int_0^t U_s\sh \,dS_s + \int_0^t K_s\,ds$. For
our purposes, it suffices to take $U_s = \1_{[0,t]}1\tensor1$ so that
the stochastic integral $\int U_s\sh\, dS_s$ is just the free Brownian
motion $S_t$, and so we state the formula only in in this special case.

\begin{proposition}[(Biane, Speicher)] \label{propItoformula} Let
$K=(K_t)_{t\ge0}$ be a self-adjoint adapted process. Let $M_0$ be self
adjoint in $L^2(\mathcal{S},\ff)$, and let $M=(M_t)_{t\ge0}$ be a
process of the form
%
\begin{equation} \label{eqItoprocess} M_t = M_0+ S_t + \int_0^t
K_s\,
ds.
\end{equation}
Let $h\in\C[X]$ be a polynomial, and let $\Delta_t$ denote the operator
$\Delta_t = \Delta_{\mu_{M_t}}$; cf. equation (\ref{eqLaplacian}). Then
%
\begin{equation} \label{eqItoformula} h(M_t) = h(M_0) + \int_0^t
\del
h(M_s)\sh \,dM_s + \frac{1}{2}\int_0^t \Delta_s h (M_s)\,ds.\vadjust{\goodbreak}
\end{equation}
\end{proposition}

\begin{remark} In equation (\ref{eqItoformula}), we are viewing the
function $\del h$ as living in $\C[X]\tensor\C[X]$ directly rather than
$\C[X,Y]$. In particular, if $h(x)=x^n$ then $\del h (X) = \sum_{k=1}^n
X^{k-1}\tensor X^{n-k}$.
\end{remark}

\begin{remark} Of course, given equation (\ref{eqItoprocess}) defining
$M_t$, the integral $\int_0^t \del h(M_s)\sh \,dM_s$ in equation (\ref{eqItoformula}) is shorthand for
\[
\int_0^t \del h(M_s)\sh \,dM_s = \int_0^t \del h(M_s)\sh \,dS_s + \int_0^t
\del h(M_s)\sh K_s\,ds,
\]
following standard conventions of stochastic calculus.
\end{remark}

We will use Proposition~\ref{propItoformula} in the calculations in
Section~\ref{sectionMalliavinestimates} below. It will be convenient
to extend the It\^o formula beyond polynomial functions $h$ for this
purpose. The canonical derivation $\del$ of equation (\ref{eqderivation}) makes sense for any $C^1$-function $h$; we restrict this
domain slightly as follows. Suppose that $h$ is the Fourier transform
of a complex measure $\nu$ on $\R$,
%
\begin{equation} \label{eqFouriertransform} h(x) = \widehat{\nu}(x)=
\int_{\R} e^{ix\xi} \nu(d\xi).
\end{equation}
By definition, a complex measure is finite, and so such functions $h$
are continuous and bounded, $h\in C_b(\R)$. In order to fit into the
It\^o framework, such functions must be $L^2$ in the appropriate sense.
In the context of equation~(\ref{eqFouriertransform}), the relevant
normalization is as follows.
\begin{definition}\label{defFourierclass} Let $h$ have a Fourier
expansion as in equation (\ref{eqFouriertransform}). Define a seminorm
$\mathscr{I}_2(h)$ on such functions $h$ by
%
\begin{equation} \label{eqFouriernorm} \mathscr{I}_2(h) = \int_{\R}
\xi^2 |\nu|(d\xi).
\end{equation}
Denote by $\mathcal{C}_2$ the set of functions $h$ with $\mathscr
{I}_2(h)<\infty$.
\end{definition}

\begin{remark} \label{rkseminorm} $\mathscr{I}_2$ is not a norm: if
$h=a\in\C$ is a constant function, then $h = \widehat{a\delta_0}$, and
$\mathscr{I}_2(h) = \int\xi^2 |a|\delta_0(d\xi) = 0$. It is easy to
check that $\mathscr{I}_2$ is a~seminorm (i.e., nonnegative and
satisfies the triangle inequality), and that its kernel consists
exactly of constant functions in $\mathcal{C}_2$. Indeed, the quotient
of $\mathcal{C}_2$ by constants is a Banach space in the descended
$\mathscr{I}_2$-norm.
\end{remark}

 Standard Fourier analysis shows that $\mathcal{C}_2\subset
C_b^2(\R)$ (bounded twice-contin\-uously-differentiable functions), where
$\mathscr{I}_2(h)$ is like a sup-norm on the second derivative\vadjust{\goodbreak} $h''$.
In particular, nonconstant polynomials are not in $\mathcal{C}_2$. For
our purposes, we are only concerned with applying polynomials to
bounded operators, meaning that we only care about their action on a
compact subset of $\R$. In fact, locally any $C^\infty$ function is in
$\mathcal{C}_2$.\vspace*{-2pt}

\begin{lemma} \label{lemmasmoothfunctionsinC2} Let $r>0$. Given
any $C^\infty$ function $h\dvtx \R\to\C$, there is a~function
$h_r\in
\mathcal{C}_2$ such that $h(x) = h_r(x)$ for $|x|\le r$.\vspace*{-2pt}
\end{lemma}

\begin{pf} Let $\psi_r$ be a $C_c^\infty$ function such that $\psi
_r(x) = 1$ for $|x|\le r$. Then~$\psi_r h$ is equal to $h$ on $[-r,r]$.
This function is $C_c^\infty$, and hence its inverse Fourier transform
$(\psi_r h)^\vee$ is in the Schwartz space of rapidly-decaying smooth
functions. Set $\nu_r(d\xi) = (\psi_r h)^{\vee}(\xi)\,d\xi$; then
$\nu
_r$ has finite absolute moments of all orders, and $h_r\equiv\widehat
{\nu_r}=\psi_r h$ is in $\mathcal{C}_2$ and is equal to $h$ on $[-r,r]$.\vspace*{-2pt}
\end{pf}

 In particular, polynomials are locally in the class $\mathcal
{C}_2$. Later we will need the following result which says that
resolvent functions are globally in $\mathcal{C}_2$.\vspace*{-2pt}
\begin{lemma} \label{lemmaresolvents} For any fixed $z$ in the upper
half-plane $\C_+$, the function $\rho_z(x) = \frac{1}{z-x}$ is in
$\mathcal{C}_2$.\vspace*{-2pt}
\end{lemma}

\begin{pf} The resolvent $\rho_z$ is the Fourier transform of the
measure $\nu_z(d\xi) = -ie^{-iz\xi}\1_{(-\infty,0]}(\xi)\,d\xi$; a
simple calculation shows that $\mathscr{I}_2(\rho_z) = 2(\Im z)^{-3}$
when $\Im z > 0$.\vspace*{-2pt}
\end{pf}

The next theorem is a technical approximation tool which will greatly
simplify some of the more intricate calculations in
Section~\ref{sectionMalliavinestimates}.\vspace*{-2pt}

\begin{theorem} \label{theorempolynomialsdense} Let $K$ be a compact
interval in $\R$. Denote by $\mathcal{C}_2^{K,P}$ the subset of
$\mathcal{C}_2$ consisting of those functions in $\mathcal{C}_2$ that
are equal to polynomials on $K$. If $h\in\mathcal{C}_2$, there is a
sequence $h_n\in\mathcal{C}_2^{K,P}$ such that:
\begin{longlist}[(2)]
\item[(1)] $\mathscr{I}_2(h_n)\to\mathscr{I}_2(h)$ as $n\to\infty$;
\item[(2)] if $\mu$ is any probability measure supported in $K$, then
$\int h_n\,d\mu\to\int h\,d\mu$.\vspace*{-2pt}
\end{longlist}
\end{theorem}

 In fact, our proof will actually construct such a sequence
$h_n$ that converges to $h$ pointwise as well, although this is not
necessary for our intended applications. The proof of
Theorem~\ref{theorempolynomialsdense} is quite technical, and
is delayed to
\hyperref[appendix]{Appendix}.

Since $\mathcal{C}_2\subset C^2(\R)$, the operator $\del$ makes perfect
sense on $\mathcal{C}_2$ (and has $L^2$-norm appropriately controlled);
we can then reinterpret the function $\del h \in C^1(\R^2)$ as an
element of $L^2_{\mathrm{loc}}(\R)\tensor L^2_{\mathrm{loc}}(\R)$ so it fits the notation
of the It\^o formula equation (\ref{eqItoformula}). It will be useful
to have a more tensor-explicit representation of the function $\del h$
for $h\in\mathcal{C}_2$ in the sequel. If $h = \widehat{\nu}$, then
%
\begin{eqnarray} \label{eqdel2} \del h(x,y) &=& \frac{h(x)-h(y)}{x-y} =
\int_0^1 h'\bigl(\alpha x + (1-\alpha)y\bigr)\,d\alpha\nonumber
\\[-8pt]
\\[-8pt]&=& \int_0^1 \int_{\R}
i\xi
e^{i\alpha\xi x} e^{i(1-\alpha)\xi y} \nu(d\xi)\,d\alpha.
\nonumber\vadjust{\goodbreak}
\end{eqnarray}
Under the standard tensor identification, we can rewrite equation (\ref{eqdel2}) as
%
\begin{equation} \label{eqdelonC2} \del h(Y) = \int_0^1 \int_{\R}
i\xi\bigl(e^{i\alpha\xi Y}\tensor e^{i(1-\alpha)\xi Y}\bigr) \nu(d\xi)\,
d\alpha
.
\end{equation}
As for the It\^o correction term in equation (\ref{eqItoformula}), two
applications of the Dominated Convergence Theorem show that the
operator $\Delta_\mu$ of Definition~\ref{defLaplacian} is well defined
on $h\in\mathcal{C}_2$ whenever $\mu$ is compactly-supported, and the
resulting function $\Delta_\mu h$ is continuous. As such, all the terms
in the It\^o formula equation (\ref{eqItoformula}) are well defined for
$h\in\mathcal{C}_2$, and standard approximations show the following.\vspace*{-2pt}
\begin{corollary} \label{corItoformula}
The It\^o formula of equation (\ref{eqItoformula}) holds for $h\in\mathcal{C}_2$.\vspace*{-2pt}
\end{corollary}

\begin{remark} \label{rkfunctionalcalculus} The evaluations of the
functions $h$, $\del h$, and $\Delta_t h$ on the noncommutative random
variables $M_0$ and $M_t$ are given sense through functional calculus;
this is possible (and routine) because $M_0$ and $M_t$ are self adjoint.\vspace*{-2pt}
\end{remark}

\subsection{Chaos expansion for biprocesses} \label{sectionchaosexpansionforbiprocesses}

Recall the multiple Wigner integrals $I^S_n$ as discussed in Section
\ref{sectionWignerintegral}. By de-emphasizing the explicit
dependence on $n$, $I^S$ can then act (linearly) on finite sums $\sum_n
f_n$ of functions $f_n\in L^2(\R_+^n)\cong L^2(\R_+)^{\tensor n}$; that
is, $I^S$ acts on the algebraic Fock space $\mathcal{F}_{\mathrm
{alg}}=\mathscr{F}_{\mathrm{alg}}(L^2(\R_+))$. Utilizing the Wigner
isometry, equation (\ref{eqWignerisometry}), this means $I^S$ extends
to a map defined on the Fock space,
%
\begin{equation} \label{eqchaosexpansion} I^S\dvtx \mathscr
{F}_0\to
L^2(\mathcal{S},\ff);
\end{equation}
here and in the sequel, $\mathscr{F}_0 = \mathscr{F}_0(L^2(\R_+))$ and
$\mathcal{S} = \mathcal{S}(L^2(\R_+))$. In fact, the map in equation (\ref{eqchaosexpansion}) is an isometric isomorphism; this is one way
to state the Wigner chaos decomposition. This extended map $I^S$ is the
inverse of the map $Y\mapsto Y\Omega$ of Proposition~\ref{propFockiso}.

For $n,m$ positive integers, define for $f\in L^2(\R_+^n)\tensor
L^2(\R
_+^m) \cong L^2(\R_+^{n+m})$ the \textit{Wigner bi-integral}
%
\begin{equation} \label{eqbi-integral}
 [I_n^S\tensor I_m^S ](f) = \int f(t_1,\ldots
,t_n;s_1,\ldots
,s_m)\,dS_{t_1}\cdots \,dS_{t_n}\tensor dS_{s_1}\cdots \,dS_{s_m}.\hspace*{-40pt}
\end{equation}
To be clear: if $f = g\tensor h$ with $g\in L^2(\R_+^n)$ and $h\in
L^2(\R_+^m),$ then $ [I_n^S\tensor I_m^S ](f) =
I^S_n(g)\tensor I_m^S(h)$; in general, $I_n^S\tensor I_m^S$ is the
$L^2$-closed linear extension of this action. Thus,
\[
I^S_n\tensor I^S_m\dvtx L^2(\R_+^n)\tensor L^2(\R_+^m)\to
L^2(\mathcal
{S}\tensor\mathcal{S},\ff\tensor\ff).
\]
The Wigner isometry [cf. equation (\ref{eqWignerisometry})] in this
context then says that if $f\in L^2(\R_+^n)\tensor L^2(\R_+^m)$ and
$g\in L^2(\R_+^{n'})\tensor L^2(\R_+^{m'})$, then
%
\begin{eqnarray} \label{eqWignerbisometry} &&\ff\tensor\ff
\bigl(
[I_{n'}^S\tensor I_{m'}^S ](g)^\ast [I_n^S\tensor
I_m^S
](f) \bigr)\nonumber
\\[-8pt]
\\[-8pt]&& \qquad =
\cases{\displaystyle \langle f,g \rangle_{L^2(\R_+^n)\tensor L^2(\R_+^m)} ,&\quad
 if $n=n'$ and $m=m'$,\vspace*{1pt}\cr\displaystyle
 0,&\quad otherwise.
}
\nonumber\vadjust{\goodbreak}
\end{eqnarray}
This ``bisometry'' allows us to put the $I^S_n\tensor I^S_m$ together
for different $n,m$ as in equation (\ref{eqchaosexpansion}), to yield
an isometric isomorphism
%
\begin{equation} \label{chaosexpansion2}
I^S\tensor I^S\dvtx \mathscr{F}_0\tensor\mathscr{F}_0\to
L^2(\mathcal
{S}\tensor\mathcal{S},\ff\tensor\ff).
\end{equation}
What's more, by taking these Hilbert spaces as the ranges of
vector-valued $L^2(\R_+)$-functions, and utilizing the isomorphism
$L^2(\R_+;\mathfrak{A}\tensor\mathfrak{B}) \cong\mathfrak
{A}\tensor\break
L^2(\R_+)\tensor\mathfrak{B}$ for given Hilbert spaces $\mathfrak
{A},\mathfrak{B}$, we have an isometric isomorphism
%
\begin{equation} \label{eqchaosexpansion3}
I^S\tensor I^S \dvtx L^2 (\R_+ ; \mathscr{F}_0\tensor
\mathscr
{F}_0 ) \to\mathscr{B}_2.
\end{equation}
 Here $\mathscr{B}_2$ denotes the $L^2$ biprocesses (cf.
Definition~\ref{defbiprocess}), in this case taking values in
$\mathcal
{S}\tensor\mathcal{S}$. If $f\in L^2 (\R_+ ; \mathscr
{F}_0\tensor
\mathscr{F}_0 )$, the bi-integral acts only on components:
$
[I^S\tensor I^S ] (f) (t) = [I^S\tensor I^S ](f_t)$.
Equation (\ref{eqchaosexpansion3}) [through the action defined in
equation (\ref{eqbi-integral})] is the \textit{Wigner chaos expansion for
$L^2$ biprocesses} in the Wigner space.

As in the classical case, adaptedness is easily understood in terms of
the chaos expansion. If $U\in\mathscr{B}_2$, it has a chaos expansion
$U = [I^S\tensor I^S](f)$ for some $f\in L^2(\R_+ ; \mathscr
{F}_0\tensor\mathscr{F}_0)$, which we may write as an orthogonal sum
\[
f\dvtx t\mapsto f_t = \sum_{n,m=0}^\infty f_t^{n,m},
\]
where $f_t^{n,m}\in L^2(\R_+^n)\tensor L^2(\R_+^m)$. Then $U$ is
adapted (in the sense of Definition~\ref{defbiprocess}) if and only if
for each $n,m$ and $t_1,\ldots,t_n,s_1,\ldots,s_m\ge0$, the kernels
$f^{n,m}_t(t_1,\ldots,t_n;s_1,\ldots,s_m)$ are adapted, meaning they
are $0$ whenever $\max\{t_1,\ldots,t_n,s_1,\ldots,s_m\}>t$. In this
case, the stochastic integral defined in equations (\ref{eqstochasticintegral}) and (\ref{eqWigner-Itoisometry}) can be succinctly expressed;
cf. \cite{BianeSpeicher}, Proposition~5.3.7. In particular, if
$f^{n,m}\in L^2(\R_+ ; L^2(\R_+^n)\tensor L^2(\R_+^m))$ is
adapted, then\looseness=-1
%
\begin{eqnarray} \label{eqSkorohod1}
&& \qquad \int[I^S\tensor I^S](f_t)\sh
\,dS_t\nonumber
\\[-8pt]
\\[-8pt] && \qquad \qquad = \int f_t^{n,m}(t_1,\ldots,t_n;s_1,\ldots,s_m)\,dS_{t_1}\cdots
\,dS_{t_n}\,dS_t\,dS_{s_1}\cdots \,dS_{s_m}.
\nonumber
\end{eqnarray}\looseness=0
This is consistent with the notation of equation (\ref{eqbi-integral});
informally, it says~that
\[
(dS_{t_1}\cdots \,dS_{t_n}\tensor dS_{s_1}\cdots \,dS_{s_m})\sh \,dS_t =
dS_{t_n}\cdots \,dS_{t_1}\,dS_t\,dS_{s_1}\cdots \,dS_{s_m}
\]
as one would expect.

\subsection{Gradient and divergence revisited} \label{sectiondivgrad2}

Both the gradient and the divergence have simple representations in
terms of the chaos expansions in
Section~\ref{sectionchaosexpansionforbiprocesses}.

\begin{proposition}[(Propositions 5.3.9 and 5.3.10 in \cite{BianeSpeicher})] \label{propgradientonchaos} The gradient operator
is densely-defined and closable in
\[
\nabla\dvtx L^2(\mathcal{S},\ff)\to\mathscr{B}_2.\vadjust{\goodbreak}
\]
Its domain $\mathscr{D}(\nabla)$, expressed in terms of the chaos
expansion for $L^2(\mathcal{S},\ff)$, is as follows. If $f = \sum_n
f_n\in\mathscr{F}_0$ with $f_n\in L^2(\R_+^n)$, and if $L^2(\mathcal
{S},\ff)\ni Y = I^S(f)$, then $Y\in\mathscr{D}(\nabla)$ if and only if
%
\begin{equation} \label{eqdomain1} \sum_{n=0}^\infty n\|f_n\|
_{L^2(\R
_+^n)}^2 < \infty.
\end{equation}
In this case, the quantity in equation (\ref{eqdomain1}) is equal to
the norm
\[
\int_{\R_+} \|\nabla_t Y\|_{L^2(\mathcal{S}\tensor\mathcal{S},\ff
\tensor
\ff)}^2\,dt = \sum_{n=0}^\infty n\|f_n\|_{L^2(\R_+^n)}^2.
\]
Moreover, the action of $\nabla$ on this domain is determined by
%
\begin{eqnarray} \label{eqgradientonintegrals}
&&  \nabla_t \biggl(\int f(t_1,\ldots,t_n)\,
dS_{t_1}\cdots
\,dS_{t_n} \biggr)\nonumber\hspace*{-35pt}
\\[-10pt]
\\[-10pt] && \qquad = \sum_{k=1}^n \int f(t_1,\ldots
,t_{k-1},t,t_{k+1},\ldots,t_n)\,dS_{t_1}\cdots \,dS_{t_{k-1}}\tensor
dS_{t_{k+1}}\cdots \,dS_{t_n}.
\nonumber\hspace*{-35pt}\vspace*{-2pt}
\end{eqnarray}
\end{proposition}

\begin{remark} \label{rkdomainN0} It is similarly straightforward to
write the domain of the free Ornstein--Uhlenbeck operator in terms of
Wigner chaos expansions. If $Y = I^S(f)$ where $f = \sum_n f_n \in
\mathscr{F}_0$, then $Y\in\mathscr{D}(N_0)$ iff $\sum_n n^2\|f_n\|
_{L^2(\R_+^n)}^2 <\infty$.\vspace*{-1pt} Likewise, $Y\in\mathscr{D}(N_0^{-1})$ iff
$f_0=0$ and $\sum_{n>0} n^{-2}\|f_n\|_{L^2(\R_+^n)}^2 <\infty$. In
particular, we see that
%
\begin{equation} \label{eqdomainsN0} \mathscr{D}(N_0) \subset
\mathscr{D}(\nabla), \qquad\mathscr{D}(\nabla)\ominus\mathrm
{image}
(I^S_0) \subset\mathscr{D}(N_0^{-1}).\vspace*{-2pt}
\end{equation}
\end{remark}

The divergence operator can also be simply described in terms of the
chaos. We could similarly describe its domain, but its action on
adapted processes is already well known, as in the classical case.\vspace*{-2pt}

\begin{proposition}[(Propositions 5.3.9 and 5.3.11 in \cite{BianeSpeicher})] \label{propdivergenceonchaos} The divergence
operator is densely defined and closable in
\[
\delta\dvtx \mathscr{B}_2\to L^2(\mathcal{S},\ff).
\]
Using the chaos expansion for biprocesses, the action of $\delta$ is
determined as follows. If $f\in L^2(\R_+ ; L^2(\R_+^n)\tensor
L^2(\R
_+^m))$, then
%
\begin{eqnarray} \label{eqdivergenceonintegrals}\qquad
&& \delta \biggl(\int f_t(t_1,\ldots,t_n;s_1,\ldots
,s_m)\,
dS_{t_1}\cdots \,dS_{t_m}\tensor dS_{s_1}\cdots \,dS_{s_m} \biggr) \nonumber
\\[-10pt]
\\[-10pt]
\qquad&& \qquad =\int f_t(t_1,\ldots,t_n;s_1,\ldots,t_m)\,dS_{t_1}\cdots
\,dS_{t_n}\,dS_t\,dS_{s_1}\cdots \,dS_{s_m}.
\nonumber
\end{eqnarray}
In particular, comparing with equation (\ref{eqSkorohod1}), if $U$ is
an adapted biprocess $U\in\mathscr{B}_2^a$, then $U\in\mathscr
{D}(\delta)$ and
\[
\delta(U) = \int U_t\sh \,dS_t.\vadjust{\goodbreak}
\]
\end{proposition}

\begin{remark} In light of the second part of Proposition~\ref{propdivergenceonchaos}, the divergence operator is also called the
\textit{free Skorohod integral}. To be more precise, as in the classical case,
there is a domain $\mathbb{L}^{1,2}$ in between $\mathscr{B}_2^a$ and
the natural domain $\mathscr{D}(\delta)$ on which $\delta$ is closable
and such that for $U\in\mathbb{L}^{1,2}$ the relation $\nabla
_t(\delta
(U)) = U_t + \delta_s(\nabla_t U_s)$ holds true. It is this restriction
of $\delta$ that is properly called the Skorohod integral.
\end{remark}

\begin{remark} \label{rkchainrule} Given a random variable $X\in
\mathscr{D}(\nabla)$, using the derivation properties of the operators
$\del_X$ (cf. Definition~\ref{deffreedifferencequotient}) and
$\nabla$, it is relatively easy to derive the following chain rule. If
$p\in\C[X]$ is a polynomial, then\looseness=1
%
\begin{equation} \label{eqchainrule}
\nabla p(X) = \del_Xp(X)\sh \nabla X.
\end{equation}\looseness=0
\end{remark}

We conclude this section with one final result. The space of $L^2$
adapted biprocesses $\mathscr{B}_2^a$ is a closed subspace of the
Hilbert space $\mathscr{B}_2$; cf. Definition~\ref{defbiprocess}.
Hence there is an orthogonal projection $\Gamma\dvtx \mathscr
{B}_2\to
\mathscr{B}_2^a$. The next result is a free version of the Clark--Ocone
formula. It can be found as \cite{BianeSpeicher}, Proposition~5.3.12.

\begin{proposition} \label{propClark-Ocone} If $X\in\mathscr
{D}(\nabla
)$, then
\[
X = \ff(X) + \delta(\Gamma\nabla X).\vspace*{6pt}
\]
\end{proposition}

\section{Quantitative bounds on the distance to the semicircular
distribution} \label{lastsection}

As described in the restricted form of Theorem~\ref{theoremMalliavinestimate} in Section~\ref{sectionintroduction}, we are primarily
concerned in this section with quantitative estimates for the following
distance function on probability distributions.

\begin{definition} \label{defdistance} Given two self-adjoint random
variables $X,Y$, define the distance
\[
d_{\mathcal{C}_2}(X,Y) = \sup\{|\ff[h(X)]-\ff[h(Y)]| \dvtx h\in
\mathcal
{C}_2, \mathscr{I}_2(h)\le1\};
\]
the class $\mathcal{C}_2$ and the seminorm $\mathscr{I}_2$ are
discussed in Definition~\ref{defFourierclass}.
\end{definition}


\begin{remark} \label{rkdistanceonlaws}
Note that we could write the definition of $d_{\mathscr{C}_2}(X,Y)$
equally well as
\[
\sup \biggl\{ \biggl|\int h\,d\mu_X - \int h\,d\mu_Y \biggr| \dvtx h\in
\mathcal
{C}_2, \mathscr{I}_2(h)\le1 \biggr\}.
\]
In this form, it is apparent that $d_{\mathcal{C}_2}(X,Y)$ only depends
on the laws $\mu_X$ and $\mu_Y$ of the random variables $X$ and $Y$. In
computing it, we are therefore free to make any simplifying assumption
about the correlations of $X$ and $Y$ that are convenient; for example,
we may assume that $X$ and $Y$ are freely independent.
\end{remark}

Lemma~\ref{lemmaresolvents} shows that resolvent functions
$\rho_z(x) = (z-x)^{-1}$ are in $\mathcal{C}_2$ for $z\in\C_+$, and in
fact that if $\Im z = 1$, then $\mathscr{I}_2(\rho_z) = 2$. Thus,
\[
d_{\mathcal{C}_2}(X,Y) \ge\frac12 \sup_{\Im z =1} |\ff
[(z-X)^{-1}] - \ff[(z-Y)^{-1}] | = \frac12 \sup_{\Im z = 1}
|G_{\mu
_X}(z)-G_{\mu_Y}(z)|,
\]
where $G_\mu(z) = \int_{\R} (z-x)^{-1} \mu(dx)$ is the Stieltjes
transform of the law $\mu$. It is a standard theorem that convergence
in law is equivalent to convergence of the Stieltjes transform on any
set with an accumulation point, and hence this latter distance metrizes
converge in law; so our stronger distance $d_{\mathcal{C}_2}$ also
metrizes convergence in law. The class $\mathcal{C}_2$ is somewhat
smaller than the space of Lipschitz functions, and so this metric is, a
priori, weaker than the Wasserstein distance (as expressed in
Kantorovich form; cf. \cite{BianeVoiculescu,Wassersteinref}). However,
as Lemma~\ref{lemmasmoothfunctionsinC2} shows, all smooth
functions are locally in $\mathcal{C}_2$; the relative strength of
$d_{\mathcal{C}_2}$ versus the Wasserstein metric is an interesting
question we leave to future investigation.

\subsection{\texorpdfstring{Proof of Theorem \protect\ref{theoremMalliavinestimate}}
{Proof of Theorem 1.10}} \label{sectionMalliavinestimates}

We begin by restating Theorem~\ref{theoremMalliavinestimate} in the
language and full generality of Section~\ref{sectionMalliavin}.

\begin{theoremMalliavin*} Let $S$ be a standard semicircular random
variable; cf. equation (\ref{eqsemicircle}). Let $F$ be self adjoint in
the domain of the gradient, $F\in\mathscr{D}(\nabla)\subset
L^2(\mathcal
{S},\ff)$, with $\ff(F)=0$. Then
%
\begin{equation} \label{eqWasserstein1} d_{\mathcal{C}_2}(F,S) \le
\frac{1}{2}\ff\tensor\ff\biggl ( \biggl|\int\nabla_s(N_0^{-1}F)\sh
(\nabla
_s F)^\ast\,ds - 1\tensor1 \biggr| \biggr).
\end{equation}
\end{theoremMalliavin*}

\begin{pf} The main idea is to connect the random variables $F$ and
$S$ through a free Brownian bridge, and control the differential along
the path using free Malliavin calculus; cf. Section~\ref{sectionMalliavin}. For $0\le t\le1$, define
%
\begin{equation} \label{eqFt} F_t = \sqrt{1-t} F + S_t,
\end{equation}
where $S_t$ is a free Brownian motion. In particular, $S_1$ has the
same law as the random variable $S$. Since $d_{\mathcal{C}_2}(F,S)$
depends only on the laws of~$F$ and~$S$ individually, for convenience
we will take $S_t$ freely independent from~$F$. Fix a function $h\in
\mathcal{C}_2$. In the proceeding calculations, it will be useful to
assume that $h$ is a polynomial; however, polynomials are not in
$\mathcal{C}_2$. Rather, fix a compact interval $K$ in $\R$ that
contains the spectrum of $F_t$ for each $t\in[0,1]$; for example, since
$\|F_t\|\le2\sqrt{t}+\sqrt{1-t}\|F\|$, we could choose $K = [-2-\|F\|
,2+\|F\|]$. For the time being, we will assume that $h$ is equal to a
polynomial on $K$; that is, we take $h\in\mathcal{C}_2^{K,P};$ cf.
Theorem~\ref{theorempolynomialsdense}.

Define $g(t) = \ff[h(F_t)]$. The fundamental theorem of calculus yields
the desired quantity,
%
\begin{equation} \label{eqFTC}
  \ff[h(S)] - \ff[h(F)] = \ff
[h(F_1)]-\ff
[h(F_0)] = g(1)-g(0) = \int_0^1 g'(t)\,dt.\hspace*{-35pt}
\end{equation}
 We can use the free It\^o formula of equation (\ref{eqItoformula}) to calculate
 the derivative~$g'(t)$. In particular,\vadjust{\goodbreak} $dF_t =
-\frac{1}{2\sqrt{1-t}}F\,dt + dS_t$, and so applying equation~(\ref{eqItoformula}) yields
%
\begin{eqnarray} \label{eqIto1}
 d[h(F_t)] &=& \del h(F_t)\sh \,dF_t + \frac{1}{2}\Delta_t
h(F_t)\,dt \nonumber
\\[-9pt]
\\[-9pt]
&=& \del h(F_t)\sh \biggl\{-\frac{1}{2\sqrt{1-t}}F\,dt + dS_t \biggr\}
+\frac{1}{2}\Delta_t h(F_t)\,dt.
\nonumber
\end{eqnarray}
 Linearity (and uniform boundedness of all terms) allows us to
exchange $\ff$ with stochastic integrals; in particular, we may write
$dg(t) = \ff (d[h(F_t)] )$. The (stochastic integral of the)
term $\del h(F_t)\sh \,dS_t$ has mean $0$, and so we are left with two terms,
%
\begin{equation} \label{eqIto2} dg(t) = \frac{1}{2} \biggl\{-\frac
{1}{\sqrt{1-t}}\ff[\del h(F_t)\sh F] + \ff[\Delta_t h(F_t)] \biggr\}\,dt.
\end{equation}
The following lemma allows us to simplify these terms.

\begin{lemma} \label{lemmacalculations} Let $X$ and $Y$ be
self-adjoint random variables.
Let $h\in\mathscr{C}_2$.
\begin{longlist}[(b)]
\item[(a)] $\ff[\del h(Y)\sh X] = \ff[h'(Y)X]$.
\item[(b)] $\ff[\Delta_{\mu_Y} h(Y)] = \ff\tensor\ff [\del
h'(Y) ]$.
\end{longlist}
\end{lemma}

\begin{pf*}{Proof of Lemma~\ref{lemmacalculations}} By assumption
$h$ takes the form $h=\widehat{\nu}$ for some complex measure $\nu$
with finite second absolute moment.
\begin{longlist}[(b)]
\item[(a)] We use the representation of equation (\ref{eqdelonC2}) for
$\del$, so that
%
\begin{eqnarray}\label{eqlemma1}
 \del h(Y)\sh X &=& \int_0^1 d\alpha\int_\R i\xi
\nu
(d\xi) \bigl(e^{i\alpha\xi Y}\tensor e^{i(1-\alpha)\xi Y}\bigr)\sh X \nonumber
\\[-9pt]
\\[-9pt] &=&
\int
_0^1 d\alpha\int_\R i\xi \nu(d\xi) e^{i\alpha\xi
Y}Xe^{i(1-\alpha
)\xi Y}.
\nonumber
\end{eqnarray}
Since $\ff$ is a trace, $\ff[e^{i\alpha\xi Y}Xe^{i(1-\alpha)\xi Y}] =
\ff[e^{i\xi Y}X]$. Taking $\ff$ of both sides of equation (\ref{eqlemma1}), the $\alpha$ integration just yields a constant $1$, and so
%
\begin{equation} \label{eqlemma2} \ff[\del h(Y)\sh X] = \int_\R
i\xi
\nu(d\xi) \ff[e^{i\xi Y}X]
= \ff \biggl[\biggl (\int_\R i\xi e^{i\xi Y} \nu(d\xi)
\biggr)X \biggr].
\end{equation}
Since $h'(x) = \int_\R i\xi e^{i\xi x} \nu(d\xi)$, this yields the result.

\item[(b)] By Definition~\ref{defLaplacian}, $\Delta_{\mu_Y} h(x) =
2\frac{d}{dx}\int_\R\del h(x,y) \mu_Y(dy)$. Using the chain rule, we
can express $\del h(x,y) = \int_0^1 h'(\alpha x + (1-\alpha)y)\,
d\alpha
$. Since $h\in C^2$ and the integrand is bounded, we can rewrite
$\Delta
_{\mu_Y}h(x)$ as
%
\begin{eqnarray}
 \Delta_{\mu_Y} h(x) &=& 2\frac{d}{dx} \int_\R\mu
_Y(dy)\int_0^1d\alpha h'\bigl(\alpha x + (1-\alpha)y\bigr) \nonumber
\\[-8pt]
\\[-10pt]
&=& \int_\R\mu_Y(dy)\int_0^1 2\alpha\,d\alpha h''\bigl(\alpha x +
(1-\alpha
)y\bigr).
\nonumber\vadjust{\goodbreak}
\end{eqnarray}
Now $h''(x) = \int_\R-\xi^2e^{i\xi x} \nu(d\xi)$, and so
%
\begin{eqnarray}
 \Delta_{\mu_Y} h(x) &=& -\int_\R\xi^2 \nu(d\xi
)\int
_0^1 2\alpha\,d\alpha\int_\R\mu_Y(dy) e^{i(1-\alpha)\xi y}
e^{i\alpha
\xi x}\nonumber
\\[-8pt]
\\[-8pt]
&=& -\int_\R\xi^2 \nu(d\xi)\int_0^1 2\alpha\,d\alpha
e^{i\alpha\xi
x} \ff\bigl[e^{i(1-\alpha)\xi Y}\bigr].
\nonumber
\end{eqnarray}
Evaluating at $x=Y$ and taking the trace, this yields
%
\begin{equation} \label{eqlemmab1} \ff[\Delta_{\mu_Y}h(Y)] =
-\int_\R
\xi^2\nu(d\xi)\int_0^1 2\alpha\,d\alpha \ff[e^{i\alpha\xi
Y}]\ff
\bigl[e^{i(1-\alpha)\xi Y}\bigr].
\end{equation}
On the other hand, following the same identification as in equation (\ref{eqdelonC2}), we have
%
\begin{equation} \label{eqlemmab2} \del h'(Y) = -\int_0^1 d\alpha
\int
_\R\xi^2 \nu(d\xi) e^{i\alpha\xi Y}\tensor e^{i(1-\alpha)\xi Y}.
\end{equation}
Taking the trace yields
%
\begin{equation} \label{eqlemmab3} \ff\tensor\ff[\del h'(Y)] =
-\int
_0^1 d\alpha\int_\R\xi^2 \nu(d\xi) \ff[e^{i\alpha\xi Y}]\ff
\bigl[e^{i(1-\alpha)\xi Y}\bigr].
\end{equation}
Subtracting equation (\ref{eqlemmab2}) from equation (\ref{eqlemmab3})
and using Fubini's theorem (justified since the modulus of the
integrand is $\le\xi^2$ which is in $L^1(\nu\times[0,1])$) yields
%
\begin{eqnarray} \label{eqlemmab4} &&\ff\tensor\ff[\del h'(Y)] -
\ff
[\Delta_{\mu_Y}h(Y)] \nonumber
\\[-8pt]
\\[-8pt]&& \qquad = \int_0^1 (2\alpha-1)\,d\alpha\int_\R\xi
^2 \nu
(d\xi) \ff[e^{i\alpha\xi Y}]\ff\bigl[e^{i(1-\alpha)\xi Y}\bigr].
\nonumber
\end{eqnarray}
Equation (\ref{eqlemmab4}) expresses the difference $\ff\tensor\ff
[\del
h'(Y)] - \ff[\Delta_{\mu_Y}h(Y)]$ as an integral of the form $\int_0^1
(2\alpha-1) \kappa(\alpha)\,d\alpha$, where $\kappa$ is a function with
the symmetry $\kappa(\alpha)=\kappa(1-\alpha)$. The substitution
$\alpha
\mapsto1-\alpha$ shows that any such integral is $0$, which yields
the result.
\end{longlist}
\upqed
\end{pf*}

We now apply Lemma~\ref{lemmacalculations} to equation (\ref{eqIto2})
with $X=F$ and $Y=F_t$; note that $\Delta_t h(F_t)$ is by definition
(cf. Proposition~\ref{propItoformula}) equal to $\Delta_{\mu_{F_t}}
h(F_t)$. Equation (\ref{eqIto2}) then becomes
%
\begin{equation} \label{eqIto3}
g'(t) = \frac12\biggl \{-\frac{1}{\sqrt{1-t}}\ff[h'(F_t)F] + \ff
\tensor\ff
[\del h'(F_t)] \biggr\}.
\end{equation}
At this point, we invoke the free Malliavin calculus of variations
(cf. Section~\ref{sectionMalliavin}) to re-express these two
terms.
For the first term, we use a standard trick to introduce conditional
expectation; by Definition~\ref{defconditionalexpectation}, $\ff
[h'(F_t)F] = \ff [F\cdot\ff[h'(F_t)|F] ]$. Since $F\in
\mathscr
{D}(\nabla)$ and $\ff(F)=0$, equation (\ref{eqdomainsN0}) shows that
$F\in\mathscr{D}(N_0^{-1})$, and so $F = \delta(\nabla N_0^{-1}F)$. Hence
%
\begin{equation} \label{eqIto4}\qquad
\ff[h'(F_t)F] = \ff[Fh'(F_t)] = \ff \{\delta(\nabla
N_0^{-1}F)\cdot
\ff[h'(F_t)|F] \}.
\end{equation}
The right-hand-side of equation (\ref{eqIto4}) is the $L^2(\mathcal
{S},\ff)$-inner-product of\break $\delta(\nabla N_0^{-1}F)$ with $\ff
[h'(F_t)|F]^\ast= \ff[\overline{h'}(F_t)|F]$ (since $F$ and $F_t$ are
self adjoint), and this random variable is in the domain $\mathscr
{D}(\nabla)$. Hence, since $\delta$ and $\nabla$ are adjoint to each
other, we have
%
\begin{eqnarray}\label{eqIto5} \ff[h'(F_t)F] &=& \langle\nabla
N_0^{-1} F, \nabla\ff[\overline{h'}(F_t)|F]\rangle_{\mathscr{B}_2}\nonumber
\\[-8pt]
\\[-8pt] &=&
\int_\R\ff\tensor\ff \bigl(\nabla_s (N_0^{-1} F) \sh(\nabla
_s\ff
[\overline{h'}(F_t)|F])^\ast \bigr)\,ds.
\nonumber
\end{eqnarray}
To be clear: $\sh$ is the product $(A_1\tensor B_1)\sh(A_2\tensor
B_2) =
(A_1A_2)\tensor(B_2B_1)$. It is easy to check that this product is
associative and
distributive, as will be needed in the following.

Recall that $F_t = \sqrt{1-t}F + S_t$ and $\overline{h'}$ is equal to a
polynomial on a compact interval $K$ which contains the spectrum of
$F_t$. Hence, $\overline{h'}(F_t)$ is a (noncommutative) polynomial in
$F$ and $S_t$. Thus, the conditional expectation $\ff[\overline
{h'}(F_t)|F]$ is a polynomial $p(F)$ in $F$. We may thus employ the
chain rule of equation (\ref{eqchainrule}) to find that, for each $s$,
%
\begin{equation} \label{eqIto55} \nabla_s\ff[\overline
{h'}(F_t)|F] =
\del_F\ff[\overline{h'}(F_t)|F]\sh\nabla_s F.
\end{equation}
Taking adjoints yields
%
\begin{equation} \label{eqIto6} (\nabla_s\ff[\overline
{h'}(F_t)|F] )^\ast= (\nabla_sF)^\ast\sh\del_F\ff[h'(F_t)|F].
\end{equation}
Now we use the intertwining property of the free difference quotient
for the sum of free random variables
with respect to conditional expectation (see,~\cite{Voiculescu4},
Proposition 2.3) and the
simple scaling property $\partial_{aX}=a^{-1}\partial_X$ (for $a\in
\CC
$) to get
%
\begin{eqnarray} \label{eqchainrule2}
\del_F\ff[h'(F_t)|F] &=& \del_F\ff\bigl[h'\bigl(\sqrt{1-t} F + S_t\bigr)|F\bigr]\nonumber\\
&=&\sqrt{1-t}
\del_{\sqrt{1-t}F}\ff\bigl[h'\bigl(\sqrt{1-t} F + S_t\bigr)|F\bigr]\nonumber
\\[-8pt]
\\[-8pt]
&=&
\sqrt{1-t} \ff\tensor\ff\bigl[\del_{\sqrt{1-t} F + S_t} h'\bigl(\sqrt
{1-t} F + S_t\bigr)|F\bigr]\nonumber\\
&=&
\sqrt{1-t} \ff\tensor\ff[\del h'(F_t)|F].
\nonumber
\end{eqnarray}

\begin{remark} It is here, and only here, that the assumption that
$S_t$ is free from the $F$ is required.
\end{remark}

Combining equation (\ref{eqchainrule2}) with equations (\ref{eqIto5})
and (\ref{eqIto6}) yields
%
\begin{eqnarray}\label{eqIto7} &&  \ff[h'(F_t)F] \nonumber\hspace*{-35pt}
\\[-8pt]
\\[-8pt]&&  \qquad = \sqrt{1-t} \ff
\tensor
\ff \biggl(\int_\R\nabla_s (N_0^{-1} F) \sh
(\nabla_sF)^\ast \,ds \sh\ff\tensor\ff[\del h'(F_t)|F] \biggr).
\nonumber\hspace*{-35pt}
\end{eqnarray}
As for the second term in equation (\ref{eqIto3}), using property (3)
of conditional expectation (cf. Definition~\ref{defconditionalexpectation}) and taking expectations, we express
%
\begin{equation} \label{eqIto8} \ff\tensor\ff[\del h'(F_t)] = \ff
\tensor\ff (\ff\tensor\ff[\del h'(F_t)|F] ).
\end{equation}
Combining equations (\ref{eqIto3}), (\ref{eqIto7}) and (\ref{eqIto8}) yields
%
\begin{eqnarray}\label{eqIto9}
g'(t) &=& -\frac{1}{2}\ff\nonumber\hspace*{-30pt}\nonumber\\
&&{}\tensor\ff \biggl\{\int\nabla
_s(N_0^{-1}F)\sh
(\nabla_sF)^\ast\,ds\sh\ff\tensor\ff[\del h'(F_t)|F]\nonumber\hspace*{-30pt}\\[-8pt]\\[-8pt]
 &&\hphantom{\tensor\ff \biggl\{}- \ff
\tensor\ff
[\del h'(F_t)|F] \biggr\}\nonumber\hspace*{-30pt} \\
&=& -\frac{1}{2}\ff\tensor\ff \biggl\{ \!\biggl (\int\nabla
_s(N_0^{-1}F)\sh
(\nabla_sF)^\ast\,ds - 1\tensor1 \!\biggr)\sh\ff\tensor\ff[\del
h'(F_t)|F] \biggr\}.
\nonumber\hspace*{-30pt}
\end{eqnarray}
Integrating with respect to $t$ and using equation (\ref{eqFTC}) gives
%
\begin{eqnarray} \label{eqIto10}
&& \ff[h(S)] - \ff[h(F)] \nonumber\hspace*{-30pt}\\
&& \qquad = -\frac{1}{2}\ff\hspace*{-30pt}\\
&& \qquad  \quad {}\tensor\ff \biggl\{\! \biggl (\int\nabla_s(N_0^{-1}F)\sh
 (\nabla_sF)^\ast\,ds - 1\tensor1\! \biggr)\sh\int_0^1
 \ff\tensor\ff[\del h'(F_t)|F]\,dt \biggr\}.\nonumber\hspace*{-30pt}
\end{eqnarray}
Applying the noncommutative $L^1$--$L^\infty$ H\"older inequality
(which holds for the product $\sh$ on the algebra $\mathcal{S}$ since
$\sh$ is really just the natural product on the algebra $\mathcal
{S}\tensor\mathcal{S}^{\mathrm{op}}$; cf. Remark~\ref{rksharp})
gives us
%
\begin{eqnarray} \label{eqIto11}
&&|\ff[h(F)] - \ff[h(S)]| \nonumber
\\
&& \qquad \le\frac{1}{2}\ff\tensor\ff \biggl\{ \biggl|\int\nabla
_s(N_0^{-1}F)\sh
(\nabla_sF)^\ast\,ds - 1\tensor1 \biggr| \biggr\}\\
&& \hphantom{\frac{1}{2}\ff\tensor}\qquad \quad {}\times\int_0^1 \bigl\|
\ff
\tensor\ff[\del h'(F_t)|F] \bigr\|_{\mathcal{S}\tensor\mathcal
{S}}\,dt.
\nonumber
\end{eqnarray}
The norm $\| \cdot \|_{\mathcal{S}\tensor\mathcal{S}}$ is the
operator ($L^\infty$)
norm on the doubled abstract Wigner space. The conditional expectation
is an
$L^\infty$-contraction [cf. property (2) in Definition~\ref{defconditionalexpectation}],
and so the second term in equation (\ref{eqIto11}) satisfies
%
\begin{equation} \label{eqIto12} \int_0^1 \bigl\|\ff\tensor\ff
[\del
h'(F_t)|F] \bigr\|_{\mathcal{S}\tensor\mathcal{S}}\,dt
\le\int_0^1 \|\del h'(F_t) \|_{\mathcal{S}\tensor
\mathcal
{S}}\,dt.\vadjust{\goodbreak}
\end{equation}
Using equation (\ref{eqlemmab2}) with $Y=F_t$, note that
%
\begin{eqnarray} \label{eqIto13}
 \|\del h'(F_t) \|_{\mathcal{S}\tensor\mathcal{S}} &=&
 \biggl\|
\int_0^1d\alpha\int_\R\xi^2 \nu(d\xi) e^{i\alpha\xi
F_t}\tensor
e^{i(1-\alpha)\xi F_t} \biggr\|_{\mathcal{S}\tensor\mathcal{S}} \nonumber
\\[-8pt]
\\[-8pt]
&\le&\int_0^1 d\alpha\int_\R\xi^2 \nu(d\xi) \|e^{i\alpha\xi
F_t}\|_\mathcal{S}\bigl\|e^{i(1-\alpha)\xi F_t}\bigr\|_{\mathcal{S}}.
\nonumber
\end{eqnarray}
Both of the norm terms in the second line of equation (\ref{eqIto13})
are equal to $1$
since $F_t$ is self adjoint. This shows that $\|\del
h'(F_t)\|_{\mathcal{S}\tensor\mathcal{S}} \le\mathscr{I}_2(h)$.
Combining this with
equations (\ref{eqIto11}) and (\ref{eqIto12}) yields
%
\begin{eqnarray} \label{eqIto14}
&&|\ff[h(F)] - \ff[h(S)]| \nonumber
\\[-8pt]
\\[-8pt]
&& \qquad \le\frac{1}{2}\ff\tensor\ff\biggl \{
\biggl|\int
\nabla_s(N_0^{-1}F)\sh(\nabla_sF)^\ast\,ds - 1\tensor1
\biggr| \biggr\}
\cdot\mathscr{I}_2(h).
\nonumber
\end{eqnarray}
Inequality (\ref{eqIto14}) is close to the desired result, but as
proved it only holds for $h\in\mathcal{C}_2^{K,P}$. Now take any
$h\in
\mathcal{C}_2$, and fix an approximating sequence $h_n\in\mathcal
{C}_2^{K,P}$ as guaranteed by Theorem~\ref{theorempolynomialsdense}.
That theorem shows that $\mathscr{I}_2(h_n)\to\mathscr{I}_2(h)$, while
\begin{eqnarray*}
|\ff[h_n(F)]-\ff[h_n(S)]| &=& \biggl|\int h_n\,d\mu_F - \int h_n\,d\mu
_S \biggr| \to \biggl|\int h\,d\mu_F - \int h\,d\mu_S \biggr|\\
 &=& |\ff
[h(F)]-\ff[h(S)]|
\end{eqnarray*}
as $n\to\infty$, since the supports of $\mu_F$ and $\mu_S$ are
contained in $K$. This shows that inequality (\ref{eqIto14}) actually
holds for all $h\in\mathcal{C}_2$, and this concludes the proof.
\end{pf}

\begin{remark} \label{rkGammavsOU}
$\!\!\!$In equation (\ref{eqIto4}),
instead of using the Ornstein--Uhlenbeck operator, we might have used
the Clark--Ocone formula (Proposition~\ref{propClark-Ocone}). Tracking
this through the remainder of the proof would yield the related estimate
%
\begin{equation} \label{eqClark-Oconeestimate} d_{\mathcal{C}_2}(F,S)
\le\frac{1}{2}\ff\tensor\ff \biggl( \biggl|\int\Gamma(\nabla
_sF)\sh(\nabla
_s F)^\ast\,ds - 1\tensor1 \biggr| \biggr).
\end{equation}
This estimate is, in many instances, equivalent to equation (\ref{eqWasserstein1}) as far as convergence to the semicircular law is
concerned, as we discuss in Section~\ref{sectiondistanceestimate};
the formulation of equation (\ref{eqWasserstein1}) is ideally suited to
prove Corollary~\ref{corWasserstein}, which is why we have chosen this
presentation.
\end{remark}

\subsection{Distance estimates} \label{sectiondistanceestimate} We
begin by proving Corollary~\ref{corWasserstein}, which we restate here
for convenience with a little more detail.

\begin{corollaryWasserstein*} Let $f\in L^2(\R_+^2)$ be mirror-symmetric
and normalized $\|f\|_{L^2(\R_+^n)}=1$, and let $S$ be a standard
semicircular random variable. Then
\[ d_{\mathcal{C}_2}(I^S_2(f),S) \le\frac{1}{2}\sqrt
{\frac{3}{2}}\|f\cont{1}f\|_{L^2(\R_+^2)} = \frac{1}{2}\sqrt{\frac
{3}{2}} \sqrt{\E[I^S_2(f)^4]-2}.
\]
\end{corollaryWasserstein*}

\begin{pf} We will utilize the estimate of Theorem~\ref{theoremMalliavinestimate}
applied to the random variable $F = I^S_2(f)$
[which is indeed centred and in the domain~$\mathscr{D}(\nabla)$].
Note, from the definition, that $N_0^{-1}F = \frac12 F$ for a double
integral. From equation (\ref{eqgradientonintegrals}), we have
%
\begin{equation} \label{eqWass1} \nabla_tF = \nabla_t I^S_2(f) =
\int
f(t,t_2) 1\tensor dS_{t_2} + \int f(t_1,t)\,dS_{t_1}\tensor1.
\end{equation}
Using the fact that $f=f^\ast$, this yields
%
\begin{equation} \label{eqWass2} (\nabla_tF)^\ast= \int f(t_2,t)
1\tensor dS_{t_2} + \int f(t,t_1)\,dS_{t_1}\tensor1.
\end{equation}
(Note: the adjoint on tensor-product operators is, as one would expect,
$(A\tensor B)^\ast= A^\ast\tensor B^\ast$, contrary to the convention
on page 379 in \cite{BianeSpeicher}.) When multiplying equations (\ref{eqWass1}) and
(\ref{eqWass2}), one must keep in mind the product
formula (\ref{eqproductcont}) for Wigner integrals; in this context of
Wigner bi-integrals, the results are
%
\begin{eqnarray} \label{eqWass31}\qquad
 && \biggl(\int f(t,s_2) 1\tensor dS_{s_2} \biggr)\sh
\biggl(\int f(t_2,t) 1\tensor dS_{t_2} \biggr)\nonumber
\\[-8pt]
\\[-8pt]
&& \qquad = \int f(t,s_2)f(t_2,t) 1\tensor dS_{t_2}\,dS_{s_2} + \int f(t,s)f(s,t)\,ds 1\tensor1,
 \nonumber\\\label{eqWass32}
&& \biggl(\int f(t,s_2) 1\tensor dS_{s_2} \biggr)\sh \biggl(\int
f(t,t_1)\,
dS_{t_1}\tensor1 \biggr)\nonumber
\\[-8pt]
\\[-8pt]
&& \qquad = \int f(t,s_2)f(t,t_1)\,dS_{t_1}\tensor dS_{s_2},
\nonumber\\\label{eqWass33}
 && \biggl(\int f(s_1,t)\,dS_{s_1}\tensor1 \biggr)\sh \biggl(\int
f(t_2,t)
1\tensor dS_{t_2} \biggr)\nonumber
\\[-8pt]
\\[-8pt]
&& \qquad =\int f(s_1,t)f(t_2,t)\,dS_{s_1}\tensor dS_{t_2},
\nonumber\\\label{eqWass34}
 && \biggl(\int f(s_1,t)\,dS_{s_1}\tensor1 \biggr)\sh \biggl(\int
f(t,t_1)\,
dS_{t_1}\tensor1 \biggr)\nonumber
\\[-8pt]
\\[-8pt]
&& \qquad = \int f(s_1,t)f(t,t_1)\,dS_{s_1}\,dS_{t_1}\tensor1
 + \int f(s,t)f(t,s)\,ds 1\tensor1.
\nonumber
\end{eqnarray}
Integrating with respect to $t$ and using the identity
$f(s,t)=\overline
{f(t,s)}$, we then have
%
\begin{eqnarray} \label{eqWass4}
 2\int\nabla_t(N_0^{-1} F)\sh(\nabla_t F)^\ast\,dt &=&
\int\!\!\int f(t,s_2)f(t_2,t)\,dt 1\tensor dS_{t_2}\,dS_{s_2}\nonumber\\
&&{} + \int\!\!\int
f(t,s_2)f(t,t_1)\,dt \,dS_{t_1}\tensor dS_{s_2} \nonumber\\
&&{}+ \int\!\!\int f(s_1,t)f(t_2,t)\,dt \,dS_{s_1}\tensor dS_{t_2}\\
&&{} + \int\!\!\int
f(s_1,t)f(t,t_1)\,dt\,dS_{s_1}\,dS_{t_1}\tensor1 \nonumber\\
&&{}+ 2\int|f(s,t)|^2\,dt\,ds 1\tensor1.\nonumber
\end{eqnarray}
Now using the normalization $\|f\|_{L^2(\R^2)} = 1$, and making use of
contraction notation (cf. Definition~\ref{defcontraction}), we have
%
\begin{eqnarray} \label{eqWass5}\qquad
 && 2\biggl (\int\nabla_t(N_0^{-1} F)\sh(\nabla_t
F)^\ast\,
dt-1\tensor1 \biggr)\nonumber\\ && \qquad =
\int\overline{f}\cont{1}\overline{f}(s_2,t_2) 1\tensor
dS_{t_2}\,dS_{s_2} + \int\overline{f}\cont{1} f(s_2,t_1)\,
dS_{t_1}\tensor dS_{s_2}
\\
&& \qquad \quad {}+ \int f\cont{1}\overline{f}(s_1,t_2)\,dS_{s_1}\tensor dS_{t_2} +
\int
f\cont{1}f(s_1,t_1)\,dS_{s_1}\,dS_{t_1}\tensor1.
\nonumber
\end{eqnarray}
We now employ Theorem~\ref{theoremMalliavinestimate}. Equation (\ref{eqWasserstein1}) states that
%
\begin{equation} \label{eqL^1Wass}  d_{\mathcal{C}_2}(F,S) \le\tfrac
{1}{2} \bigl\|\langle\nabla(N_0^{-1}F)\sh(\nabla F)^{\ast},\1_{\R
_+}\rangle_{L^2(\R_+)} -1\tensor1 \bigr\|_{L^1(\mathcal{S}\tensor
\mathcal{S},\ff\tensor\ff)}.\hspace*{-40pt}
\end{equation}
In any $W^\ast$-probability space, $\| \cdot \|_{L^1}\le\| \cdot
 \|
_{L^2}$; we will estimate the $L^2(\mathcal{S}\tensor\mathcal{S},\allowbreak \ff
\tensor\ff)$ norm. It is useful to relabel the indices in equation (\ref{eqWass5}) and group them according to different orders of (bi)chaos;
the right-hand side of that equation is equal~to
\begin{eqnarray*} && \int\overline{f}\cont{1}\overline{f}(t_2,t_1)
1\tensor dS_{t_1}\,dS_{t_2} \\
&& \qquad {}+ \int [\overline{f}\cont{1}f(t_2,t_1)+f\cont{1}\overline
{f}(t_1,t_2) ]\,dS_{t_1}\tensor dS_{t_2} \\
&& \qquad {}+ \int f\cont{1}f(t_1,t_2)\,dS_{t_1}\,dS_{t_2}\tensor1.
\end{eqnarray*}
A simple calculation using the fact that $f=f^\ast$ shows\vspace*{-1.5pt} that
$\overline{f}\cont{1}f(t_2,t_1) = \overline{f}\cont{1}f (t_1,t_2)$. The
three integrals above are in orthogonal orders of chaos; employing the
Wigner bisometry~\ref{eqWignerbisometry}, we have
%
\begin{eqnarray} \label{eqWass6}
&&4 \bigl\|\langle\nabla(N_0^{-1}F)\sh(\nabla F)^{\ast},\1_{\R
_+}\rangle
_{L^2(\R_+)} -1\tensor1 \bigr\|_{L^2(\mathcal{S}\tensor\mathcal
{S},\ff
\tensor\ff)}^2\nonumber\hspace*{-35pt}
\\[-8pt]
\\[-8pt]
&& \qquad = \|\overline{f}\cont{1}\overline{f}\|_{L^2(\R_+^2)}^2
+ \| \overline{f}\cont{1}f + f\cont{1}\overline{f} \|_{L^2(\R
_+^2)}^2 +
\| f\cont{1} f\|_{L^2(\R_+^2)}^2.
\nonumber\hspace*{-35pt}
\end{eqnarray}
Another simple calculation, again using the identity $f(s,t)=\overline
{f(t,s)}$, shows that
%
\begin{equation} \label{eqWass71}  \|f\cont{1}f\|_{L^2(\R_+^2)}^2
= \|
\overline{f}\cont{1}\overline{f}\|_{L^2(\R_+)^2}^2 = \int_{\R
_+^3} dt\,
ds_1\,ds_2 |f(t,s_1)|^2|f(t,s_2)|^2,\hspace*{-40pt}
\end{equation}
while
%
\begin{equation} \label{eqWass72}  \|f\cont{1}\overline{f}\|
_{L^2(\R
_+^2)}^2 = \|\overline{f}\cont{1}f\|_{L^2(\R_+)^2}^2 = \int_{\R_+^3}
dt\,ds_1\,ds_2 f(t,s_1)^2\overline{f(t,s_2)}^2.\hspace*{-35pt}
\end{equation}
Hence $\|f\cont{1}\overline{f}\|_{L^2(\R_+^2)} = \|\overline
{f}\cont
{1}f\|_{L^2(\R_+^2)} \le\|f\cont{1}f\|_{L^2(\R_+)}$. Using the
triangle inequality in equation (\ref{eqWass6}) then gives us the estimate
%
\begin{eqnarray} \label{eqWass8}
 && \bigl\|\langle\nabla(N_0^{-1}F)\sh(\nabla F)^{\ast},\1_{\R
_+}\rangle
_{L^2(\R_+)} -1\tensor1 \bigr\|_{L^2(\mathcal{S}\tensor\mathcal
{S},\ff
\tensor\ff)}^2\nonumber
\\[-8pt]
\\[-8pt]&& \qquad \le\tfrac32\|f\cont{1}f\|_{L^2(\R_+^2)}^2,
\nonumber
\end{eqnarray}
and so equation (\ref{eqL^1Wass}) and the ensuing discussion imply
%
\begin{equation} \label{eqWass9} d_{\mathcal{C}_2}(F,S) =
d_{\mathcal
{C}_2}(I^S_s(f),S) \le\frac12\sqrt{\frac32}\|f\cont{1}f\|_{L^2(\R
_+^2)}.
\end{equation}
Now, as calculated in equation (\ref{eq4thmomentcont3}) (in this
instance with $n=2$),
%
\begin{equation} \label{eqWass10} \ff(I^S_2(f)^4) = 2 + \|f\cont
{1}f\|
_{L^2(\R_+^2)}^2.
\end{equation}
Equations (\ref{eqWass9}) and (\ref{eqWass10}) together conclude the proof.
\end{pf}

\begin{remark} \label{rkonlydoubleintegrals} At first glance it
might seem that calculations like those in the proof of Corollary~\ref{corWasserstein} could be employed to prove similar quantitative
results for Wigner integrals $I^S_n$ of arbitrary order $n\ge2$. Note,
however, that the mirror symmetry of $f$ was used in different ways at
several points in the above proof. In practice, if one tries to
generalize these techniques to~$I^S_3$, in fact $f$ must be \textit{fully
symmetric}. The range of $I^S_n$ on fully symmetric functions is a very
small subspace of the full $n$th Wigner chaos, and so we do not have
quantitative bounds for generic higher integrals.
\end{remark}

\begin{remark} As a quick illustration, we use the first inequality in
Corollary~\ref{corWasserstein} to refine Corollary~\ref{corBreuerMajor}
in the case $n=2$ and the random variables $X_k$ are freely
independent $S(0,1)$ random variables; in particular, $\rho(k) =
\delta
_{k0}$. In this case, one can take these random variables to be such\vadjust{\goodbreak}
that $X_k = S_{k+1} - S_{k}$, $k\geq0$, so that $V_m = I_2^S(f_m)$, with
\[
f_m(x,y) = \frac{1}{\sqrt{m}} \sum_{k=0}^{m-1} \1_{(k,k+1]}(x)\1
_{(k,k+1]}(y).
\]
Elementary computations now yield $ \|f_m\stackrel{1}{\frown} f_m\|
_{L^2(\mathbb{R}_+^2)} = 1/\sqrt{m}$, and therefore
\[
d_{\mathcal{C}_2}(V_m,S) \leq\frac12\sqrt{\frac{3}{2m} },
\]
which is consistent with usual Berry--Esseen estimates.
\end{remark}

In light of Theorem~\ref{thm4thmomentsemicircle}, the proof of
Corollary~\ref{corWasserstein} shows that convergence of the quantity
on the right-hand side of equation (\ref{eqWasserstein1}) to $0$ is
equivalent to convergence of $F$ to $S$ in law, at least in the case of
double Wigner integrals. We conclude this paper with a collection of
other equivalences, stated in terms of the gradient operator, in the
class of double Wigner integrals; whether they hold for higher orders,
or more generally on the domain $\mathscr{D}(\nabla)$, is left as an
open question for further investigation. To simplify matters, we
restrict to the real case for the following.

\begin{theorem} \label{thmequivalencesforI2} Let $(f_k)_{k\in\N}$
be a sequence of fully symmetric functions in $L^2(\R_+^2)$, each
normalized $\|f_k\|_{L^2(\R_+)^2}=1$, and set $F_k = I^S_2(f_k)$. Then
for each $k$
\[
F_k = \delta(\nabla N_0^{-1}F_k) = \delta(\Gamma\nabla F_k).
\]
Moreover, the following four conditions are equivalent:
\begin{longlist}[(3)]
\item[(1)] $F_k$ converges in law to the standard semicircular
distribution $S(0,1)$;
\item[(2)] ${\int\nabla_t(N_0^{-1}F_k)\sh(\nabla_t
F_k)^\ast\,dt}$ converges to $1\tensor1$ in $L^2(\mathcal{S}\tensor
\mathcal{S},\ff\tensor\ff)$;
\item[(3)] ${\int\Gamma(\nabla_tF_k)\sh(\nabla_t
F_k)^\ast
\,dt}$ converges to $1\tensor1$ in $L^2(\mathcal{S}\tensor\mathcal
{S},\ff\tensor\ff)$;
\item[(4)] ${\int\langle\hspace*{-1.5pt}\langle\Gamma\nabla
_tF_k,\Gamma
\nabla_tF_k\rangle\hspace*{-1.5pt}\rangle\,dt}$ converges to $1$ in $L^2(\mathcal
{S},\ff)$.
\end{longlist}
\end{theorem}

 The pairing $\langle\hspace*{-1.5pt}\langle \cdot , \cdot \rangle\hspace*{-1.5pt}
\rangle
\dvtx (\mathcal{S}\tensor\mathcal{S})^2\to\mathcal{S}$ is defined by
$\langle\hspace*{-1.5pt}\langle X,Y\rangle\hspace*{-1.5pt}\rangle= (1_{\mathcal{S}}\tensor\ff
)[X\sh
Y]$. For example,
\[
\biggl\langle\hspace*{-1.5pt}\biggl\langle\int f(t_1) 1\tensor dS_{t_1}, \int g(t_2) 1\tensor
dS_{t_2} \biggr\rangle\hspace*{-1.5pt}\biggr\rangle= \int f(t_1)g(t_2)\,dS_{t_2}\,dS_{t_1} + \int
f(t)g(t)\,dt,
\]
where we have used the product formula of equation (\ref{eqproductcont}). On the other hand,
$\langle\hspace*{-1.5pt}\langle1\tensor
dS_{t_1},dS_{t_2}\tensor1\rangle\hspace*{-1.5pt}\rangle= 0$ since $\ff(dS_1)=0$.

\begin{pf} Equations (\ref{eqWass6}) and (\ref{eqWass10}) in the
proof of Corollary~\ref{corWasserstein} show that, in the case that
$f$ is real-valued,
\[
 \biggl\|\int\nabla_t(N_0^{-1}F)\sh(\nabla_t F)\,dt - 1\tensor
1 \biggr\|
_{L^2(\mathcal{S}\tensor\mathcal{S},\ff\tensor\ff)}^2 =\frac
{3}{2}\bigl(\ff
(F^4)-2\bigr),
\]
where $F=I_2^S(f)$. In light of Theorem~\ref{thm4thmomentsemicircle},
this proves the equivalence (1)${}\Longleftrightarrow{}$(2).
The bound~\ref{eqClark-Oconeestimate} shows that (3)${}\Longrightarrow
{}$(1), and so to prove the equivalence of (1) and (3) it suffices (due
to Theorem~\ref{thm4thmomentcont}) to prove that the condition
$f_k\cont{1}f_k\to0$ implies (3). To that end, we adopt the standard
notation $x\vee y = \max\{x,y\}$ and $x\wedge y = \min\{x,y\}$. The
following identity is easily proved:\looseness=-1
%
\begin{eqnarray} \label{eqequiv1} &&\int_0^\infty\biggl (\int_{x\vee
y}^\infty f_k(x,t)f_k(y,t)\,dt \biggr)^2\,dx\,dy\nonumber
\\[-8pt]
\\[-8pt]&& \qquad = \int_0^\infty
\biggl(\int
_0^{x\wedge y} f_k(x,t)f_k(y,t)\,dt \biggr)^2\,dx\,dy.
\nonumber
\end{eqnarray}\looseness=0
The following equivalence was proved in \cite{nunugio}.
%
\begin{eqnarray} \label{eqequiv2} &\displaystyle f_k\cont{1} f_k\to0 \qquad \mbox{in
}L^2(\R_+^2) \quad \mbox{iff}&\nonumber
\\[-8pt]
\\[-8pt]
&\displaystyle \int_0^\infty \biggl(\int
_0^{x\wedge
y} f_k(x,t)f_k(y,t)\,dt \biggr)^2\,dx\,dy\to0.&\nonumber
\end{eqnarray}
Note that equations (\ref{eqequiv1}) and (\ref{eqequiv2}) imply that,
if $f_k\cont{1}f_k\to0$, then the three functions
%
\begin{eqnarray} \label{eqvanish} &&\int_0^{x\wedge y}
f_k(x,t)f_k(y,t)\,
dt, \qquad\int_{x\vee y}^\infty f_k(x,t)f_k(y,t)\,dt,\nonumber
\\[-8pt]
\\[-8pt] &&\int
_{x\vee y}^{x\wedge y} f_k(x,t)f_k(y,t)\,dt
\nonumber
\end{eqnarray}
each vanish in the limit. Note also that the action of $\Gamma$ on the
biprocess $\nabla_t F_k$ is, as in the classical case, to restrict
stochastic integrals to the interval $[0,t]$.
%
\begin{equation} \label{eqequiv3}\qquad
\Gamma\nabla_t F_k = \int
f_k(t,t_2)\1_{t_2\le t} 1\tensor dS_{t_2} + \int f_k(t_1,t)\1_{t_1\le
t}\,dS_{t_1}\tensor1.
\end{equation}
The present symmetry assumptions on $f_k$, imply that $(\nabla
F_k)^\ast= \nabla F_k$. Proceeding with calculations like those in the
proof of Corollary~\ref{corWasserstein}, using the symmetry and
$L^2$-normalization of $f_k$, we then have
%
\begin{eqnarray} \label{eqequiv4}
\int\Gamma\nabla_t F_k\sh(\nabla_t F_k)^\ast\,dt - 1\tensor1 &=&
\int \biggl(\int_{t_2}^\infty f_k(t,t_1)f_k(t,t_2)\,dt \biggr)
1\tensor
dS_{t_1}\,dS_{t_2} \nonumber\\
&&{}+ \int \biggl(\int_{t_2}^\infty f_k(t,t_1)f_k(t,t_2)\,dt \biggr)\,
dS_{t_1}\tensor dS_{t_2}\nonumber
\\[-8pt]
\\[-8pt]
&&{}+ \int \biggl(\int_{t_1}^\infty f_k(t_1,t)f_k(t_2,t)\,dt \biggr)\,
dS_{t_1}\tensor dS_{t_2} \nonumber\\
&&{}+ \int \biggl(\int_{t_1}^\infty f_k(t_1,t)f(t_2,t)\,dt \biggr)\,
dS_{t_1}\,dS_{t_2}\tensor1.
\nonumber
\end{eqnarray}
Using Fubini's theorem, we can calculate that the $L^2$ norm of each of
the four terms in equation (\ref{eqequiv4}) is given by
\begin{eqnarray*}
&&\int_{\R_+^2} dx\,dy \biggl(\int_x^\infty f_k(t,x)f_k(t,y)\,dt
\biggr)^2 \\[-2pt]
&& \qquad = \int_0^\infty\,dx\int_0^x\,dy \biggl(\int_{x\vee y}^\infty
f_k(t,x)f_k(t,y)\,dt \biggr)^2\\[-2pt]
&& \qquad \quad {}
+ \int_0^\infty\,dy\int_0^y\,dx \biggl(\int_{x\wedge y}^\infty
f_k(t,x)f_k(t,y)\,dt \biggr)^2.
\end{eqnarray*}
Hence, if $f_k\cont{1}f_k\to0$, then equation (\ref{eqvanish}) shows
that each of these two terms vanishes in the limit. This proves the
implication (3)${}\Longrightarrow{}$(1).

 For the final equivalence, we use the explicit representation
\begin{eqnarray*}
&&\langle\hspace*{-1.5pt}\langle\Gamma\nabla_tF_k,\Gamma\nabla_t
F_k\rangle\hspace*{-1.5pt}\rangle\\[-2pt]
&& \qquad = 2\int f_k(t,v)^2\1_{v\le t}\,dv + \int f_k(t_1,t)f_k(t_2,t)\1
_{t_1\le t}\1_{t_2\le t}\,dS_{t_1}\,dS_{t_2}.
\end{eqnarray*}
Integrating with respect to $t$ and using equations (\ref{eqequiv1})
and (\ref{eqequiv2}) as above proves the equivalence
(1)${}\Longleftrightarrow{}$(4).
\end{pf}

\begin{remark} As demostrated in \cite{BianeSpeicher}, Theorem~4.12,
the quantity\break $\int\langle\hspace*{-1.5pt}\langle\Gamma\nabla_t F_k, \Gamma\nabla_t
F_k\rangle\hspace*{-1.5pt}\rangle \,dt$ in condition (4) of Theorem~\ref{thmequivalencesforI2} can be interpreted as the ``quadratic variation''
 of an
appropriate free Brownian martingale. Note that quadratic variations
play a crucial role in the original proof of Theorem~\ref{thmNP1}, as
originally given in \cite{nunugio}.
\end{remark}

\begin{remark} Once again, one might expect that calculations like
those above would show the equivalence of items (1)--(4) in Theorem
\ref{thmequivalencesforI2} for any order of chaos (higher than $1$), as
was proved in the classical case in \cite{nunugio}. In principle, this
may be possible for fully symmetric kernels $f$, but in orders $\ge3$
of Wigner chaos, such kernels span only a tiny subspace of all
stochastic integrals. Indeed, it is an interesting open question if a
counter-example to these equivalences can be found in the third chaos;
until now, the authors have not been able to find one, but suspect that
Theorem~\ref{thmequivalencesforI2} does not generally hold in the
free context.
\end{remark}

\begin{appendix}
\setcounter{equation}{0}
\section*{\texorpdfstring{Appendix: Proof of Theorem \lowercase{\protect\ref{theorempolynomialsdense}}}
{Appendix: Proof of Theorem 3.20}} \label{appendix}

We break the proof into four steps. First we show that it is sufficient
to consider only those $h\in\mathcal{C}_2$ that arise as Fourier
transforms of compactly-supported measures, in Lemma~\ref{lemmaapprox1}. Next we reduce to those $h$ that are Fourier transforms of measures
with a smooth, compactly-supported density, in Lemma~\ref{lemmaapprox2}. In Lemma~\ref{lemmaRudin}, we show (following
\cite{Rudin}, Theorem~7.26) that\vadjust{\goodbreak} there is a polynomial approximate identity on any
symmetric compact interval. Finally, we use this approximate identity
locally to approximate any smoothly-arising $h$ by local polynomials on
the Fourier side in Lemma~\ref{lemmaapprox3}, completing the proof.
The proof will actually show that a space smaller than~$\mathcal
{C}_2^{K,P}$ is appropriately dense: the local polynomials may be
assumed to live in the Schwartz space $\mathcal{S}(\R)$ of
rapidly-decaying smooth functions.\vspace*{-2pt}

\begin{lemma} \label{lemmaapprox1} Let $h\in\mathcal{C}_2$. There
exists a sequence of compactly-supported complex measures $\nu_n$ such
that, setting $h_n = \widehat{\nu}_n$:
\begin{longlist}[(2)]
\item[(1)] $\mathscr{I}_2(h_n)\to\mathscr{I}_2(h)$;
\item[(2)] if $\mu$ is any finite measure, then $\int h_n\,d\mu\to
\int
h\,d\mu$.\vspace*{-2pt}
\end{longlist}
\end{lemma}

\begin{pf} Let $h = \widehat{\nu}$ where $\nu$ is a complex measure
satisfying $\int\xi^2 |\nu|(d\xi)<\infty$. Let $\nu_n(d\xi) = \1
_{|\xi
|\le n}\nu(d\xi)$, and take $h_n = \widehat{\nu_n}$. Then
%
\begin{equation} \label{eqlemmaapprox11} \mathscr{I}_2(h_n) =
\int
_{-n}^n \xi^2 |\nu|(d\xi).
\end{equation}
Since $h\in\mathcal{C}_2$, the function $\xi\mapsto\xi^2$ is in
$L^1(|\nu|)$; hence, by the dominated convergence theorem, the
integrals in equation (\ref{eqlemmaapprox11}) converge to $\int\xi
^2 |\nu|(d\xi)=\mathscr{I}_2(h)$ as desired. Now, for any $x\in\R$,
%
\begin{equation} \label{eqlemmaapprox12}
|h_n(x)-h(x)| = \biggl|\int e^{ix\xi}\bigl(\1_{|\xi|\le n}-1\bigr) \nu(d\xi
)
\biggr| \le\int\1_{|\xi|>n} |\nu|(d\xi).
\end{equation}
The integrand $\1_{|\xi|>n}$ converges pointwise to $0$ and is bounded,
so since $|\nu|$ is a finite measure, the dominated convergence theorem
shows that $h_n\to h$ pointwise. Finally, note also that $\|h_n\|
_{L^\infty} \le\int|\nu_n| \le\int|\nu| <\infty$, and so since~$\mu$
is a finite measure, one more application of the Dominated
Convergence Theorem shows that $\int h_n\,d\mu\to\int h\,d\mu$ as desired.\vspace*{-2pt}
\end{pf}

\begin{lemma} \label{lemmaapprox2} Let $h\in\mathcal{C}_2$ with
$h=\widehat{\nu}$ for some compactly-supported complex measure $\nu$.
There exists a sequence of smooth $\C$-valued functions $\psi_n\in
C_c^\infty$ such that, setting $h_n = \widehat{\psi}_n$:
\begin{longlist}[(2)]
\item[(1)] $\mathscr{I}_2(h_n)\to\mathscr{I}_2(h)$;
\item[(2)] if $\mu$ is any finite measure then $\int h_n\,d\mu\to
\int
h\,d\mu$.\vspace*{-2pt}
\end{longlist}
\end{lemma}

\begin{pf} Let $\phi\in C_c^\infty$ be a nonnegative smooth
compactly supported function, such that $\int\phi(\xi)\,d\xi= 1$. Let
$\phi_n(\xi) = n \phi(\xi/n)$. Define $\psi_n = \phi_n\ast\nu
$; then
$\psi_n\to\nu$ weakly. Note that $\operatorname{supp} \phi_n\subset
\operatorname
{supp} \phi$. Since $\nu$ is compactly supported, there is thus a
single compact interval $K$ that contains the supports of $\psi_n$ for
all $n$ along with the support of $\nu$; moreover, the functions $\psi
_n$ are all smooth since $\phi_n$ is smooth. Set $h_n = \widehat{\psi
}_n$. Hence,
%
\begin{equation} \label{eqlemmaapprox21}
\mathscr{I}_2(h_n) = \int_K \xi^2 |\psi_n(\xi)|\,d\xi\to\int_K
\xi^2
|\nu|(d\xi) = \mathscr{I}_2(h),
\end{equation}
where the convergence follows from the weak convergence of $|\psi_n|$
to $|\nu|$ and the continuity\vadjust{\goodbreak} of $\xi\mapsto\xi^2$ on the compact set
$K$. For the second required convergence, we use Fubini's theorem,
%
\begin{eqnarray} \label{eqlemmaapprox22}
 \int h_n(x) \mu(dx) &=& \int\widehat{\psi}_n(x) \mu
(dx) = \int\mu(dx) \int e^{ix\xi} \psi_n(\xi)\,d\xi\nonumber
\\[-9pt]
\\[-9pt]
&=& \int\psi_n(\xi)\,d\xi\int e^{ix\xi} \mu(dx) = \int\widehat
{\mu
}(\xi)\psi_n(\xi)\,d\xi,
\nonumber
\end{eqnarray}
where the application of Fubini's theorem is justified by the fact that
the function $(x,\xi)\mapsto e^{ix\xi}\psi_n(\xi)$ is in $L^1(\mu
\times
d\xi)$ since $\psi_n\in L^1(d\xi)$ and $\mu$ is a~finite measure. The
function $\widehat{\mu}$ is continuous and bounded since $\mu$ is
finite, and so since $\psi_n\to\nu$ weakly and $\operatorname{supp} \psi
_n\subset K$ for each $n$,
%
\begin{equation} \label{eqlemmaapprox23}
\int\widehat{\mu}(\xi)\psi_n(\xi)\,d\xi= \int_K \widehat{\mu
}(\xi)\psi
_n(\xi)\,d\xi\to\int_K \widehat{\mu}(\xi) \nu(d\xi).
\end{equation}
The function $(x,\xi)\mapsto e^{ix\xi}$ is in $L^1(\mu\times|\nu|)$
since both are finite measures, and so we may apply Fubini's theorem
again to find that
%
\begin{eqnarray} \label{eqlemmaapprox24}
\int_K \widehat{\mu}(\xi) \nu(d\xi) &=& \int\nu(d\xi) \int
e^{ix\xi}
\mu(dx) = \int\mu(dx)\int e^{ix\xi} \nu(d\xi)\nonumber
\\[-9pt]
\\[-9pt] &=& \int\widehat
{\nu
}(x) \mu(dx),
\nonumber
\end{eqnarray}
where the first equality uses the fact that $\operatorname{supp} \nu
\subseteq K$. Equations (\ref{eqlemmaapprox22})--(\ref{eqlemmaapprox24})
combine to show that $\int h_n\,d\mu\to\int h\,d\mu$, as required.\vspace*{-2pt}
\end{pf}

\begin{lemma} \label{lemmaRudin} Let $r>0$. There is a sequence of
real polynomials $q_n$ such that, for any function $f$ continuous on
$\R
$ and equal to $0$ outside of $[-r,r]$, the functions
%
\begin{equation} \label{eqlemmaapprox31}
f_n(x) = \int_{-r}^{r} f(x-t)q_n(t)\,dt = \bigl (f\ast\bigl(q_n\1
_{[-r,r]}\bigr) \bigr)(x)
\end{equation}
are polynomials that converge uniformly to $f$ on $[-r,r]$.\vspace*{-2pt}
\end{lemma}

\begin{pf} This is proved in \cite{Rudin}, Theorem~7.26, in the case
$r=1$ with polynomials $c_n(1-x^2)^n$ for appropriate normalization
constants $c_n$. Rudin only states (and uses) the uniform convergence
on $[0,1]$, but it is easy to check that the proof yields uniform
convergence on $[-1,1]$. Rescaling the polynomials
%
\begin{equation} \label{eqBernstein} q_n(x) = \frac
{c_n}{r^{2n+1}}(r^2-x^2)^n
\end{equation}
gives us the desired result. To be clear: the functions $f_n$ in
equation (\ref{eqlemmaapprox31}) are polynomials due to the following
change of variables:
%
%
%
%
%
\begin{eqnarray} \label{eqlemmaapprox32} f_n(x) &=& \int_{-r}^r
f(x-t)q_n(t)\,dt = \int_{x-r}^{x+r} f(x-t)q_n(t)\,dt\nonumber
\\[-9pt]
\\[-9pt] &=& \int_{-r}^r
f(t)q_n(x+t)\,dt,
\nonumber\vadjust{\goodbreak}
\end{eqnarray}
where the second equality is justified by the fact that $f(x-t)=0$
unless $t\in[x-r,x+r]$.
%
\end{pf}

\begin{lemma} \label{lemmaapprox3} Let $h\in\mathcal{C}_2$ with
$h=\widehat{\psi}$ for some $\psi\in C_c^\infty$. Let $K\subset\R
$ be a
compact interval. There exists a sequence $\psi_n$ of functions in the
Schwartz space $\mathcal{S}(\R)$ such that the functions $h_n =
\widehat
{\psi}_n$ are in $\mathcal{C}_2^{K,P}$, and:
\begin{longlist}[(2)]
\item[(1)] $\mathscr{I}_2(h_n)\to\mathscr{I}_2(h)$;
\item[(2)] if $\mu$ is a finite measure supported in $K$ then $\int
h_n\,d\mu\to\int h\,d\mu$.
\end{longlist}
\end{lemma}

\begin{pf} Choose $r > \sup\{|x|\dvtx x\in K\}$. Let $\phi\in
C_c^\infty$
be nonnegative, with support contained in $[-r,r]$, such that $\phi
(x)=1$ for $x\in K$ (which is possible since $K$ is strictly contained
in $[-r,r]$). For convenience, set $p_n = q_n\1_{[-r,r]}$ where $q_n$
is the Bernstein polynomial of equation (\ref{eqBernstein}). Define
%
\begin{equation} \label{eqpsin}
\psi_n = \psi- [\widehat{\psi}\cdot\phi^2 ]^\vee+
\bigl[[(\widehat{\psi}\phi)\ast p_n]\cdot\phi \bigr]^\vee.
\end{equation}
Note: for a Schwartz function $\gamma\in\mathcal{S}(\R)$, the function
$\gamma^\vee=\check{\gamma}$ denotes the inverse Fourier transform of
$\gamma$,
\[
\gamma^\vee(\xi) = \check{\gamma}(\xi) = \frac{1}{2\pi}\int
e^{-ix\xi}
\gamma(x)\,dx.
\]
Since $\widehat{\psi}\phi\in C_c^\infty$, the convolution with
$p_n$ is
well defined and $C^\infty$; cutting off with $\phi$ again yields a
$C_c^\infty$ function, and so the inverse Fourier transform is a
Schwartz function. Similarly, $\phi^2$ is $C_c^\infty$ and $\widehat
{\psi}\in\mathcal{S}(\R)$, so their product is a Schwartz function, as
is its inverse Fourier transform. Thus, $\psi_n\in\mathcal{S}(\R)$. Now
we compute
\[
\widehat{\psi}_n = \widehat{\psi}-\widehat{\psi}\cdot\phi^2 +
[(\widehat
{\psi}\phi)\ast p_n]\cdot\phi
= \widehat{\psi}\cdot(1-\phi^2) + [(\widehat{\psi}\phi)\ast
p_n]\cdot
\phi.
\]
Since $\phi(x)^2 = 1$ for $x\in K$, we have $\widehat{\psi}_n(x) =
[(\widehat{\psi}\phi)\ast p_n](x)$ for $x\in K$. Since the function
$f=\widehat{\psi}\phi$ is continuous and equal to $0$ outside of
$[-r,r]$, equations (\ref{eqlemmaapprox31}) and (\ref{eqlemmaapprox32}) show that $\widehat{\psi}_n$ is a polynomial on $K$. Moreover,
$\psi_n$ is rapidly decaying and smooth, so $\int\xi^2|\psi_n(\xi
)|\,
d\xi<\infty$. Thus $h_n=\widehat{\psi}_n \in\mathcal{C}_2^{K,P}$ as
required. We must now verify conditions (1) and (2) of the lemma.

First, we compute that
%
\begin{equation} \label{eqinvFourier} \qquad \psi_n(x)-\psi(x) = \frac
{1}{2\pi}\int e^{-i\xi x} \bigl[[(\widehat{\psi}\phi)\ast p_n](\xi
)-\widehat{\psi}(\xi)\phi(\xi) \bigr]\phi(\xi)\,d\xi.
\end{equation}
Following this we make the straightforward estimate
%
\begin{eqnarray} \label{eqestimate1}
 |\psi_n(x)-\psi(x)| &\le&\frac{1}{2\pi}\int
|[(\widehat{\psi}\phi)\ast p_n](\xi)-\widehat{\psi}(\xi)\phi
(\xi)
|\phi(\xi)\,d\xi\nonumber
\\[-8pt]
\\[-8pt]
& =& \frac{1}{2\pi}\int_{-r}^r |[(\widehat{\psi}\phi)\ast
p_n](\xi
)-\widehat{\psi}(\xi)\phi(\xi) |\phi(\xi)\,d\xi.
\nonumber
\end{eqnarray}
Lemma~\ref{lemmaRudin} shows that $(\widehat{\psi}\phi)\ast p_n$
converges to $\widehat{\psi}\phi$ uniformly on $[-r,r]$. Hence, since
the integrand in equation (\ref{eqestimate1}) converges to $0$
uniformly on the (compact) domain of integration, it follows that $\psi
_n(x)\to\psi(x)$ for each~$x$.

We must now show that $\mathscr{I}_2(h_n)\to\mathscr{I}_2(h)$ (recall
that $h_n = \widehat{\psi}_n$ and $h=\widehat{\psi}$). This will follow
from the stronger claim that $\mathscr{I}_2(h_n-h)\to0$, which we now
show to be true. We compute as follows.
\[
\mathscr{I}_2(h_n-h) = \int|\psi_n(\xi)-\psi(\xi)|\xi^2\,d\xi=
\int
g_n(\xi) \frac{d\xi}{1+\xi^2},
\]
where $g_n(\xi) = \xi^2(1+\xi^2)|\psi_n(\xi)-\psi(\xi)|$. We
make this
transformation so we can use the finite measure $\upsilon(d\xi) =
d\xi
/(1+\xi^2)$ in the following estimates. Since $\psi_n\to\psi$
pointwise, it follows that $g_n\to0$ pointwise. In order to use
a~uniform integrability condition, we wish to bound the $L^2(\upsilon
)$-norm of $g_n$. To that end, we compute
%
\begin{equation} \label{eqL21} \hspace*{35pt} \|g_n\|_{L^2(\upsilon)}^2 = \int
g_n(\xi)^2 \upsilon(d\xi) = \int|\psi_n(\xi)-\psi(\xi)|^2 \xi
^4(1+\xi
^2)^2\cdot\frac{d\xi}{1+\xi^2}.
\end{equation}
Now, referring to equation (\ref{eqinvFourier}), $\psi_n - \psi=
\check
{\vartheta}_n$ where $\vartheta_n = [(\widehat{\psi}\phi)\ast
p_n]\cdot
\phi- \widehat{\psi}\phi^2$. Simplifying equation (\ref{eqL21}) yields
\begin{eqnarray*}
\|g_n\|_{L^2(\upsilon)}^2 &=& \int|\check{\vartheta}_n(\xi)|^2 \xi
^4(1+\xi^2)\,d\xi\le\int|\check{\vartheta}_n(\xi)|^2 \xi
^2(1+\xi
^2)^2\,d\xi\\
&=& \int|\xi(1+\xi^2)\check{\vartheta}_n(\xi)|^2\,d\xi.
\end{eqnarray*}
Since $\xi^k\check{\vartheta}_n(\xi) = (-i)^k (\vartheta
_n^{(k)})^\vee
(\xi)$ for $k\in\N$, this simplifies to
\[
\|g_n\|_{L^2(\upsilon)}^2 \le\int|(\vartheta_n')^\vee(\xi
)+(\vartheta
_n''')^{\vee}(\xi)|^2\,d\xi.
\]
That is, $\|g_n\|_{L^2(\upsilon)} \le\|(\vartheta_n')^\vee+
(\vartheta
_n''')^\vee\|_{L^2(\R)} \le\|(\vartheta_n')^\vee\|_{L^2(\R)} + \|
(\vartheta_n''')^\vee\|_{L^2(\R)} = \|\vartheta_n'\|_{L^2(\R)} + \|
\vartheta_n'''\|_{L^2(\R)}$, where we have used Parseval's identity in
the last equality. We now must compute some derivatives. Using the fact
that $(\gamma\ast p)' = \gamma'\ast p$ whenever $\gamma$ and $p$ are
functions whose convolution is well defined and $\gamma$ is $C^1$, we have
%
\begin{eqnarray}
\label{eqderivative1} \quad \vartheta_n' &=&\bigl ((\widehat{\psi}\phi)'\ast
p_n\bigr)\cdot\phi+ \bigl((\widehat{\psi}\phi)\ast p_n\bigr)\cdot\phi' -
(\widehat
{\psi}\phi^2)', \\
\label{eqderivative3} \quad \vartheta_n''' &=& \bigl((\widehat{\psi}\phi
)'''\ast
p_n\bigr)\cdot\phi+ 3\bigl((\widehat{\psi}\phi)''\ast p_n\bigr)\cdot\phi' +
3\bigl((\widehat{\psi}\phi)'\ast p_n\bigr)\cdot\phi''\nonumber
\\[-8pt]
\\[-8pt] &&{}+ \bigl((\widehat{\psi
}\phi
)\ast p_n\bigr)\cdot\phi''' - (\widehat{\psi}\phi^2)'''.
\nonumber
\end{eqnarray}
The functions $\widehat{\psi}\phi$ and $\widehat{\psi}\phi^2$ are both
in $C_c^\infty$, and so there is a constant~$A$ so that $\|(\widehat
{\psi}\phi)^{(k)}\|_{L^2(\R)} \le A$ and $\|(\widehat{\psi}\phi
^2)^{(k)}\|_{L^2(\R)}\le A$ for $0\le k\le3$. Since $\phi\in
C_c^\infty
$, there is a constant $B$ so that $\|\phi^{(k)}\|_{L^\infty(\R)}\le B$
for $0\le k\le3$. Using Young's convolution inequality $\|\gamma\ast
p\|_{L^2(\R)} \le\|\gamma\|_{L^2(\R)}\|p\|_{L^1(\R)}$,\vadjust{\goodbreak} equation~(\ref{eqderivative1}) gives us
\begin{eqnarray*}
\|\vartheta_n'\|_{L^2(\R)} &\le& B\|(\widehat{\psi}\phi)'\ast p_n\|
_{L^2(\R)} + B\|(\widehat{\psi}\phi)\ast p_n\|_{L^2(\R)} + A \\
&\le& B\|(\widehat{\psi}\phi)'\|_{L^2(\R)}\|p_n\|_{L^1(\R)} + B\|
(\widehat{\psi}\phi)\|_{L^2(\R)}\|p_n\|_{L^1(\R)} + A \\
& \le& BA + BA + A,
\end{eqnarray*}
where we use the normalization $\|p_n\|_{L^1(\R)} = 1$. A similar
calculation using equation (\ref{eqderivative3}) shows that
\[
\|\vartheta_n'''\|_{L^2(\R)} \le8BA + A.
\]
Hence, we have
%
\begin{equation} \label{eqsup}
  \sup_n \|g_n\|_{L^2(\upsilon)} \le\sup_n \bigl(\|\vartheta_n'\|
_{L^2(\R
)} + \|\vartheta_n'''\|_{L^2(\R)} \bigr) \le10BA+2A <\infty.\hspace*{-35pt}
\end{equation}
This allows us to conclude the proof as follows. For any $M>0$, we have
\[
\mathscr{I}_2(h_n-h) = \int g_n\,d\upsilon= \int g_n\1_{0\le g_n\le
M}\,
d\upsilon+ \int g_n\1_{g_n>M}\,d\upsilon.
\]
The first integrand is bounded above by $M$, and since $\upsilon$ is a
finite measure, the constant $M$ is in $L^1(\upsilon)$. Hence, since we
have already shown that $g_n\to0$ pointwise, we conclude that the
first integral converges to $0$ using the dominated convergence
theorem. For the second integral, notice that on the domain $\{g_n>M\}$
the function $g_n/M$ is $\ge1$, and so
\[
\int g_n\1_{g_n>M}\,d\upsilon\le\int g_n\cdot\frac{g_n}{M}\1
_{g_n>M}\,d\upsilon\le\frac{1}{M}\int g_n^2\,d\upsilon\le\frac
{1}{M}\sup_n \|g_n\|_{L^2(\upsilon)}^2.
\]
Since this is true for any $M$, by taking $M\to\infty$ while $n\to
\infty
$ we have $\mathscr{I}_2(h_n-h)\to0$ as desired.

 Finally, since $\mu$ is supported in $K$ and $\phi=1$ on $K$,
\begin{eqnarray*} \int h_n\,d\mu&=& \int_K \widehat{\psi}_n\,d\mu=
\int
_K (1-\phi^2)\widehat{\psi}\,d\mu+ \int_K\bigl((\widehat{\psi}\phi
)\ast
p_n\bigr)\cdot\phi\,d\mu\\
&=& \int_K (\widehat{\psi}\phi)\ast p_n\,d\mu.
\end{eqnarray*}
By construction $(\widehat{\psi}\phi)\ast p_n \to\widehat{\psi
}\phi$
(uniformly) on $K$, and also $\|(\widehat{\psi}\phi)\ast p_n\|
_{L^\infty
} \le\|\widehat{\psi}\phi\|_{L^\infty}\|p_n\|_1 = \|\widehat{\psi
}\phi
\|_{L^\infty}<\infty$. Since $\mu$ is a finite measure, the dominated
convergence theorem therefore shows that
\[
\int h_n\,d\mu= \int_K (\widehat{\psi}\phi)\ast p_n\,d\mu\to\int_K
\widehat{\psi}\phi\,d\mu= \int_K \widehat{\psi}\,d\mu= \int h\,
d\mu.
\]
This completes the proof.
\end{pf}
\end{appendix}

\section*{Acknowledgment} The authors wish to thank Bruce Driver
for many useful conversations.\vadjust{\goodbreak}


\printaddresses

\end{document}